\documentclass{siamart190516}
\usepackage[utf8]{inputenc}
\usepackage{soul}
\usepackage{todonotes}
\usepackage{enumitem}
\usepackage{amsfonts}
\usepackage{amsmath}
\usepackage{algorithm}
\usepackage{algpseudocode}% http://ctan.org/pkg/algorithmicx
\usepackage{hyperref}
\usepackage{multirow}
\usepackage{multicol}
\usepackage{subcaption}
\usepackage[capitalize]{cleveref}
\usepackage{epsfig}
\usepackage{amssymb}
\usepackage{graphicx}
\usepackage{grffile}
\usepackage{appendix}
\usepackage{placeins}
\usepackage{bm}
\usepackage{algpseudocode}
\usepackage{enumitem}

\newcommand{\ba}[1]{\begin{array}{#1}}
	
	\newcommand{\ea}{\end{array}}

\newcommand{\beq}[1]{\begin{equation}\label{#1}}
\newcommand{\eeq}{\end{equation}}

\usepackage{color}

\newcommand{\mat}[1]{\mathbf{#1}}
\newcommand{\dbar}[1]{\overline{{#1}}}
\newcommand{\ts}[1]{\bm{\mathcal{#1}}}
\newcommand{\ddd}{\cdot\cdot\cdot}
\newcommand{\midx}[1]{\bm{#1}}

\makeatletter
\renewcommand{\ALG@beginalgorithmic}{\small\setlength{\ALG@tlm}{-0.3em}}
\makeatother

\allowdisplaybreaks

\newcommand{\ignore}[1]{}
\newcommand{\ylrev}[1]{{\color{black}{#1}}}

\usepackage{tikz}
\tikzstyle{cnode} = [draw, circle,scale=.6]
\tikzstyle{level 1} = [level distance=.35\textwidth, sibling distance=1\textwidth]
\tikzstyle{level 2} = [level distance=.35\textwidth, sibling distance=.45\textwidth]
\tikzstyle{level 3} = [level distance=.3\textwidth, sibling distance=.25\textwidth]
\usetikzlibrary{decorations.pathreplacing,patterns}

\title{A Linear-complexity Tensor Butterfly Algorithm for Compressing High-dimensional Oscillatory Integral Operators}

\author{P. Michael Kielstra\thanks{\ylrev{Department of Mathematics}, University of California, Berkeley, CA, USA.\newline Email: 
		\texttt{pmkielstra@berkeley.edu}} \and Tianyi Shi\thanks{\ylrev{Applied Mathematics and Computational Research Division}, Lawrence Berkeley National Laboratory, Berkeley, CA, USA.\newline Email: 
		\texttt{\{tianyishi,liuyangzhuan\}@lbl.gov}} \and Hengrui Luo \thanks{Department of Statistics, Rice University, Houston, TX, USA. Email: {\tt hrluo@rice.edu}} \and Jianliang Qian \thanks{Department of Mathematics and Department of CMSE, Michigan State University, East Lansing, MI, USA. Email: 
		\texttt{jqian@msu.edu}} \and Yang Liu\footnotemark[2]}

\begin{document}
	\maketitle
	
	\begin{abstract}
		This paper presents a multilevel tensor compression algorithm called tensor butterfly algorithm for efficiently representing large-scale and high-dimensional oscillatory integral operators, including Green's functions for wave equations and integral transforms such as Radon transforms and Fourier transforms. The proposed algorithm leverages a tensor extension of the so-called complementary low-rank property of existing matrix butterfly algorithms. The algorithm partitions the discretized integral operator tensor into subtensors of multiple levels, and factorizes each subtensor at the middle level as a \ylrev{Tucker-type} interpolative decomposition, whose factor matrices are formed in a multilevel fashion. For a $d$-dimensional \ylrev{($d>1$)} integral operator discretized into a $2d$-mode tensor with $n^{2d}$ entries, the overall CPU time and memory requirement scale as $O(n^d)$, in stark contrast to the $O(n^d\log n)$ complexity of existing matrix algorithms such as matrix butterfly algorithms and fast Fourier transforms (FFT), where $n$ is the number of points per direction. When comparing with other tensor algorithms such as quantized tensor train (QTT), the proposed algorithm also shows superior CPU and memory performance for tensor contraction. Remarkably, the tensor butterfly algorithm can efficiently model high-frequency Green's function interactions between two unit cubes, each spanning 512 wavelengths per direction, which \ylrev{represent problems of scale over $512\times$ larger than that} existing \ylrev{butterfly} algorithms \ylrev{can handle, with the same amount of computation resources}. On the other hand, for a problem representing 64 wavelengths per direction, which is the largest size existing algorithms can handle, our tensor butterfly algorithm exhibits 200x speedups and $30\times$ memory reduction comparing with existing ones. Moreover, the tensor butterfly algorithm also permits $O(n^d)$-complexity FFTs and Radon transforms up to $d=6$ dimensions.            
	\end{abstract}
	
	

	\section{Introduction}
	
	Oscillatory integral operators (OIOs), such as Fourier transforms and Fourier integral operators \cite{hormander1994fourier,candes2007fast}, are critical computational and theoretical tools for many scientific and engineering applications, such as signal and image processing, inverse problems and imaging, computer vision, quantum mechanics, and analyzing and solving partial differential equations (PDEs). The development of accurate and efficient algorithms for computing OIOs has profound impacts on the evolution of the pertinent research areas including, perhaps mostly remarkably, the invention of the fast Fourier transform (FFT) by Cooley and Tukey in 1965 and the invention of the fast multipole method (FMM) by Greengard and Rokhlin in 1987, both of which were listed among the ten most significant algorithms discovered in the 20th century. Among existing analytical and algebraic methods for OIOs, butterfly algorithms \cite{Eric_1994_butterfly,liu2021butterfly,li2015butterfly,Yingzhou_2017_IBF,Pang2020IDBF} represent an emerging class of multilevel matrix decomposition algorithms that have been proposed for Fourier transforms and Fourier integral operators \cite {Candes_butterfly_2009,Lexing_SFT_2009,Haizhao_2018_Phase}, special function transforms \cite{Tygert_2010_spherical,bremer2020SHT,Oneil_2010_specialfunction}, fast iterative \cite{michielssen_multilevel_1996,Eric_1994_butterfly,Yang_2020_BFprecondition} and direct \cite{Han_2013_butterflyLU,Liu_2017_HODBF,Han_2017_butterflyLUPEC,Han_2018_butterflyLUDielectric,Sadeed2022VIE,liu2023detecting,sheng2023domain} solution of surface and volume integral equations for wave equations, high-frequency Green's function ansatz for inhomogeneous wave equations \cite{liu2023fast,liu2022comparative,luqiabur16}, direct solution of PDE-induced sparse systems \cite{liu2021sparse,claus2023sparse}, and machine learning for inverse problems \cite{khoo2019switchnet,li2022wide}. The (matrix) butterfly algorithms leverage the so-called complementary low-rank (CLR) property of the matrix representation of OIOs after proper row/column permutation. The CLR states that \ylrev{any submatrix with contiguous row and column index sets exhibits a low numerical rank if the number of the submatrix entries approximately equals the matrix size. These ranks are known as} the butterfly ranks, which stay constant irrespective of the matrix sizes. This permits a multilevel sparse matrix decomposition requiring $O(n\log n)$ factorization time, application time, and storage units with $n$ being the matrix size.      
	
	Despite their low asymptotic complexity, the matrix butterfly algorithms oftentimes exhibit relatively large prefactors, i.e., constant but high butterfly ranks, particularly for higher-dimensional OIOs. Examples include Green's functions for 3D high-frequency wave equations \cite{Sadeed2022VIE,liu2023fast}, 3D Radon transforms for linear inverse problems \cite{deans2007radon}, 6D Fourier–Bros–Iagolnitzer transforms for Wigner equations \cite{de2017wigner,PhysRev.40.749}, 6D Fourier transforms in diffusion magnetic resonance imaging \cite{cheng2015joint} and plasma physics \cite{delzanno2015multi}, 4D space-time transforms in quantum field theories \cite{peskin2018introduction,Maas:2007uv}, and multi-particle Green's functions in quantum chemistry \cite{fetter2012quantum}. For these high-dimensional OIOs, the computational advantage of the matrix butterfly algorithms over other existing algorithms becomes significant only for very large matrices. 
	
	More broadly speaking, for large-scale multi-dimensional scientific data and operators, tensor algorithms are typically more efficient than matrix algorithms. Popular low-rank tensor compression algorithms include CANDECOMP/PARAFAC \cite{harshman1970foundations}, Tucker \cite{de2000multilinear}, hierarchical Tucker \cite{hackbusch2009new}, tensor train (TT) \cite{oseledets2011tensor}, and tensor network \cite{cichocki2016tensor} decomposition algorithms. See references \cite{kolda2009tensor,grasedyck2013literature} for a more complete review of available tensor formats and their applications. When applied to the representation of high-dimensional integral operators, tensor algorithms often leverage additional translational- or scaling-invariance property to achieve superior compression performance, including solution of quasi-static wave equations \cite{wang2020voxcap,wang2021supervoxhenry,giannakopoulos2019memory,corona2017tensor}, elliptic PDEs \cite{biagioni2012numerical,hackbusch2008tensor}, many-body Schr\"odinger equations \cite{PhysRevB.105.165116}, and quantum Fourier transforms (QFTs) \cite{chen2024direct}. That being said, most existing tensor decomposition algorithms will break down for OIOs due to their incapability to exploit the oscillatory structure of these operators; therefore, new tensor algorithms are called for.     
	
	In this paper, we propose a linear-complexity, low-prefactor tensor decomposition algorithm for large-scale and high-dimensional OIOs. This new tensor algorithm, henceforth dubbed the tensor butterfly algorithm, leverages the intrinsic CLR property of high-dimensional OIOs more effectively than the matrix butterfly algorithm, which is enabled by additional tensor properties such as translational  invariance of free-space Green's functions and dimensional separability of Fourier transforms. The algorithm partitions the OIO tensor into subtensors of multiple levels, and factorizes each subtensor at the middle level as a \ylrev{Tucker-type} interpolative decomposition, whose factor matrices are further constructed in a nested fashion. For a $d$-dimensional OIO (assuming a constant \ylrev{$d>1$}) discretized as a $2d$-mode tensor with $n$ being the size per mode, the factorization time, application time, and storage cost scale as $O(n^d)$, and the resulting tensor factors have small multi-linear ranks. This is in stark contrast both to the $O(n^d\log n)$ scaling of existing matrix algorithms such as matrix butterfly algorithms and FFTs, and to the super-linear scaling of existing tensor algorithms. We mention that the linear complexity of the factorization time in our proposed algorithm is achieved via a simple random entry evaluation scheme, assuming that any arbitrary entry can be computed in $O(1)$ time. We remark that, for 3D high-frequency wave equations, the proposed tensor butterfly algorithm can handle \ylrev{$512\times$ larger discretized Green's function tensors} than existing \ylrev{butterfly} algorithms \ylrev{using the same amount of computation resources}; on the other hand, for the largest sized tensor that can be handled by existing \ylrev{butterfly} algorithms, our tensor butterfly algorithm is $200\times$ faster than existing ones. Moreover, we claim that the tensor butterfly algorithm instantiates the first linear-complexity implementation of high-dimensional FFTs for arbitrary input data.

	\subsection{Related Work}\label{sec:relatedwork}
	
	\textit{Multi-dimensional butterfly algorithms} represent a version of matrix butterfly algorithms designed for high-dimensional OIOs \cite{li2018multidimensional,chen2020multidimensional}. Instead of the traditional binary tree partitioning of the matrix rows/columns \cite{Eric_1994_butterfly}, these algorithms can be viewed as a modern version of \cite{michielssen_multilevel_1996} that permits quadtree and octree partitioning of the matrix rows/columns, which have been demonstrated on 2D and 3D OIOs. For a general $d$-dimensional OIO, the $d$-dimensional tree partitioning leads to a butterfly factorization with a $d$-fold reduction in the number of levels compared to the binary tree partitioning. \ylrev{That said, the binary tree based butterfly algorithms are easier to implement and exhibit very competitive overall costs comparing with the multi-dimensional butterfly algorithms.} We note that both the multi-dimensional and binary tree-based butterfly algorithms are still matrix-based algorithms that scale as $O(n^d\log n)$, as opposed to the proposed tensor algorithm that scales as $O(n^d)$.
	
	\textit{Quantized tensor train (QTT) algorithms}, or simply TT algorithms, are tensor algorithms well-suited for very high-dimensional integral operators. They have been proposed to compress volume integral operators \cite{corona2017tensor} arising from quasi-static wave equations and static PDEs with $O(\log n)$ memory and CPU complexities. However, for high-frequency wave equations, the QTT rank scales proportionally to the wave number \cite{corona2017tensor} leading to deteriorated CPU and memory complexities (see our numerical results in Section \cref{sec:example}). Moreover, QTT has been proposed for computing FFT and QFT with $O(\log n)$ memory and CPU complexities \cite{chen2024direct}. However, after obtaining the QTT-compressed formats of both the volume-integral operator and the Fourier transform, the CPU complexity for contracting such a QTT compressed operator with arbitrary (i.e., non QTT-compressed) input data scales super-linearly. In contrast, our algorithm yields a linear CPU and memory complexity for the contraction operation.

	\subsection{Contents}
	In what follows, we first review the matrix low-rank decomposition and butterfly decomposition algorithms in \cref{sec:matrix_algorithm}. In \cref{sec:tuckerID}, we introduce the \ylrev{Tucker-type} interpolative decomposition algorithm as the building block for the proposed tensor butterfly algorithm detailed in \cref{sec:butterflytensor}. The multi-linear butterfly ranks for a few special cases are analyzed in \cref{sec:rank} and the complete complexity analysis is given in \cref{sec:cc}. \Cref{sec:example} shows a variety of numerical examples, including Green's functions for wave equations, Radon transforms, and uniform and non-uniform discrete Fourier transforms, to demonstrate the performance of matrix butterfly, tensor butterfly, Tucker and QTT algorithms.    
	
	\subsection{Notations}
	Given a scalar-valued function $f(x)$, its integral transform is defined as \begin{equation}g(x)=\int_y K(x,y)f(y) dy\end{equation} with an integral kernel $K(x,y)$. The indexing of a matrix $\mat{K}$ is denoted by $\mat{K}(i,j)$ or $\mat{K}(t,s)$, where $i,j$ are indices and $t,s$ are index sets. We use $\mat{K}^T$ to denote the transpose of matrix $\mat{K}$. For a sequence of matrices $\mat{K}_1,\ldots,\mat{K}_n$, the matrix product is \begin{equation}\prod_{i=1}^{n}\mat{K}_i=\mat{K}_1\mat{K}_2\ldots\mat{K}_n,\end{equation} the vertical stacking (assuming the same column dimension) is \begin{equation}[\mat{K}_i]_{i}=[\mat{K}_1; \mat{K}_2; \ldots;\mat{K}_n],\end{equation} and \begin{equation}\text{diag}_i(\mat{K}_i)=\text{diag}(\mat{K}_1,\mat{K}_2,\ldots,\mat{K}_n)\end{equation} is a block diagonal matrix with $\mat{K}_i$ being the diagonal blocks. Given an $L$-level binary-tree partitioning $\mathcal{T}_{t}$ of an index set $t=\{1,2,\ddd,n\}$, any node $\tau$ at each level is a subset of $t$. The parent and children of $\tau$ are denoted by $p_\tau$ and $\tau^{c}$ ($c=1,2$), respectively, and $\tau=\tau^1\cup\tau^2$.    
	
	A multi-index $\midx{i}=(i_1,\ddd,i_d)$ is a tuple of indices, and similarly a multi-set $\midx{\tau}=(\tau_1,\tau_2,\ddd,\tau_d)$ is a tuple of index sets. We define \begin{equation}\midx{\tau}_{k\leftarrow t} = (\tau_1,\tau_2,\ddd,\tau_{k-1},t,\tau_{k+1},\tau_{k+2},\ddd,\tau_{d}).\end{equation} Given a tuple of nodes (i.e. a multi-set) $\midx{\tau}=(\tau_1,\tau_2,\ddd,\tau_d)$ and a multi-index $\midx{c}=(c_1,c_2,\ddd,c_d)$ with $c_i\in\{1,2\}$, the children of $\midx{\tau}$ are denoted $\midx{\tau}^{\midx{c}}=({\tau_{1}}^{c_1},{\tau_{2}}^{c_2},\ddd,{\tau_{d}}^{c_d})$ and the parents of $\tau_i$, $i=1,2,\ddd,d$ can be simply written as $\midx{p}_{\midx{\tau}}=(p_{\tau_1},p_{\tau_2},\ddd,p_{\tau_d})$. Similar to the above-described notations, we can replace the index $i$ in $[\mat{K}_i]_{i}$ and $\text{diag}_i(\mat{K}_i)$ with an index set $\tau$, a multi-index $\midx{c}$, or a multi-set $\midx{\tau}$ assuming certain predefined index ordering.       
	
	Given complex-valued (or real-valued) functions $f(x)$ of $d$ variables and integral operators $K(x,y)$, the tensor representations of their discretizations are respectively denoted by $\ts{F}\in \mathbb{C}^{n_1\times n_2\times\ddd\times n_d}$ and $\ts{K}\in \mathbb{C}^{m_1\times m_2\times\ddd\times m_d\times n_1\times n_2\times\ddd\times n_d}$, where $n_1,\cdots,n_d$ and $m_1,\cdots,m_d$ are sizes of discretizations for the corresponding variables. In this paper, we use \textit{matricization} to denote the reshaping of $\ts{K}$ into a $(\Pi_{k} m_k)\times (\Pi_{k} n_k)$ matrix, and the reshaping of $\ts{F}$ into a $(\Pi_{k} n_k)\times 1$ matrix. The entries of $\ts{F}$ and $\ts{K}$ are denoted by $\ts{F}(\midx{i})$ (or equivalently $\ts{F}(i_1,i_2,\ddd,i_d)$) and $\ts{K}(\midx{i},\midx{j})$, respectively. Similarly the subtensors are denoted by $\ts{F}(\midx{\tau})$ (or equivalently $\ts{F}(\tau_1,\tau_2,\ddd,\tau_d)$) and $\ts{K}(\midx{\tau},\midx{\nu})$. %The concatenation of two $d$-mode tensors $\ts{F}$ and $\ts{G}$ along dimension $k$ is denoted $\begin{bmatrix}\ts{F}, \ts{G}\end{bmatrix}_{\mathrm{dim}=k}$. 
	
	Given a $d$-mode tensor $\ts{F}\in \mathbb{C}^{n_1\times n_2\times\ddd\times n_d}$, the mode-$j$ unfolding is denoted by $\mat{F}^{(j)}\in \mathbb{C}^{(\Pi_{k\neq j} n_k)\times n_j}$, the mode-$j$ tensor-matrix product of $\ts{F}$ with a matrix $\mat{X}\in\mathbb{C}^{m\times n_j}$ is denoted by $\ts{Y}=\ts{F}\times_j \mat{X}$, or equivalently $\mat{Y}^{(j)}=\mat{F}^{(j)}\mat{X}^T$.

	\section{Review of Matrix Algorithms}\label{sec:matrix_algorithm} We consider a $d$-dimensional OIO kernel $K(x,y)$ with $x,y\in\mathbb{R}^d$ discretized on point pairs $x^{i}$ and $y^{j}$, $i=1,2,...,(m_1 m_2 \ddd  m_d)$, $j=1,2,...,(n_1 n_2 \ddd  n_d)$, where $i$ (and similarly $j$) is the flattening of the corresponding multi-index $\midx{i}$. Such a discretization can be represented as a matrix $\mat{K}\in \mathbb{C}^{(m_1 m_2 \ddd  m_d)\times (n_1  n_2 \ddd  n_d)}$. When it is clear in the context, we assume that $m_k=n_k=n$ for $k=1,\ldots,d$. Throughout this paper, we assume that $\mat{K}$ (and its tensor representation) is never fully formed, but instead a function is provided to evaluate any matrix (or tensor) entry in $O(1)$ time. Next we review matrix compression algorithms for $\mat{K}$ including low-rank and butterfly algorithms.  
	
	\subsection{Interpolative Decomposition}\label{sec:matID}	
	The interpolative decomposition (ID) algorithm \cite{halko_finding_2011,libwoomarroktyg07} is a matrix compression technique that constructs a low-rank decomposition whose factors contain original entries of the matrix. More specifically, consider the matrix $\mat{K}(\tau,\nu)\in \mathbb{C}^{m\times n}$, $\tau=\{1,2,\ldots,m\}$, $\nu=\{1,2,\ldots,n\}$, the column ID of $\mat{K}$ (the index sets $\tau$ and $\nu$ are omitted for clarity in context) is
	\begin{equation}\label{eqn:col-ID}
	\mat{K} \approx \mat{K}(:, \overline{\nu}) {\mat{V}},
	\end{equation}    
	where the skeleton matrix $\mat{K}(:, \overline{\nu})$ contains $r$ skeleton columns indexed by $\overline{\nu} \subseteq \nu$ and the interpolation matrix $\mat{V}$ has bounded entries. Here the numerical rank $r$ is chosen such that \begin{equation}\| \mat{K} - \mat{K}(:,\overline{\nu}) {\mat{V}}\|_F^2 \leqslant O(\epsilon^2)\|\mat{K}\|_F^2\end{equation} for a prescribed relative tolerance $\epsilon$. In practice, the column ID can be computed via rank-revealing QR decomposition with a relative tolerance $\epsilon$ \cite{libwoomarroktyg07}. Similarly, the row ID of the matrix $\mat{K}$ reads 
	\begin{equation}\label{eqn:row-ID}
	\mat{K} \approx {\mat{U} \mat{K}(\overline{\tau},:)},
	\end{equation}    
	where the skeleton matrix $\mat{K}(\overline{\tau},:)$ contains $r$ skeleton rows indexed by $\overline{\tau} \subseteq \tau$ and the interpolation matrix $\mat{U}$ has bounded entries. The row ID can be simply computed by the column ID of $\mat{K}^T$. Combining the column and row ID in \cref{eqn:col-ID} and \cref{eqn:row-ID} gives 
	\begin{equation}\label{eqn:hybrid-ID}
	\mat{K} \approx {\mat{U}} \mat{K}(\overline{\tau}, \overline{\nu}) {\mat{V} }.
	\end{equation}    
	It is straightforward to note that the memory and CPU complexities of ID scale as $O(nr)$ and $O(n^2r)$, respectively. The CPU complexity can be reduced to $O(nr^2)$ when properly selected proxy rows in \cref{eqn:col-ID} and columns in \cref{eqn:row-ID} are used in the rank-revealing QR. Common strategies of choosing proxy rows/columns (henceforth called \textit{proxy index} strategies) for integral operators include evenly spaced or uniform random samples, and more generally the use of Chebyshev nodes and proxy surfaces (where new rows $K(\overline{x},y^j)$ other than original rows of $\mat{K}$ are used with $\overline{x}$ denoting the proxies). However, for large OIOs \ylrev{(e.g., Green's functions of high-frequency wave equations  discretized with a small number of points per wavelength)}, the rank $r$ depends on the size $n$ of the matrix; consequently, ID is not an efficient compression algorithm. Next, we review the matrix butterfly algorithm capable of achieving quasi-linear memory and CPU complexities for OIOs.  
	
	\subsection{Matrix Butterfly Algorithm}\label{sec:matBF}		
	\ylrev{For reasons discussed in \cref{sec:relatedwork}, we only consider the binary tree based matrix butterfly algorithm as the reference algorithm for the proposed tensor butterfly algorithm throughout this paper.} Let $t^0=\{1,2,\ddd,m\}$ and $s^0=\{1,2,\ddd,n\}$. \ylrev{Without loss of generality, we assume that $m=n$.} The $L$-level butterfly representation of the discretized OIO $\mat{K}(t^0,s^0)$ is based on two binary trees, $\mathcal{T}_{{t^0}}$ and $\mathcal{T}_{s^0}$, and the CLR property of the OIO takes the following form: at any level $0\leq l \leq L$, for any node $\tau$ at level $l$ of $\mathcal{T}_{{t^0}}$ and any node $\nu$ at level $L-l$ of $\mathcal{T}_{s^0}$, the subblock $\mat{K}(\tau,\nu)$ is numerically low-rank with rank $r_{\tau,\nu}$ bounded by a small number $r$ called the butterfly rank \cite{liu2021butterfly,Yingzhou_2017_IBF,li2015butterfly,Pang2020IDBF}.
	
	For any subblock $\mat{K}(\tau,\nu)$, \ylrev{the} ID in \cref{eqn:hybrid-ID} permits 
	\begin{equation}\label{eqn:IDtaunu}
	\mat{K}(\tau,\nu) \approx \mat{U}_{\tau,\nu} \mat{K}(\overline{\tau}, \overline{\nu}) \mat{V}_{\tau,\nu},
	\end{equation}
	where the skeleton rows and columns are indexed by $\overline{\tau}$ and $\overline{\nu}$, respectively. It is worth noting that given a node $\nu$, the selection of skeleton columns $\overline{\nu}$ depends on the node $\tau$. However, the notation $\bar{\cdot}$ does not reflect the dependency when it is clear in the context. \ylrev{By CLR, there are at most $r$ skeleton rows and columns.} 
	
	Without loss of generality, we assume that $L$ is an even number so that $L^c=L/2$ denotes the middle level. At levels $l=0,\ldots,L^c$, the interpolation matrices $\mat{V}_{\tau,\nu}$ are computed as follows:   
	
	At level $l=0$, ${\mat{V}}_{\tau,\nu}$ are explicitly formed. While at level $0<l\leq L^c$, they are represented in a nested fashion. To see this, consider a node pair $(\tau,\nu)$ at level $l>0$ and let $\nu^1,\nu^2$ and $p_\tau$ be the children and parent of $\nu$ and $\tau$, respectively. Let $s$ be the ancestor of $\nu$ at level $L^c$ of $\mathcal{T}_{s^0}$ and let $\mathcal{T}_{s}$ denote the subtree rooted at $s$.

	\ylrev{By \cref{eqn:hybrid-ID}}, we have  
	\begin{align}\label{eqn:transfer_computation}
	\mat{K}(\tau,\nu) &=\begin{bmatrix}
	\mat{K}(\tau,{\nu^1}) & \mat{K}(\tau,{\nu^2})\end{bmatrix}\nonumber\\
	&\approx\begin{bmatrix}
	\mat{K}(\tau,\overline{\nu^1}) & \mat{K}(\tau,\overline{\nu^2})\end{bmatrix}\begin{bmatrix}
	\mat{V}^s_{p_\tau, \nu^1} & \\
	& \mat{V}^s_{p_\tau, \nu^2}
	\end{bmatrix}\\
	&\approx\mat{K}(\tau,\overline{\nu})\mat{W}^s_{\tau,\nu}\begin{bmatrix}
	\mat{V}^s_{p_\tau, \nu^1} & \\
	& \mat{V}^s_{p_\tau, \nu^2}
	\end{bmatrix}.
	\end{align}
	Here $\mat{W}^s_{\tau,\nu}$ and $\overline{\nu}$ are the interpolation matrix and skeleton columns from the ID of $\mat{K}(\tau,\overline{\nu^1}\cup\overline{\nu^2})$, respectively. $\mat{W}_{\tau,\nu}$ is henceforth referred to as the transfer matrix for $\nu$ in the rest of this paper. \ylrev{By CLR, $\mat{W}_{\tau,\nu}$ is of sizes at most $r\times2r$.} Note that we have added an additional superscript $s$ to $\mat{V}_{p_\tau, \nu^c}$ and $\mat{W}_{\tau, \nu}$, for notation convenience in the later context. From \cref{eqn:transfer_computation}, it is clear that the interpolation matrix $\mat{V}^s_{\tau,\nu}$ can be expressed in terms of its parent $p_\tau$'s and children $\nu^1,\nu^2$'s  interpolation matrices as 
	\begin{equation}\label{eqn:nested_basis}
	\mat{V}^s_{\tau,\nu} =
	\mat{W}^s_{\tau,\nu}\begin{bmatrix}
	\mat{V}^s_{p_\tau, \nu^1} & \\
	& \mat{V}^s_{p_\tau, \nu^2}
	\end{bmatrix}.
	\end{equation}
	Note that the interpolation matrices $\mat{V}^s_{\tau,\nu}$ at level $l=0$ and transfer matrices $\mat{W}^s_{\tau,\nu}$ at level $0<l\leq L^c$ do not require the column ID on the full subblocks $\mat{K}(\tau,\nu)$ and $\mat{K}(\tau,\overline{\nu^1}\cup\overline{\nu^2})$, which would lead to at least an  $O(mn)$ computational complexity. 
	
	In practice, one can select $O(r_{\tau,\nu})$ proxy rows $\hat{\tau}\subset\tau$ to compute $\mat{V}^s_{\tau,\nu}$ and $\mat{W}^s_{\tau,\nu}$ via ID as:
	\begin{align}\label{eqn:proxy}
	\mat{K}(\hat{\tau},\nu) &\approx \mat{K}(\hat{\tau}, \overline{\nu}) \mat{V}^s_{\tau,\nu},~~l=0,\\
	\mat{K}(\hat{\tau},\overline{\nu^1}\cup\overline{\nu^2})&\approx\mat{K}(\hat{\tau},\overline{\nu})\mat{W}^s_{\tau,\nu},~~0<l\leq L^c. 
	\end{align} 
	The viable choices for proxy rows have been discussed in several existing papers \cite{liu2023fast,Pang2020IDBF,Sadeed2022VIE,Candes_butterfly_2009}. 
	
	At levels $l=L^c,\ldots,L$, the interpolation matrices $\mat{U}_{\tau,\nu}$ are computed by performing similar operations on $\mat{K}^T$. We only provide their expressions here and omit the redundant explanation. Let $t$ be the ancestor of $\nu$ at level $L^c$ of $\mathcal{T}_{t^0}$ and let $\mathcal{T}_t$ be the subtree rooted at $t$. At level $l=L$, $\mat{U}^t_{\tau,\nu}$ are explicitly formed. At level $L^c\leq l<L$, only the transfer matrices $\mat{P}^t_{\tau,\nu}$ are computed from the column ID of $\mat{K}^T(\nu,\overline{\tau^1}\cup\overline{\tau^2})$ satisfying
	\begin{equation}\label{eqn:nested_basis_P}
	\mat{U}^t_{\tau,\nu} =
	\begin{bmatrix}
	\mat{U}^t_{\tau^1, p_\nu} & \\
	& \mat{U}^t_{\tau^2, p_\nu}
	\end{bmatrix}\mat{P}^t_{\tau,\nu}.
	\end{equation}  
	
	Combining \cref{eqn:IDtaunu}, \cref{eqn:nested_basis} and \cref{eqn:nested_basis_P}, the matrix butterfly decomposition can be expressed for each node pair $(t,s)$ at level $L^c$ of $\mathcal{T}_{t^0}$ and $\mathcal{T}_{s^0}$ as 
	\begin{align}
	\mat{K}(t,s) \approx \dbar{\mat{U}}^t\bigg(\prod_{l=1}^{L^c}\dbar{\mat{P}}^{t,s}_{l}\bigg)\mat{K}(\overline{t},\overline{s})\bigg(\prod_{l=L^c}^{1}\dbar{\mat{W}}^{t,s}_{l}\bigg)\dbar{\mat{V}}^s.
	\label{eqn:butterfly_mat}
	\end{align}
	Here, $\overline{t}$ and $\overline{s}$ represent the skeleton rows and columns of the ID of $\mat{K}(t,s)$. %With a slight abuse of notations, let $\nu_1,\nu_2,\ldots,\nu_{2^{l}}$ denote the nodes at level $l$ of $\mathcal{T}_{s^0}$ or $\mathcal{T}_{s}$ (Recall that $s^0$ represents the root of $\mathcal{T}_{s^0}$ and $s$ represents a node at level $L^c$ of $\mathcal{T}_{s^0}$). Similarly, let $\tau_1,\tau_2,\ldots,\tau_{2^{l}}$ denote the nodes at level $l$ of $\mathcal{T}_{t^0}$ or $\mathcal{T}_{t}$. 
	The interpolation factors $\dbar{\mat{U}}^t$ and $\dbar{\mat{V}}^s$ in \cref{eqn:butterfly_mat} are 
	\begin{align}
	\dbar{\mat{U}}^t&=\mathrm{diag}_{\tau}(\mat{U}^t_{\tau,s^0})
	\label{eqn:butterfly_factors_U},~~\tau \mathrm{~at~level~} L^c \mathrm{~of~} \mathcal{T}_t,\\
	\dbar{\mat{V}}^s&=\mathrm{diag}_{\nu}(\mat{V}^s_{t^0,\nu})
	\label{eqn:butterfly_factors_V},~~\nu \mathrm{~at~level~} L^c \mathrm{~of~} \mathcal{T}_s, 
	\end{align}
	and the transfer factors $\dbar{\mat{P}}_l^{t,s}$ and $\dbar{\mat{W}}_l^{t,s}$ for $l=1$, $\ldots$, $L^c$ consist of transfer matrices $\mat{W}^s_{{\tau},\nu}$ and $\mat{P}^s_{{\tau},\nu}$:
%\begin{align}
%\dbar{\mat{W}}^{t,s}_{l}&=\mathrm{diag}_{\tau}\Big(\begin{bmatrix}
%\text{diag}_{\nu}(\mat{W}^s_{{\tau}^c,\nu}) 
%\end{bmatrix}_{c}\Big), \begin{aligned}~~\tau &\mathrm{~at~level~} l-1 \mathrm{~of~} \mathcal{T}_{t^0}, \mathrm{~and~} t \subseteq \tau, \\ \nu &\mathrm{~at~level~} L^c-l \mathrm{~of~} \mathcal{T}_{s};\end{aligned}\label{eqn:butterfly_factors_W}\\
%(\dbar{\mat{P}}^{t,s}_{l})^T&=\mathrm{diag}_{\nu}\Big(\begin{bmatrix}
% \text{diag}_{\tau}\big((\mat{P}^t_{\tau,{\nu}^c})^T\big) 
% \end{bmatrix}_{c}\Big), \begin{aligned}~~\tau &\mathrm{~at~level~} L^c-l \mathrm{~of~} \mathcal{T}_{t},\\ \nu &\mathrm{~at~level~} l-1 \mathrm{~of~} \mathcal{T}_{s^0},\mathrm{~and~} s \subseteq \nu.\end{aligned}\label{eqn:butterfly_factors_P}
%\end{align}
\begin{align}
	\dbar{\mat{W}}^{t,s}_{l}&=\text{diag}_{\nu}(\mat{W}^s_{{\tau},\nu}), \begin{aligned}~~\tau &\mathrm{~at~level~} l \mathrm{~of~} \mathcal{T}_{t^0}, \mathrm{~and~} t \subseteq \tau, \\ \nu &\mathrm{~at~level~} L^c-l \mathrm{~of~} \mathcal{T}_{s};\end{aligned}\label{eqn:butterfly_factors_W}\\
	(\dbar{\mat{P}}^{t,s}_{l})^T&=\text{diag}_{\tau}\big((\mat{P}^t_{\tau,{\nu}})^T\big), \begin{aligned}~~\tau &\mathrm{~at~level~} L^c-l \mathrm{~of~} \mathcal{T}_{t},\\ \nu &\mathrm{~at~level~} l \mathrm{~of~} \mathcal{T}_{s^0},\mathrm{~and~} s \subseteq \nu.\end{aligned}\label{eqn:butterfly_factors_P}
	\end{align}\ylrev{For the ease of comparison with the tensor butterfly algorithm in \cref{sec:butterflytensor}, we list some notations of the matrix butterfly algorithm in \cref{tab:notation}.}

	\begin{table}[t]
		\centering
		\begin{tabular}{|c c c|}
			\hline
			Meaning & Matrix butterfly & Tensor butterfly\\ 
			\hline        
			Butterfly rank & $r_m$ & $r_t$\\
			Set/multi-set & $\tau,\nu$ & $\midx{\tau}, \midx{\nu}$\\
			$k^{th}$ set of multi-set & - & $\tau_k,\nu_k$\\
			Parent set/multi-set &$p_\tau$ &$\midx{p}_{\midx{\tau}}$ \\
			Children set/multi-set &$\tau^{c}$ & $\midx{\tau}^{\midx{c}}$\\         
			Root-level set/multi-set & $t^0,s^0$ & $\midx{t}^0, \midx{s}^0$ \\         
			Mid-level set/multi-set & $t,s$ & $\midx{t}, \midx{s}$ \\
			Binary tree & $\mathcal{T}_{{t^0}}$, $\mathcal{T}_{s^0}$& $\mathcal{T}_{{t^0_k}}$, $\mathcal{T}_{s^0_k}$\\
			Cardinality of leaf nodes& $C_b^d$ & $C_b$\\   
			Cardinality of root nodes& $n^d$ & $n$\\            
			Mid-level submatrix/subtensor & $\mat{K}(\overline{t},\overline{s})$ & $\ts{K}(\overline{\midx{t}},\overline{\midx{s}})$ \\
			Interpolation matrix & $\mat{V}^s_{\tau,\nu}$,$\mat{U}^t_{\tau,\nu}$ & $\mat{V}^{\midx{s},k}_{\midx{\tau},\nu}$, $\mat{U}^{\midx{t},k}_{\tau,\midx{\nu}}$ \\
			Transfer matrix & $\mat{W}^s_{\tau,\nu}$,$\mat{P}^t_{\tau,\nu}$ & $\mat{W}^{\midx{s},k}_{\midx{\tau},\nu}$, $\mat{P}^{\midx{t},k}_{\tau,\midx{\nu}}$ \\
			Interpolation factor &$\dbar{\mat{U}}^{t}$, $\dbar{\mat{V}}^{s}$ & $\dbar{\mat{U}}^{\midx{t},k}$, $\dbar{\mat{V}}^{\midx{s},k}$\\
			Transfer factor &$\dbar{\mat{P}}_l^{t,s}$, $\dbar{\mat{W}}_l^{t,s}$ & $\dbar{\mat{P}}_l^{\midx{t},\midx{s},k}$, $\dbar{\mat{W}}_l^{\midx{t},\midx{s},k}$\\
			\hline
		\end{tabular}
		\caption{\ylrev{Notation comparison of the matrix butterfly algorithm in \cref{sec:matBF} and the tensor butterfly algorithm in \cref{sec:butterflytensor}. Note that the subscript $k$ in $\tau_k,\nu_k$,  in the tensor notations of the interpolation/transfer matrix and interpolation/transfer factor for dimension $k$, is dropped for simplicity throughout this paper.}}\label{tab:notation}
	\end{table}

	The CPU and memory requirement for computing the matrix butterfly decomposition can be briefly analyzed as follows. Note that we only need to analyze the costs for $\mat{V}^s_{\tau,\nu}$, $\mat{W}^s_{\tau,\nu}$ and $\mat{K}(\overline{t},\overline{s})$ as those for $\mat{U}^t_{\tau,\nu}$ and $\mat{P}^t_{\tau,\nu}$ are similar. By the CLR assumption, we assume that $r_{\tau,\nu}\leq r, \forall \tau,\nu$ for some constant $r$. Thanks to the use of the proxy rows and columns, the computation of one individual $\mat{V}^s_{\tau,\nu}$ and $\mat{W}^s_{\tau,\nu}$ by ID only operates on $O(r)\times O(r)$ matrices, hence its memory and CPU requirements are $O(r^2)$ and $O(r^3)$, respectively. In total, there are $O(2^{L^c})$ middle-level nodes $s$ each having $O(2^{L^c})$ numbers of $\mat{V}^s_{\tau,\nu}$ and $O(L^c2^{L^c})$ numbers of $\mat{W}^s_{\tau,\nu}$. Similarly, each $\mat{K}(\overline{t},\overline{s})$ requires $O(r^2)$ CPU and memory costs, and there are in total $O(2^L)$ middle-level node pairs $(t,s)$. These numbers sum up to the overall $O(nr^2\log n)$ memory and $O(nr^3\log n)$ CPU complexities for matrix butterfly algorithms.  
	
	For $d$-dimensional discretized OIOs $\mat{K}\in \mathbb{C}^{(m_1 m_2 \ddd  m_d)\times (n_1  n_2 \ddd  n_d)}$ with $m_k=n_k=n$, we can assume that $n=C_b2^L$ with some constant $C_b$. For the above-described binary-tree-based butterfly algorithm, the leaf nodes of the trees are of size $C_b^d$ and this leads to a $dL$-level butterfly factorization. The memory and CPU complexities for this algorithm become $O(dn^dr^2\log n)$ and $O(dn^dr^3\log n)$, respectively. On the other hand, the multi-dimensional tree-based butterfly algorithm \cite{li2018multidimensional,chen2020multidimensional} leads to a $L$-level factorization with $O(2^dn^dr^2\log n)$ memory and $O(2^dn^dr^3\log n)$ CPU complexities. \ylrev{In this paper, we only use the binary-tree-based algorithm as the baseline matrix butterfly algorithm.} Despite their quasi-linear complexity for high-dimensional OIOs, the butterfly rank $r$ is constant but high, leading to very large prefactors of these binary and multi-dimensional tree-based algorithms. In the following, we turn to tensor decomposition algorithms to reduce both the prefactor and asymptotic scaling of matrix butterfly algorithms. \ylrev{The proposed tensor decomposition explores additional tensor compressibility of high-dimensional OIOs such as translational invariance of free-space Green's functions and dimensional separability of Fourier transforms. As will be clear in the next section, the prefactor (dependent on the butterfly rank) can be reduced by leveraging Tucker decomposition for tensorization of the middle-level submatrices $\mat{K}(\overline{t},\overline{s})$ of \cref{eqn:butterfly_mat}. The Tucker decomposition is further factored out along each dimension in a nested fashion by simultaneously moving along the binary tree of that dimension and $d$ binary trees of other dimensions. As a result, the number of transfer matrices becomes dominant only towards the middle level $L^c$, leading to a factor of $\log n$ reduction in the asymptotic complexity.}    
	
	\section{Proposed Tensor Algorithms}\label{sec:tensor_algorithm}
	In this section, we assume that the $d$-dimensional discretized OIO in \cref{sec:matrix_algorithm} is directly represented as a $2d$-mode tensor $\ts{K}\in \mathbb{C}^{m_1\times m_2\times\ddd\times m_d\times n_1\times n_2\times\ddd\times n_d}$. We first extend the matrix ID algorithm in \cref{sec:matID} to its tensor variant, which serves as the building block for the proposed tensor butterfly algorithm.  
	
	\begin{figure}[btp!]
		\centering
		\vspace{-7.5pt}
		\includegraphics[width=0.9\linewidth]{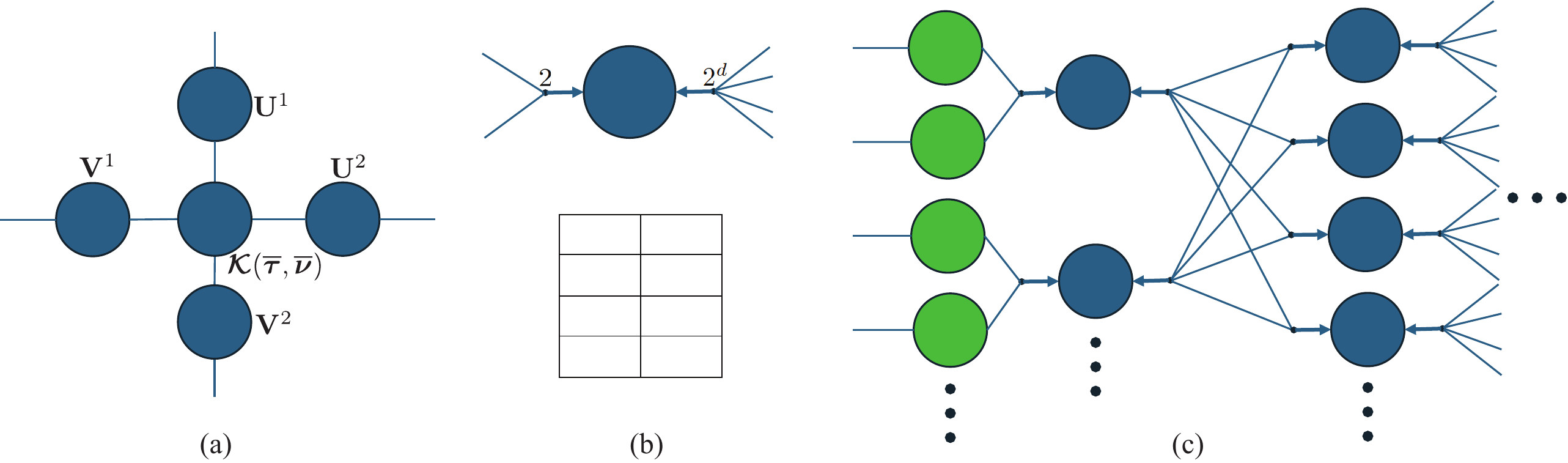}
		\vspace{-5pt}
		\caption{Tensor diagrams for (a) the Tucker-ID decomposition of a 4-mode tensor, and (b) the matrix partitioner corresponding to a $2^d\times 2$ partitioning with $d=2$ used in the tensor butterfly decomposition of a 2$d$-mode tensor, such as $\begin{bmatrix}\mat{W}^{\midx{s},k}_{\midx{\tau}^{\midx{c}},\nu}\end{bmatrix}_{\midx{c}}$ in \cref{eqn:tensor_butterfly_factors_W} for fixed $\midx{s}, \midx{\tau}, k$ and $\nu$, or $\begin{bmatrix}\mat{P}^{\midx{t},k}_{\tau,\midx{\nu}^{\midx{c}}} 
			\end{bmatrix}_{\midx{c}}$ in \cref{eqn:tensor_butterfly_factors_P} for fixed $\midx{t}, \midx{\nu}, k$ and $\tau$. \ylrev{Here, each of the row and column dimensions is connected to a partitioning node. Each partitioning node has a parent edge with an arrow pointing to the dimension to be partitioned, and several children edges connected to the parent edge. The weight of the parent edge (i.e., the number of columns or rows of the matrix) equals the sum of the weights of the children edges.} (c) The tensor diagram involving blocks $\mat{V}^{\midx{s},k}_{\midx{t}^0,\nu}$ (in green) and blocks $\begin{bmatrix}\mat{W}^{\midx{s},k}_{\midx{\tau}^{\midx{c}},\nu}\end{bmatrix}_{\midx{c}}$ (in blue) for fixed $\midx{s}$ and $k$ for the tensor butterfly decomposition of a 2$d$-mode tensor.}
		\label{fig:tensor_diagram_part}
		%	\vspace{-20pt}
	\end{figure}	
	
	\subsection{\ylrev{Tucker-type} Interpolative Decomposition}\label{sec:tuckerID}	
	Given the 2$d$-mode tensor $\ts{K}(\midx{\tau},\midx{\nu})$ with $\tau_{k}=\{1,2,\ldots,m_k\}$ and $\nu_{k}=\{1,2,\ldots,n_k\}$ for $k=1,\ldots,d$, the proposed \ylrev{Tucker-type} decomposition compresses each dimension independently via the column ID of the unfolding of $\ts{K}$ along the $k$-th  dimension, 
	\begin{equation}\label{eqn:unfoldingID}
	\mat{K}^{(k)} \ylrev{\approx} \mat{K}^{(k)}(:,\overline{\tau_k})\mat{U}^k, ~~~ \mat{K}^{(d+k)} \ylrev{\approx} \mat{K}^{(d+k)}(:,\overline{\nu_k})\mat{V}^k,~~ k=1,\ldots,d, 
	\end{equation}  
	where $\mathbf{K}^{(k)}\in\mathbb{C}^{(\prod_{j\neq k}n_j)\times n_k}$ is the mode-$k$ unfolding, or equivalently 
	\begin{equation}\label{eqn:tensorID}
	\ts{K}\ylrev{\approx} \ts{K}(\midx{\tau}_{k\leftarrow \overline{\tau_k}},\midx{\nu})\times_k\mat{U}^k, ~~~ \ts{K}\ylrev{\approx} \ts{K}(\midx{\tau},\midx{\nu}_{k\leftarrow \overline{\nu_k}})\times_{d+k}\mat{V}^k,~~ k=1,\ldots,d. 
	\end{equation} 
	Here, $\overline{\tau_k}$ and $\overline{\nu_k}$ denote the skeleton indices along modes $k$ and $d+k$ of $\ts{K}$, respectively, while $\midx{\tau}_{k\leftarrow \overline{\tau_k}}$ and $\midx{\nu}_{k\leftarrow \overline{\nu_k}}$ denote multi-sets that replace $\tau_k$ and $\nu_k$, respectively, with $\overline{\tau_k}$ and $\overline{\nu_k}$. Combining \cref{eqn:tensorID} for all dimensions yields the following proposed Tucker-type decomposition, 
	\begin{equation}\label{eqn:tuckerID}
	\ts{K} \ylrev{\approx} \ts{K}(\overline{\midx{\tau}},\overline{\midx{\nu}})\bigg(\prod_{k=1}^{d}\times_k\mat{U}^k\bigg)\bigg(\prod_{k=1}^{d}\times_{d+k}\mat{V}^k\bigg), 
	\end{equation}  
	where, $\overline{\midx{\tau}}=(\overline{\tau_1},\overline{\tau_2},\ldots,\overline{\tau_d})$, $\overline{\midx{\nu}}=(\overline{\nu_1},\overline{\nu_2},\ldots,\overline{\nu_d})$, the core tensor $\ts{K}(\overline{\midx{\tau}},\overline{\midx{\nu}})$ is a subtensor of $\ts{K}$, and $\mat{U}^k$ and $\mat{V}^k$ are the factor matrices for modes $k$ and $d+k$, respectively. 
	
	See \cref{fig:tensor_diagram_part}(a) for the tensor diagram of \cref{eqn:tuckerID} for a 4-mode tensor, \ylrev{which has the same diagram as other existing Tucker decompositions such as high-order singular value decompositions (HOSVD) \cite{de2000multilinear}. However, unlike HOSVD} that leads to orthonormal factor matrices, the proposed decomposition leads to factor matrices with bounded entries and the core tensor with the original tensor entries. 
	Therefore, the proposed decomposition is named \ylrev{Tucker-type} interpolative decomposition (Tucker-ID). It is worth noting that there exist several interpolative tensor decomposition algorithms \cite{cai2021mode,mahoney2006tensor,malik2020fast,saibaba2016hoid,minster2020randomized}. However they either use original tensor entries in the factor matrices (instead of the core tensor) \cite{mahoney2006tensor,saibaba2016hoid,cai2021mode} or rely on a different tensor diagram \cite{malik2020fast}. \ylrev{Note that the structure-preserving decomposition in \cite{minster2020randomized} is similar to Tucker-ID but relies on sketching instead of proxy indices for the construction.} As will be seen in \cref{sec:butterflytensor}, the Tucker-ID algorithm is a unique and essential building block of the tensor butterfly algorithm.  
	
	\ylrev{Just like HOSVD, one can easily show that if the approximations in \cref{eqn:unfoldingID} hold true up to a predefined relative compression tolerance $\epsilon$ as 
		\begin{align}
		||\mat{K}^{(k)} - \mat{K}^{(k)}(:,\overline{\tau_k})\mat{U}^k||_F\leq \epsilon||\ts{K}||_F,~~ k=1,\ldots,d,\nonumber\\
		||\mat{K}^{(d+k)} {-} \mat{K}^{(d+k)}(:,\overline{\nu_k})\mat{V}^k||_F\leq \epsilon||\ts{K}||_F,~~ k=1,\ldots,d,
		\end{align} 
		then the Tucker-ID of \cref{eqn:tuckerID} satisfies 
		\begin{equation}
		\bigg|\bigg|\ts{K} {-} \ts{K}(\overline{\midx{\tau}},\overline{\midx{\nu}})\bigg(\prod_{k=1}^{d}\times_k\mat{U}^k\bigg)\bigg(\prod_{k=1}^{d}\times_{d+k}\mat{V}^k\bigg)\bigg|\bigg|_F \leq\epsilon\sqrt{2d}||\ts{K}||_F.    
		\end{equation}
	}
	
	The memory and CPU complexities of Tucker-ID can be briefly analyzed as follows. Assuming that $m_k=n_k=n$ and $\max_{k}|\overline{\tau}_k|$=$\max_{k}|\overline{\nu}_k|=r$ is a constant (we will discuss the case of non-constant $r$ in \cref{sec::compare}), the memory requirement is simply $O(drn+r^{2d})$, where the first and second term represent the storage units for the factor matrices and the core tensor, respectively. The CPU cost for naive computation of Tucker-ID is $O(drn^{2d}+r^{2d})$, where the first term represents the cost of rank-revealing QR of the unfolding matrices in \cref{eqn:unfoldingID}, and the second term represents the cost forming the core tensor $\ts{K}(\overline{\midx{\tau}},\overline{\midx{\nu}})$. In practice, however, the unfolding matrices do not need to be fully formed and one can leverage the idea of proxy rows in \cref{sec:matBF} to reduce the cost for computing the factor matrices to $O(dnr^{2d})$. We will explain this in more detail in the context of the proposed tensor butterfly decomposition algorithm.  
	
	Just like the matrix ID algorithm, Tucker-ID is also not suitable for representing large-sized OIOs as the rank $r$ depends on the size $n$. That said, the Tucker-ID rank is typically significantly smaller than the matrix ID rank, as it exploits more compressibility properties across dimensions by leveraging e.g. translational-invariance or dimensional-separability properties of OIOs; see \cref{sec:rank} for a few of such examples. In what follows, we use Tucker-ID as the building block for constructing a linear-complexity tensor butterfly decomposition algorithm for large-sized OIOs.       
	
	\subsection{Tensor Butterfly Algorithm}\label{sec:butterflytensor}

	Consider a 2d-mode OIO tensor $\ts{K}(\midx{t}^0,\midx{s}^0)$ with $\midx{t}^0=(t^0_1,t^0_1,\ldots,t^0_d)$, $\midx{s}^0=(s^0_1,s^0_1,\ldots,s^0_d)$, $t^0_k=\{1,2,\ldots,m_k\}$, $s^0_k=\{1,2,\ldots,n_k\}$, $k=1,2,\ldots,d$. Without loss of generality, we assume that $m_k= n_k = n$. We further assume that each $t^0_k$ (and $s^0_k$) is binary partitioned with a tree $\mathcal{T}_{t^0_k}$ (and $\mathcal{T}_{s^0_k}$) of $L$ levels for $k=1,2,\ldots,d$.  
	
	To start with, we first define the tensor CLR property as follows: 
	\begin{itemize}[leftmargin=*,noitemsep]
		\item For any level $ 0\leq l \leq L^c$, any multi-set $\midx{\tau}=(\tau_1,\tau_2,\ldots,\tau_d)$ with $\tau_i, i\leq d$ at level $l$ of $\mathcal{T}_{{t^0_i}}$, any multi-set $\midx{s}=(s_1,s_2,\ldots,s_d)$ with $s_i, i\leq d$ at level $L^c$ of $\mathcal{T}_{s^0_i}$, any mode $1 \leq k \leq d$, and any node $\nu$ at level $L^c-l$ of $\mathcal{T}_{s_k}$, the mode-$(d+k)$ unfolding of the subtensor $\ts{K}(\midx{\tau},\midx{s}_{k\leftarrow \nu})$ is numerically low-rank (with rank bounded by $r$), permitting an ID via \cref{eqn:tensorID}:
		\begin{equation}\label{eqn:tensorID_BF_V}
		\ts{K}(\midx{\tau},\midx{s}_{k\leftarrow \nu})\approx\ts{K}(\midx{\tau},\midx{s}_{k\leftarrow \overline{\nu}})\times_{d+k}\mat{V}^{\midx{s},k}_{\midx{\tau},\nu}.
		\end{equation}    	
		\item For any level $ 0\leq l \leq L^c$, any multi-set $\midx{\nu}=(\nu_1,\nu_2,\ldots,\nu_d)$ with $\nu_i, i\leq d$ at level $l$ of $\mathcal{T}_{{s^0_i}}$, any multi-set $\midx{t}=(t_1,t_2,\ldots,t_d)$ with $t_i, i\leq d$ at level $L^c$ of $\mathcal{T}_{t^0_i}$, any mode $1 \leq k \leq d$, and any node $\tau$ at level $L^c-l$ of $\mathcal{T}_{t_k}$, the mode-$k$ unfolding of the subtensor $\ts{K}(\midx{t}_{k\leftarrow \tau},\midx{\nu})$ is numerically low-rank (with rank bounded by $r$), permitting an ID via \cref{eqn:tensorID}: 
		\begin{equation}\label{eqn:tensorID_BF_U}
		\ts{K}(\midx{t}_{k\leftarrow \tau},\midx{\nu})\approx\ts{K}(\midx{t}_{k\leftarrow \overline{\tau}},\midx{\nu})\times_k\mat{U}^{\midx{t},k}_{\tau,\midx{\nu}}.
		\end{equation}    	
	\end{itemize} 
	In essence, the tensor CLR in \cref{eqn:tensorID_BF_V} and \cref{eqn:tensorID_BF_U} investigates the unfolding of judiciously selected subtensors rather than the matricization used in the matrix CLR. Moreover, the tensor CLR requires fixing $d-1$ modes of the $2d$-mode subtensors to be of size $O(\sqrt{n})$ while changing the remaining $d+1$ modes with respect to $l$. Therefore each ID computation can operate on larger subtensors compared to the matrix CLR. In \cref{sec:rank} we provide two examples, namely a free-space Green's function tensor and a high-dimensional Fourier transform,  to explain why the tensor CLR is valid, and in \cref{sec:cc} we will see that the tensor CLR essentially reduces the quasilinear complexity of the matrix butterfly algorithm to linear complexity. Here, assuming that the tensor CLR holds true, we describe the tensor butterfly algorithm. We note that there may be alternative ways to define the tensor CLR different from \cref{eqn:tensorID_BF_V} and \cref{eqn:tensorID_BF_U}, and we leave that as a future work. \ylrev{To avoid notation confusion, we list some notations of the tensor butterfly algorithm in \cref{tab:notation}.}
	
	In what follows, we focus on the computation of $\mat{V}^{\midx{s},k}_{\midx{\tau},\nu}$ (corresponding to the mid-level multi-set $\midx{s}$), as $\mat{U}^{\midx{t},k}_{\tau,\midx{\nu}}$ (corresponding to the mid-level multi-set $\midx{t}$) can be computed in a similar fashion. At level $l=0$, $\mat{V}^{\midx{s},k}_{\midx{\tau},\nu}$ are explicitly formed. At level $0<l\leq L^c$, they are represented in a nested fashion. Let $\midx{p}_{\midx{\tau}}=(p_{\tau_1}, p_{\tau_2}, \ldots, p_{\tau_d})$ consist of parents of $\midx{\tau}=(\tau_1,\tau_2,\ldots,\tau_d)$ in \cref{eqn:tensorID_BF_V}. 
	
	By the tensor CLR property, we have  
	\begin{align}\label{eqn:transfer_computation_tensor_W}
	\ts{K}(\midx{\tau},\midx{s}_{k\leftarrow {\nu}}) &\approx\ts{K}(\midx{\tau},\midx{s}_{k\leftarrow \overline{\nu^1}\cup \overline{\nu^2}})\times_{d+k}\begin{bmatrix}
	\mat{V}^{\midx{s},k}_{\midx{p}_{\midx{\tau}},\nu^1} & \\
	& \mat{V}^{\midx{s},k}_{\midx{p}_{\midx{\tau}},\nu^2}
	\end{bmatrix}\nonumber\\
	&\approx\ts{K}(\midx{\tau},\midx{s}_{k\leftarrow \overline{\nu}})\times_{d+k}\Bigg(\mat{W}^{\midx{s},k}_{\midx{\tau},\nu}\begin{bmatrix}
	\mat{V}^{\midx{s},k}_{\midx{p}_{\midx{\tau}},\nu^1} & \\
	& \mat{V}^{\midx{s},k}_{\midx{p}_{\midx{\tau}},\nu^2}
	\end{bmatrix}\Bigg).
	\end{align}
	
	Comparing \cref{eqn:transfer_computation_tensor_W} and \cref{eqn:tensorID_BF_V}, one realizes that the interpolation matrix $\mat{V}^{\midx{s},k}_{\midx{\tau},\nu}$ is represented as the product of the transfer matrix $\mat{W}^{\midx{s},k}_{\midx{\tau},\nu}$ and $ \text{diag}_c(\mat{V}^{\midx{s},k}_{\midx{p}_{\midx{\tau}},\nu^c})$. Here, the transfer matrix $\mat{W}^{\midx{s},k}_{\midx{\tau},\nu}$ is computed as the interpolation matrix of the column ID of the mode-$(d+k)$ unfolding of $\ts{K}(\midx{\tau},\midx{s}_{k\leftarrow \overline{\nu^1}\cup \overline{\nu^2}})$. As mentioned in \cref{sec:tensor_algorithm}, in practice one never forms the unfolding matrix in full, but instead considers the unfolding of $\ts{K}(\hat{\midx{\tau}},\hat{\midx{s}}_{k\leftarrow \overline{\nu^1}\cup \overline{\nu^2}})$, where $\hat{\midx{\tau}}=(\hat{\tau}_1,\hat{\tau}_2,\ldots,\hat{\tau}_d)$ and $\hat{\midx{s}}=(\hat{s}_1,\hat{s}_2,\ldots,\hat{s}_d)$; here $\hat{\tau}_i$ and $\hat{s}_i$ consist of $O(r)$ judiciously selected indices along modes $i$ and $d+i$, respectively. Note that $\hat{s}_k$ is never used as it is replaced by $\overline{\nu^1}\cup \overline{\nu^2}$ in \cref{eqn:transfer_computation_tensor_W}. The same proxy index strategy can be used to obtain $\mat{V}^{\midx{s},k}_{\midx{\tau},\nu}$ at the level $l=0$. For each $\mat{W}^{\midx{s},k}_{\midx{\tau},\nu}$ or $\mat{V}^{\midx{s},k}_{\midx{\tau},\nu}$, its computation requires $O(r^{2d+1})$ CPU time.

	Similarly in \cref{eqn:tensorID_BF_U}, $\mat{U}^{\midx{t},k}_{\tau,\midx{\nu}}$ is explicitly formed at $l=0$ and constructed via the transfer matrix $\mat{P}^{\midx{t},k}_{\tau,\midx{\nu}}$ at level $0 < l \leq L^c$:
	\begin{align}\label{eqn:transfer_computation_tensor_P}
	\ts{K}(\midx{t}_{k\leftarrow \tau},\midx{\nu}) &\approx \ts{K}(\midx{t}_{k\leftarrow \overline{\tau^1}\cup \overline{\tau^2}},\midx{\nu})\times_{k}\begin{bmatrix}
	\mat{U}^{\midx{t},k}_{\tau^1,\midx{p}_{\midx{\nu}}} & \\
	&\mat{U}^{\midx{t},k}_{\tau^2,\midx{p}_{\midx{\nu}}}
	\end{bmatrix}\nonumber\\
	&\approx\ts{K}(\midx{t}_{k\leftarrow \overline{\tau}},\midx{\nu})\times_{k}\Bigg(\mat{P}^{\midx{t},k}_{\tau,\midx{\nu}}\begin{bmatrix}
	\mat{U}^{\midx{t},k}_{\tau^1,\midx{p}_{\midx{\nu}}} & \\
	&\mat{U}^{\midx{t},k}_{\tau^2,\midx{p}_{\midx{\nu}}}
	\end{bmatrix}\Bigg).
	\end{align}

	Putting together \cref{eqn:tensorID_BF_V}, \cref{eqn:tensorID_BF_U}, \cref{eqn:transfer_computation_tensor_W} and \cref{eqn:transfer_computation_tensor_P}, the proposed tensor butterfly decomposition can be expressed, for any multi-set $\midx{t}=(t_1,t_2,\ldots,t_d)$ with $t_i$ at level $L^c$ of $\mathcal{T}_{t^0_i}$ and any multi-set $\midx{s}=(s_1,s_2,\ldots,s_d)$ with $s_i$ at level $L^c$ of $\mathcal{T}_{s^0_i}$, by forming a Tucker-ID for the ($\midx{t},\midx{s}$) pair:
	\begin{equation}\label{eqn:butterfly_tensor}
	\ts{K}(\midx{t},\midx{s}) \approx \ts{K}(\overline{\midx{t}},\overline{\midx{s}})\bigg(\prod_{k=1}^{d}\times_k\bigg(\prod_{l=L^c}^{1}\dbar{\mat{P}}^{\midx{t},\midx{s},k}_{l}\dbar{\mat{U}}^{\midx{t},k}\bigg)\!\!\bigg)\bigg(\prod_{k=1}^{d}\times_{d+k}\bigg(\prod_{l=L^c}^{1}\dbar{\mat{W}}^{\midx{t},\midx{s},k}_{l}\dbar{\mat{V}}^{\midx{s},k}\bigg)\!\!\bigg).
	\end{equation}

	Here, $\overline{\midx{t}}$ and $\overline{\midx{s}}$ represent the skeleton indices of the Tucker-ID of $\ts{K}(\midx{t},\midx{s})$. 
	The interpolation factors $\dbar{\mat{U}}^{\midx{t},k}$ and $\dbar{\mat{V}}^{\midx{s},k}$ in \cref{eqn:butterfly_tensor} are:
	\begin{align}
	\dbar{\mat{U}}^{\midx{t},k}&=\mathrm{diag}_{\tau}(\mat{U}^{\midx{t},k}_{\tau,\midx{s}^0})
	\label{eqn:tensor_butterfly_factors_U},~~\tau \mathrm{~at~level~} L^c \mathrm{~of~} \mathcal{T}_{t_k},\\
	\dbar{\mat{V}}^{\midx{s},k}&=\mathrm{diag}_{\nu}(\mat{V}^{\midx{s},k}_{\midx{t}^0,\nu})
	\label{eqn:tensor_butterfly_factors_V},~~\nu \mathrm{~at~level~} L^c \mathrm{~of~} \mathcal{T}_{s_k}, 
	\end{align}
	and the transfer factors $\dbar{\mat{P}}_l^{\midx{t},\midx{s},k}$ and $\dbar{\mat{W}}_l^{\midx{t},\midx{s},k}$ for $l=1$, $\ldots$, $L^c$ are:
%\begin{align}
%\dbar{\mat{P}}^{\midx{t},\midx{s},k}_{l}&=\mathrm{diag}_{\midx{\nu}}\Big(\begin{bmatrix}
%\text{diag}_{\tau}(\mat{P}^{\midx{t},k}_{\tau,\midx{\nu}^{\midx{c}}}) 
%\end{bmatrix}_{\midx{c}}\Big), \begin{aligned}~~\tau &\mathrm{~at~level~} L^c-l \mathrm{~of~} \mathcal{T}_{t_k},\\ \nu_i &\mathrm{~at~level~} l-1 \mathrm{~of~} \mathcal{T}_{s^0_i}, s_i \subseteq \nu_i, i\leq d;\end{aligned}\label{eqn:tensor_butterfly_factors_P}\\
%\dbar{\mat{W}}^{\midx{t},\midx{s},k}_{l}&=\mathrm{diag}_{\midx{\tau}}\Big(\begin{bmatrix}
%\text{diag}_{\nu}(\mat{W}^{\midx{s},k}_{\midx{\tau}^{\midx{c}},\nu}) 
%\end{bmatrix}_{\midx{c}}\Big), \begin{aligned}~~\tau_i &\mathrm{~at~level~} l-1 \mathrm{~of~} \mathcal{T}_{t^0_i}, t_i \subseteq \tau_i, i\leq d,\\ \nu &\mathrm{~at~level~} L^c-l \mathrm{~of~} \mathcal{T}_{s_k}.\end{aligned}\label{eqn:tensor_butterfly_factors_W}
%\end{align}
\begin{align}
	\dbar{\mat{P}}^{\midx{t},\midx{s},k}_{l}&=\text{diag}_{\tau}(\mat{P}^{\midx{t},k}_{\tau,\midx{\nu}}), \begin{aligned}~~\tau &\mathrm{~at~level~} L^c-l \mathrm{~of~} \mathcal{T}_{t_k},\\ \nu_i &\mathrm{~at~level~} l \mathrm{~of~} \mathcal{T}_{s^0_i}, s_i \subseteq \nu_i, i\leq d;\end{aligned}\label{eqn:tensor_butterfly_factors_P}\\
	\dbar{\mat{W}}^{\midx{t},\midx{s},k}_{l}&=\text{diag}_{\nu}(\mat{W}^{\midx{s},k}_{\midx{\tau},\nu}), \begin{aligned}~~\tau_i &\mathrm{~at~level~} l \mathrm{~of~} \mathcal{T}_{t^0_i}, t_i \subseteq \tau_i, i\leq d,\\ \nu &\mathrm{~at~level~} L^c-l \mathrm{~of~} \mathcal{T}_{s_k}.\end{aligned}\label{eqn:tensor_butterfly_factors_W}
	\end{align}
	
	One can verify that when $d=1$, the tensor butterfly algorithm \cref{eqn:butterfly_tensor} reduces to the matrix butterfly algorithm \cref{eqn:butterfly_mat}. But when $d>1$, the tensor butterfly algorithm has a distinct algorithmic structure so that the corresponding  computational complexity can be significantly reduced compared with the matrix butterfly algorithm. Detailed computational complexity analysis is provided in \cref{sec:cc}.  
	
	To better understand the structure of the tensor butterfly in \cref{eqn:butterfly_tensor}, \cref{eqn:tensor_butterfly_factors_U}, \cref{eqn:tensor_butterfly_factors_V}, \cref{eqn:tensor_butterfly_factors_P},  and \cref{eqn:tensor_butterfly_factors_W}, we describe its tensor diagram here. We first create the tensor diagram for a \textit{matrix partitioner} as shown in \cref{fig:tensor_diagram_part}(b), which represents a $2^d\times2$ block partitioning of a matrix such as $\begin{bmatrix}\mat{W}^{\midx{s},k}_{\midx{\tau}^{\midx{c}},\nu}\end{bmatrix}_{\midx{c}}$ in \cref{eqn:tensor_butterfly_factors_W} for fixed $\midx{s}, \midx{\tau}, k$ and $\nu$, or $\begin{bmatrix}\mat{P}^{\midx{t},k}_{\tau,\midx{\nu}^{\midx{c}}} \end{bmatrix}_{\midx{c}}$ in \cref{eqn:tensor_butterfly_factors_P} for fixed $\midx{t}, \midx{\nu}, k$ and $\tau$. \ylrev{In \cref{fig:tensor_diagram_part}(b), each of the row and column dimensions is connected to a partitioning node. The row partitioning node has a parent edge with an arrow pointing to the row dimension to be partitioned, and 2 children edges connected to the parent edge. Similarly, the column partitioning node has a parent edge with an arrow pointing to the column dimension to be partitioned, and $2^d$ children edges connected to the parent edge. The weight of the parent edge (i.e., the number of columns and rows of the matrix) equals the sum of the weights of the children edges.} The diagram in \cref{fig:tensor_diagram_part}(c) shows the connectivity for all $\mat{V}^{\midx{s},k}_{\midx{t}^0,\nu}$ (the green circles) and $\begin{bmatrix}\mat{W}^{\midx{s},k}_{\midx{\tau}^{\midx{c}},\nu}\end{bmatrix}_{\midx{c}}$ (the blue circles) for fixed $\midx{s}$ and $k$. The multiplication or contraction of all matrices in \cref{fig:tensor_diagram_part}(c) results in $\mat{V}^{\midx{s},k}_{\midx{t},s_k}$ for all mid-level multi-sets $\midx{t}$, which are of course not explicitly formed.

	\begin{figure}[!tp]
		\centering
		\vspace{-7.5pt}
		\includegraphics[width=\linewidth]{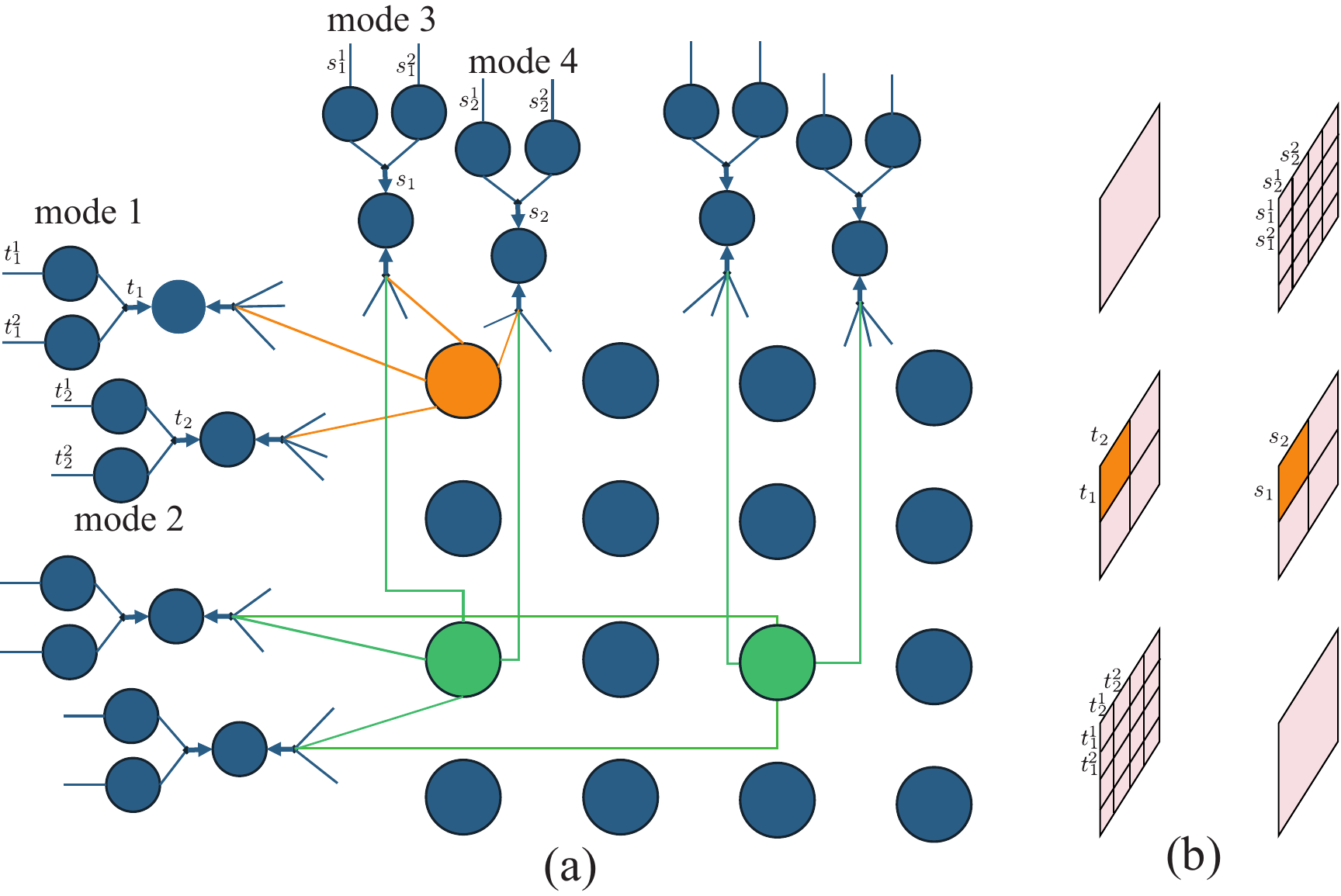}
		\vspace{-5pt}
		\caption{(a) Tensor diagram for the tensor butterfly decomposition of $L=2$ levels of a 4-mode OIO tensor representing (b) high-frequency Green's function interactions between parallel facing 2D unit squares. Only the full connectivity regarding three middle-level node pairs is shown (the two green circles and one orange circle in (a)). The orange circle in (a) represents the core tensor $\ts{K}(\overline{\midx{t}},\overline{\midx{s}})$ for a mid-level pair $(\midx{t},\midx{s})$ with $\midx{t}=(t_1,t_2)$, $\midx{s}=(s_1,s_2)$ highlighted in orange in (b).}
		\label{fig:tensor_diagram}
		%	\vspace{-20pt}
	\end{figure}	
	
	As an example, consider an OIO representing the free-space Green's function interaction between two parallel facing unit square plates in \cref{fig:tensor_diagram}. The tensor is $\ts{K}({\midx{i}},{\midx{j}})=K(x^{\midx{i}},y^{\midx{j}})=\frac{\mathrm{exp}(-\mathrm{i}\omega\rho)}{\rho}$ where $x^{\midx{i}}=(\frac{i_1}{n},\frac{i_2}{n},0)$, $y^{\midx{j}}=(\frac{j_1}{n},\frac{j_2}{n},1)$, $\rho=|x^{\midx{i}}-y^{\midx{j}}|$ and $\omega$ is the wavenumber. Here $1$ represents the distance between the two plates. Consider an $L$=2-level tensor butterfly decomposition, with a total of 16 middle-level multi-set pairs. Let ($\midx{t},\midx{s}$) denote one middle-level multi-set pair with $\midx{t}=(t_1,t_2)$ and $\midx{s}=(s_1,s_2)$ as highlighted in orange in \cref{fig:tensor_diagram}(b). Their children are $t_1^1,t_2^1,t_1^2,t_2^2$ and $s_1^1,s_2^1,s_1^2,s_2^2$. Leveraging the representations in \cref{fig:tensor_diagram_part}(b)-(c), the full diagram for $\ts{K}(\midx{t},\midx{s})$ consists of one 4-mode tensor $\ts{K}(\overline{\midx{t}},\overline{\midx{s}})$ (highlighted in orange in \cref{fig:tensor_diagram}(a)), one transfer matrix per mode, and two factor matrices per mode. In addition, we plot the full connectivity for two other multi-set pairs (highlighted in green in \cref{fig:tensor_diagram}(a)). It is important to note that the factor matrices and transfer matrices are shared among the multi-set pairs. 
	
	The proposed tensor butterfly algorithm is fully described in \cref{alg:tensor_bf_const} for a $2d$-mode tensor $\ts{K}\in \mathbb{C}^{m_1\times m_2\times\ddd\times m_d\times n_1\times n_2\times\ddd\times n_d}$, which  consists of three steps: (1) computation of $\mat{V}^{\midx{s},k}_{\midx{\tau},\nu}$ and $\mat{W}^{\midx{s},k}_{\midx{\tau},\nu}$ starting at Line \ref{line:const_1}, (2) computation of $\mat{U}^{\midx{t},k}_{\tau,\midx{\nu}}$ and $\mat{P}^{\midx{t},k}_{\tau,\midx{\nu}}$ starting at Line \ref{line:const_2}, and (3) computation of $\ts{K}(\overline{\midx{t}},\overline{\midx{s}})$ starting at Line \ref{line:const_3}. We note that, after each $\ts{K}(\overline{\midx{t}},\overline{\midx{s}})$ is formed, we leverage floating-point compression tools such as the ZFP software \cite{zfp} to further compress it.

	Once $\ts{K}$ is compressed, any input tensor $\ts{F}\in \mathbb{C}^{n_1\times n_2\times\ddd\times n_d\times n_v}$ can contract with it to compute $\ts{G}=\ts{K}\times_{d+1,d+2,\ldots,2d}\ts{F}$. It is clear to see that the contraction is equivalent to matrix-matrix multiplication $\mat{G}=\mat{K}\mat{F}$, where $\mat{G} \in \mathbb{C}^{\prod_{k}m_k\times n_v}$, $\mat{K} \in \mathbb{C}^{\prod_{k}m_k\times \prod_{k}n_k}$, and $\mat{F} \in \mathbb{C}^{\prod_{k}n_k\times n_v}$ are matricizations of $\ts{G}$, $\ts{K}$ and $\ts{F}$, respectively, and $n_v$ is the number of columns of $\mat{F}$. The contraction algorithm is described in \cref{alg:tensor_bf_apply} which consists of three steps: 
	\begin{enumerate}[leftmargin=*,noitemsep,label=(\arabic*)]
		\item Contraction with $\mat{V}^{\midx{s},k}_{\midx{\tau},\nu}$ and $\mat{W}^{\midx{s},k}_{\midx{\tau},\nu}$. For each level $l=0,1,\ldots,L^c$, one notices that, since the contraction operation for each multi-set $\midx{\tau}$ with $\tau_i$ at level $l$ of $\mathcal{T}_{t^0_i}$ and the middle-level multi-set $\midx{s}$ is independent of each other, one needs a separate tensor $\ts{F}_{\midx{\tau},\midx{s}}$ to store the contraction result for each multi-set pair $(\midx{\tau},\midx{s})$. $\ts{F}_{\midx{\tau},\midx{s}}$ can be computed by mode-by-mode contraction with the factor matrices $\dbar{\mat{V}}^{\midx{s},k}$ for $l=0$ (Line \ref{line:contract_factor}) and the transfer matrices $\text{diag}_{\nu}(\mat{W}^{\midx{s},k}_{\midx{\tau},\nu})$ for $l>0$ (Line \ref{line:contract_transfer}). 
		\item Contraction with $\ts{K}(\overline{\midx{t}},\overline{\midx{s}})$ at the middle level. Tensors at the middle level $\ts{F}_{\midx{t},\midx{s}}$ are contracted with each subtensor $\ts{K}(\overline{\midx{t}},\overline{\midx{s}})$ separately, resulting in tensors $\ts{G}_{\midx{t},\midx{s}} = \ts{K}(\overline{\midx{t}},\overline{\midx{s}})\times_{d+1,d+2,\ldots,2d}\ts{F}_{\midx{t},\midx{s}}$. 
		\item Contraction with $\mat{U}^{\midx{t},k}_{\tau,\midx{\nu}}$ and $\mat{P}^{\midx{t},k}_{\tau,\midx{\nu}}$. As Step (1), for each level $l=L^c,L^{c-1},\ldots,0$, the contraction operation for each multi-set $\midx{\nu}$ with $\nu_i$ at level $l$ of $\mathcal{T}_{s^0_i}$ and middle-level multi-set $\midx{t}$ is independent. At level $l>0$, the contribution of tensors $\ts{G}_{\midx{t},{\midx{\nu}}}$ is accumulated into $\ts{G}_{\midx{t},\midx{p}_{\midx{\nu}}}$ (Line \ref{line:contract_P}); at level $l=0$, the contraction results are stored in the final output tensor $\ts{G}(\midx{t},1:n_v)$ (Line \ref{line:contract_U}).     
	\end{enumerate}

	\begin{algorithm}[pth!]
		\caption{Construction algorithm for the tensor butterfly decomposition of a $2d$-mode tensor $\ts{K}\in \mathbb{C}^{m_1\times m_2\times\ddd\times m_d\times n_1\times n_2\times\ddd\times n_d}$}
		\label{alg:tensor_bf_const}
		\raggedright
		\textbf{Input:} A function to evaluate a $2d$-mode tensor $\ts{K}(\midx{i},\midx{j})$ for arbitrary multi-indices $(\midx{i},\midx{j})$, binary partitioning trees of $L$ levels $\mathcal{T}_{t^0_k}$ and $\mathcal{T}_{s^0_k}$ with roots $t^0_k=\{1,2,\ldots,m_k\}$ and $s^0_k=\{1,2,\ldots,n_k\}$, a relative compression tolerance $\epsilon$. \\
		\textbf{Output:} Tensor butterfly decomposition of $\ts{K}$: (1) $\mat{V}^{\midx{s},k}_{\midx{\tau},\nu}$ at $l=0$ and $\mat{W}^{\midx{s},k}_{\midx{\tau},\nu}$ at $ 1\leq l \leq L^c$ of $k\leq d$ for multi-set $\midx{\tau}$ with node $\tau_i$ at level $l$ of $\mathcal{T}_{{t^0_i}}$, multi-set $\midx{s}$ with node $s_i$ at level $L^c$ of $\mathcal{T}_{s^0_i}$, and node $\nu$ at level $L^c-l$ of subtree $\mathcal{T}_{s_k}$, (2) $\mat{U}^{\midx{t},k}_{\tau,\midx{\nu}}$ at $l=0$ and $\mat{P}^{\midx{t},k}_{\tau,\midx{\nu}}$ at $ 1\leq l \leq L^c$ of $k\leq d$ for multi-set $\midx{\nu}$ with node $\nu_i$ at level $l$ of $\mathcal{T}_{{s^0_i}}$, multi-set $\midx{t}$ with node $t_i$ at level $L^c$ of $\mathcal{T}_{t^0_i}$, and node $\tau$ at level $L^c-l$ of subtree $\mathcal{T}_{t_k}$, and (3) subtensors $\ts{K}(\overline{\midx{t}},\overline{\midx{s}})$ at $l=L^c$. 
		\begin{algorithmic}[1]
			\State (1) Compute $\mat{V}^{\midx{s},k}_{\midx{\tau},\nu}$ and $\mat{W}^{\midx{s},k}_{\midx{\tau},\nu}$:\label{line:const_1}
			\For{level $l=0,\ldots,L^c$} 
			\For{multi-set $\midx{s}=(s_1,\ldots,s_d)$ with $s_i$ at level $L^c$ of $\mathcal{T}_{s^0_i}$} 	
			\For{multi-set $\midx{\tau}=(\tau_1,\tau_2,\ldots,\tau_d)$ with $\tau_i$ at level $l$ of $\mathcal{T}_{{t^0_i}}$}
			\For{mode index $k=1,\ldots,d$} 	
			\For{node $\nu$ at level $L^c-l$ of $\mathcal{T}_{s_k}$}
			\If{$l = 0$}\Comment{Use \cref{eqn:tensorID_BF_V} with proxies $\hat{\midx{\tau}}$, $\hat{\midx{s}}$ and tolerance $\epsilon$}
			\State Compute $\mat{V}^{\midx{s},k}_{\midx{\tau},\nu}$ and $\overline{\nu}$ via mode-$(d+k)$ unfolding of $\ts{K}(\hat{\midx{\tau}},\hat{\midx{s}}_{k\leftarrow \nu})$
			\Else \Comment{Use \cref{eqn:transfer_computation_tensor_W} with proxies $\hat{\midx{\tau}}$, $\hat{\midx{s}}$ and tolerance $\epsilon$}
			\State Compute $\mat{W}^{\midx{s},k}_{\midx{\tau},\nu}$ and $\overline{\nu}$ via mode-$(d+k)$ unfolding of $\ts{K}(\hat{\midx{\tau}},\hat{\midx{s}}_{k\leftarrow \overline{\nu^1}\cup \overline{\nu^2}})$
			\EndIf 
			\EndFor
			\EndFor
			\EndFor
			\EndFor	
			\EndFor			
			\State (2) Compute $\mat{U}^{\midx{t},k}_{\tau,\midx{\nu}}$ and $\mat{P}^{\midx{t},k}_{\tau,\midx{\nu}}$:\label{line:const_2}
			\For{level $l=0,\ldots,L^c$}	
			\For{multi-set $\midx{t}=(t_1,\ldots,t_d)$ with $t_i$ at level $L^c$ of $\mathcal{T}_{t^0_i}$} 	
			\For{multi-set $\midx{\nu}=(\nu_1,\nu_2,\ldots,\nu_d)$ with $\nu_i$ at level $l$ of $\mathcal{T}_{{s^0_i}}$}
			\For{mode index $k=1,\ldots,d$} 	
			\For{node $\tau$ at level $L^c-l$ of $\mathcal{T}_{t_k}$}
			\If{$l = 0$}\Comment{Use \cref{eqn:tensorID_BF_U} with proxies $\hat{\midx{t}}$, $\hat{\midx{\nu}}$ and tolerance $\epsilon$}
			\State Compute $\mat{U}^{\midx{t},k}_{\tau,\midx{\nu}}$ and $\overline{\tau}$ via mode-$k$ unfolding of $\ts{K}(\hat{\midx{t}}_{k\leftarrow \tau},\hat{\midx{\nu}})$
			\Else \Comment{Use \cref{eqn:transfer_computation_tensor_P} with proxies $\hat{\midx{t}}$, $\hat{\midx{\nu}}$ and tolerance $\epsilon$}
			\State Compute $\mat{P}^{\midx{t},k}_{\tau,\midx{\nu}}$ and $\overline{\tau}$ via mode-$k$ unfolding of $\ts{K}(\hat{\midx{t}}_{k\leftarrow \overline{\tau^1}\cup \overline{\tau^2}},\hat{\midx{\nu}})$
			\EndIf 
			\EndFor
			\EndFor
			\EndFor
			\EndFor	
			\EndFor		
			\State (3) Compute $\ts{K}(\overline{\midx{t}},\overline{\midx{s}})$:\label{line:const_3}
			\For{multi-set $\midx{s}=(s_1,\ldots,s_d)$ with $s_i$ at level $L^c$ of $\mathcal{T}_{s^0_i}$} 	
			\For{multi-set $\midx{t}=(t_1,\ldots,t_d)$ with $t_i$ at level $L^c$ of $\mathcal{T}_{t^0_i}$} 	
			\State Compute $\ts{K}(\overline{\midx{t}},\overline{\midx{s}})$ and ZFP compress it
			\EndFor	
			\EndFor			
		\end{algorithmic}
	\end{algorithm}

	\begin{algorithm}[pth!]
		\caption{Contraction algorithm for a tensor butterfly decomposition with an input tensor}
		\label{alg:tensor_bf_apply}
		\raggedright
		\textbf{Input:} The tensor butterfly decomposition of a $2d$-mode tensor $\ts{K}\in \mathbb{C}^{m_1\times m_2\times\ddd\times m_d\times n_1\times n_2\times\ddd\times n_d}$, and a (full) $d+1$-mode input tensor $\ts{F}\in \mathbb{C}^{n_1\times n_2\times\ddd\times n_d\times n_v}$ where $n_v$ denotes the number of columns of $\mat{F}^{(d+1)}$.\\
		\textbf{Output:} The $d+1$-mode output tensor $\ts{G}=\ts{K}\times_{d+1,d+2,\ldots,2d}\ts{F}$ where $\ts{G}\in \mathbb{C}^{m_1\times m_2\times\ddd\times m_d\times n_v}$. 
		\begin{algorithmic}[1]
			\State (1) Multiply with $\mat{V}^{\midx{s},k}_{\midx{\tau},\nu}$ and $\mat{W}^{\midx{s},k}_{\midx{\tau},\nu}$:
			\For{level $l=0,\ldots,L^c$} 
			\For{multi-set $\midx{s}=(s_1,s_2\ldots,s_d)$ with $s_i$ at level $L^c$ of $\mathcal{T}_{s^0_i}$} 	
			\For{multi-set $\midx{\tau}=(\tau_1,\tau_2,\ldots,\tau_d)$ with $\tau_i$ at level $l$ of $\mathcal{T}_{{t^0_i}}$}
			\If{$l = 0$}
			\State $\ts{F}_{\midx{\tau},\midx{s}} = \ts{F}(\midx{s},1:n_v)\prod_{k=1}^{d}\times_k\dbar{\mat{V}}^{\midx{s},k}$\label{line:contract_factor}
			\Else
			\State $\ts{F}_{\midx{\tau},\midx{s}} = \ts{F}_{\midx{p}_{\midx{\tau}},\midx{s}}\prod_{k=1}^{d}\times_k\text{diag}_{\nu}(\mat{W}^{\midx{s},k}_{\midx{\tau},\nu})$ \Comment{$ \nu \mathrm{~at~level~} L^c-l \mathrm{~of~} \mathcal{T}_{s_k}$}	\label{line:contract_transfer}
			\EndIf
			\EndFor
			\EndFor	
			\EndFor		
			\State (2) Contract with $\ts{K}(\overline{\midx{t}},\overline{\midx{s}})$:	
			\For{multi-set $\midx{t}=(t_1,t_2\ldots,t_d)$ with $t_i$ at level $L^c$ of $\mathcal{T}_{t^0_i}$} 	
			\For{multi-set $\midx{s}=(s_1,s_2\ldots,s_d)$ with $s_i$ at level $L^c$ of $\mathcal{T}_{s^0_i}$} 	
			\State ZFP decompress $\ts{K}(\overline{\midx{t}},\overline{\midx{s}})$ and compute $\ts{G}_{\midx{t},\midx{s}} = \ts{K}(\overline{\midx{t}},\overline{\midx{s}})\times_{d+1,d+2,\ldots,2d}\ts{F}_{\midx{t},\midx{s}}$
			\EndFor	
			\EndFor				
			\State (3) Multiply with $\mat{U}^{\midx{t},k}_{\tau,\midx{\nu}}$ and $\mat{P}^{\midx{t},k}_{\tau,\midx{\nu}}$:
			\For{level $l=L^c,\ldots,0$} 
			\For{multi-set $\midx{t}=(t_1,t_2,\ldots,t_d)$ with $t_i$ at level $L^c$ of $\mathcal{T}_{t^0_i}$} 	
			\For{multi-set $\midx{\nu}=(\nu_1,\nu_2,\ldots,\nu_d)$ with $\nu_i$ at level $l$ of $\mathcal{T}_{{s^0_i}}$}
			\If{$l = 0$} \Comment{Compute and return $\ts{G}$}
			\State $\ts{G}(\midx{t},1:n_v) =\ts{G}_{\midx{t},\midx{\nu}}\prod_{k=1}^{d}\times_k\dbar{\mat{U}}^{\midx{t},k}$\label{line:contract_U}
			\Else
			\State $\ts{G}_{\midx{t},\midx{p}_{\midx{\nu}}} \mathrel{+}= \ts{G}_{\midx{t},\midx{\nu}}\prod_{k=1}^{d}\times_k\text{diag}_{\tau}(\mat{P}^{\midx{t},k}_{\tau,\midx{\nu}})$ \Comment{$ \tau \mathrm{~at~level~} L^c-l \mathrm{~of~} \mathcal{T}_{t_k}$}\label{line:contract_P}	
			\EndIf
			\EndFor
			\EndFor	
			\EndFor		
			
		\end{algorithmic}
	\end{algorithm}

	\subsubsection{Rank Estimate}\label{sec:rank}					
\begin{figure}[!tp]
		\centering
		\vspace{-7.5pt}
		\includegraphics[width=0.9\linewidth]{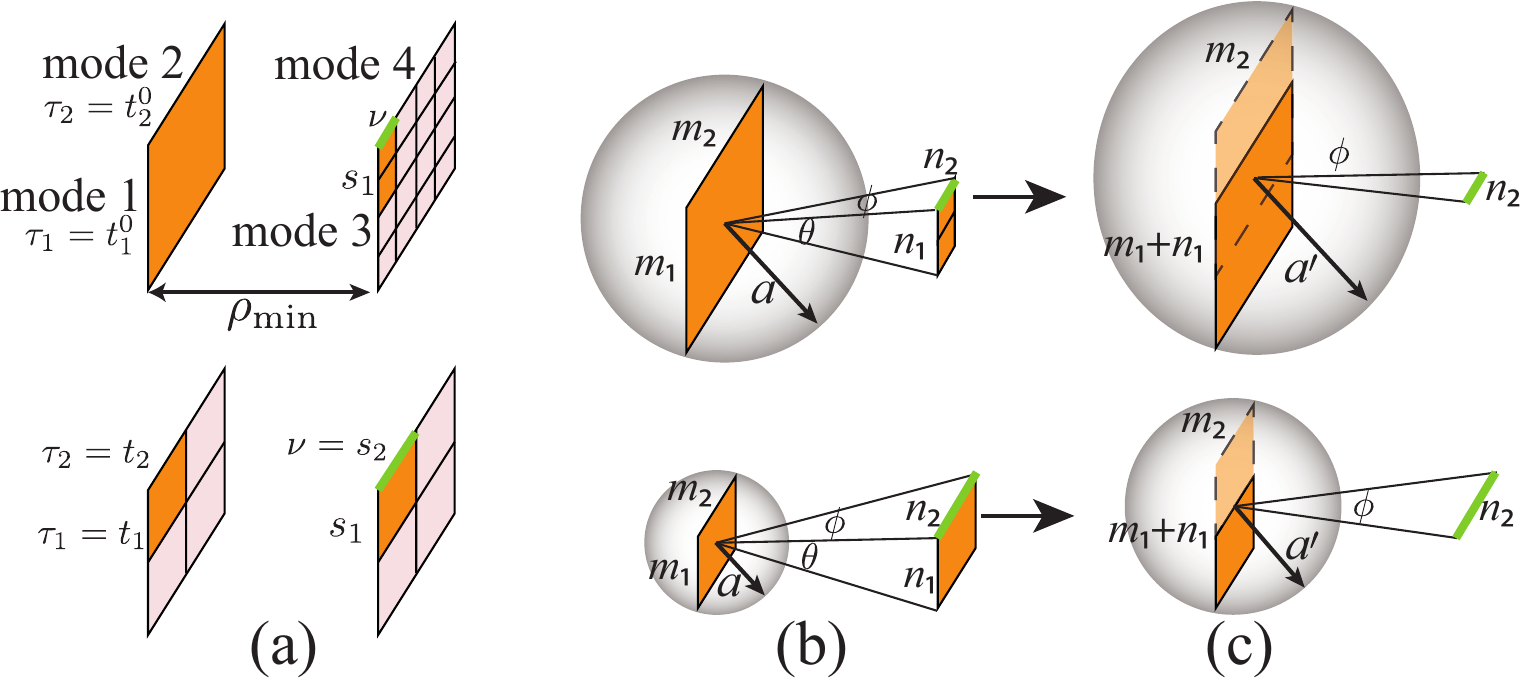}
		\vspace{-5pt}
		\caption{Illustration of the tensor CLR property with $L=2$ for a 4-mode tensor representing free-space Green's function interactions between parallel facing unit square plates. (a) The target and source domains are partitioned at $l=0$ (top) and $l=L^c=1$ (bottom) with a multi-set pair $(\midx{\tau},\midx{s}_{k\leftarrow \nu})$ highlighted in orange for the skeletonization along mode 4. The sizes of the nodes are $|\tau_1|=m_1$,$|\tau_2|=m_1$, $|s_1|=n_1$ and $|\nu|=n_2$. (b) Illustration of the rank of the matricization of  $\ts{K}(\midx{\tau},\midx{s}_{k\leftarrow \nu})$ \ylrev{used in the matrix butterfly algorithm. Here $a$ is the radius of the sphere enclosing the target domain of physical sizes ${m_1}\delta_x\times m_2\delta_y$. $\theta\approx\frac{n_1}{n\rho_{\min}}$, $\phi\approx\frac{n_2}{n\rho_{\min}}$, and the product $\theta\phi$ represents the solid angle covered by the source domain as seen from the center of the target domain.} (c) Illustration of the rank of the mode-4 unfolding of $\ts{K}(\midx{\tau},\midx{s}_{k\leftarrow \nu})$ \ylrev{used in the tensor butterfly algorithm. Here, $a'$ is the radius of the sphere enclosing the enlarged target domain. The source domain is reduced to a line segment of length ${n_2}\delta_y$.}}
		\label{fig:rank_estimate}
		%	\vspace{-20pt}
	\end{figure}

	In this subsection, we use two specific high-dimensional examples, namely high-frequency free-space Green's functions for wave equations and uniform discrete Fourier transforms (DFTs) to investigate the matrix and tensor CLR properties, and compare the matrix and tensor butterfly ranks $r_m$ and $r_t$, respectively. For the Green's function example, the tensor CLR property is a result of matrix CLR and translational invariance, and $r_t$ is much smaller than $r_m$; for the DFT example, the tensor CLR property is a result of matrix CLR and dimensionality separability, and $r_t$ is exactly the same as $r_m$ of 1D DFTs. For more-general OIOs, such as analytical and numerical Green's functions for inhomogeneous media, Radon transforms, non-uniform DFTs, and general Fourier integral operators, rigorous rank analysis is non-trivial and we rely on numerical experiments in \cref{sec:example} to demonstrate the efficacy of the tensor butterfly algorithm.      
	\paragraph{High-frequency Green's functions} We use an example similar to the one used in \cref{sec:butterflytensor}. Consider an OIO representing the free-space Green's function interaction between two parallel-facing 
	unit-square plates. The $n\times n\times n \times n$ tensor is \begin{equation}\ts{K}({\midx{i}},{\midx{j}})=K(x^{\midx{i}},y^{\midx{j}})=\frac{\mathrm{exp}(-\mathrm{i}\omega\rho)}{\rho},\end{equation} where $x^{\midx{i}}=(\frac{i_1}{n},\frac{i_2}{n},0)$, $y^{\midx{j}}=(\frac{j_1}{n},\frac{j_2}{n},\rho_{\min})$, $\omega$ is the wavenumber, and $\rho=|x^{\midx{i}}-y^{\midx{j}}|$. Here $\rho_{\min}$ represents the distance between the two plates assumed to be sufficiently large. In the high-frequency setting, $n=C_p\omega$ with a constant $C_p$ independent of $n$ and $\omega$, and the grid size is $\delta_x=\delta_y=\frac{1}{n}$ per dimension. It has been well studied \cite{Eric_1994_butterfly,michielssen_multilevel_1996,engquist2007fast,bucci1989degrees} that for any multi-set pair $(\midx{\tau},\midx{\nu})$ (\ylrev{assuming that each set of the multi-set $\midx{\tau}$ or $\midx{\nu}$ contains contiguous indices}) leading to a subtensor $\ts{K}(\midx{\tau},\midx{\nu})$ of sizes $m_1\times m_2\times n_1 \times n_2$ with $m_i, n_i\leq n$, the numerical rank of its matricization $\mat{K}\in\mathbb{C}^{m_1m_2\times n_1n_2}$ can be estimated as 
	\begin{equation}
	r_m \approx \omega^2a^2\theta\phi\ylrev{+\Delta_\epsilon} \approx \frac{\omega^2a^2n_1n_2}{n^2\rho_{\min}^2}\ylrev{+\Delta_\epsilon}.\label{eqn:dof}
	\end{equation}  
	Here $a$ is the radius of the sphere enclosing the target domain of physical sizes ${m_1}\delta_x\times m_2\delta_y$. $\theta\approx\frac{n_1}{n\rho_{\min}}$, $\phi\approx\frac{n_2}{n\rho_{\min}}$, and the product $\theta\phi$ represents the solid angle covered by the source domain as seen from the center of the target domain. Note that $\frac{\omega a}{\rho_{\min}}$ approximately represents the Nyquist sampling rate per direction needed in the source domain. \ylrev{The $\epsilon$-dependent term $\Delta_\epsilon=O(\log \epsilon^{-1})$ according to analysis in \cite{Eric_1994_butterfly,michielssen_multilevel_1996}.} The matrix and tensor butterfly ranks can be estimated as follows:
	\begin{itemize}[leftmargin=*,noitemsep]
		\item \textit{Matrix butterfly rank:} Consider a matrix butterfly factorization of matricization of $\ts{K}$. By design, for any node pair at each level, $m_1n_1=m_2n_2=C_bn$, where $C_b^2$ represents the size of the leaf nodes. Therefore, the matrix butterfly rank can be estimated from \cref{eqn:dof} as  
		\begin{equation}
		r_m \approx \frac{C_b^2}{2C_p^2\rho_{\min}^2}\ylrev{+\Delta_\epsilon}.\label{eqn:r_m}
		\end{equation}   
		Here we have assumed $a=\frac{m_1}{\sqrt{2}n}$. Note that $r_m$ is a constant independent of $n$, and therefore the matrix CLR property holds true. 
		\item \textit{Tensor butterfly rank:} Consider an $L$-level tensor butterfly factorization of $\ts{K}$. We just need to check the tensor rank, e.g., the rank of the mode-4 unfolding of the corresponding subtensors at Step (1) of \cref{alg:tensor_bf_const}, as the unfolding for the other modes can be investigated in a similar fashion. \Cref{fig:rank_estimate}(a) shows an example of $L=2$, where the target and source domains are partitioned at $l=0$ (top) and $l=L^c=1$ (bottom) at Step (1) of \cref{alg:tensor_bf_const}. Consider a multi-set pair $(\midx{\tau},\midx{s}_{k\leftarrow \nu})$ with $k=4$ required by the tensor CLR property in \cref{eqn:tensorID_BF_V}. \Cref{fig:rank_estimate}(a) highlights in orange one multi-set pair at $l=0$ (top) and one multi-set pair at $l=L^c$ (bottom). Mode 4 is highlighted in \ylrev{green (in all subfigures of \cref{fig:rank_estimate})}, which needs to be skeletonized by ID. By \cref{eqn:dof}, the rank of the matricization of $\ts{K}(\midx{\tau},\midx{s}_{k\leftarrow \nu})$ is no longer a constant as the tensor butterfly algorithm needs to keep $n_1=|s_1|=n/2^{L^c}$ (see \cref{fig:rank_estimate}(b)). However, due to translational invariance of the free-space Green's function, i.e., $K(x^{\midx{i}},y^{\midx{j}})=K(\tilde{x},\tilde{y})$, where $\tilde{x}=(0,\frac{i_2}{n},0)$, $\tilde{y}=(\frac{j_1-i_1}{n},\frac{j_2}{n},\rho_{\min})$, the mode-4 unfolding of $\ts{K}(\midx{\tau},\midx{s}_{k\leftarrow \nu})$ is the matrix representing the Green's function interaction between an enlarged target domain of sizes ${(m_1+n_1)}\delta_x\times m_2\delta_y$ and a source line segment of length ${n_2}\delta_y$. Therefore its rank (hence the tensor rank) can be estimated as 
		\begin{equation}
		r_t \approx \omega a'\phi\ylrev{+\Delta_\epsilon} \approx \frac{\omega a'n_2}{n\rho_{\min}}\ylrev{+\Delta_\epsilon}\leq \frac{\sqrt{2}C_b}{C_p\rho_{min}}\ylrev{+\Delta_\epsilon},\label{eqn:r_t}
		\end{equation}  
		where $a'$ is the radius of the sphere enclosing the enlarged target domain and $\frac{\omega a'}{\rho_{\min}}$ approximately represents the Nyquist sampling rate on the source line segment. The last inequality is a result of $a'\approx\frac{m_1+n_1}{\sqrt{2}n}\leq\frac{\sqrt{2}m_1}{n}$ and $m_2n_2=C_bn$. Here, the critical condition $n_1\leq m_1$ is a direct result of the setup of the tensor CLR in \cref{eqn:tensorID_BF_V}: $l\leq L^c$ and $n_1=|s_1|=n/2^{L^c}$ (i.e., $s_1$ is fixed as the middle level set as $l$ changes). One can clearly see from \cref{eqn:r_t} that $r_t$ is independent of $n$, and thus the tensor CLR property holds true.
	\end{itemize}  
	
	We remark that the tensor butterfly rank $r_t$ in \cref{eqn:r_t} is significantly smaller than the matrix butterfly rank $r_m$ in \cref{eqn:r_m} with $r_t \approx 2\sqrt{r_m}$. One can perform similar analysis of $r_m$ and $r_t$ for different geometrical settings, such as a pair of well-separated 3D unit cubes, or a pair of co-planar 2D unit-square plates. We leave these exercises to the readers.

	\paragraph{Discrete Fourier Transform} Our second example is the high-dimensional discrete Fourier transform (DFT) defined by \begin{equation}\ts{K}({\midx{i}},{\midx{j}})=\mathrm{exp}(2\pi\mathrm{i}x^{\midx{i}}\cdot y^{\midx{j}})\end{equation} with $x^{\midx{i}} = (i_1-1,i_2-1,\ldots, i_d-1)$ and $y^{\midx{j}} = (\frac{j_1-1}{n},\frac{j_2-1}{n},\ldots, \frac{j_d-1}{n})$. We first notice that, since \begin{equation}\mathrm{exp}(2\pi\mathrm{i}x^{\midx{i}}\cdot y^{\midx{j}})=\prod_{k=1}^{d}\mathrm{exp}\left(\frac{2\pi\mathrm{i}(i_k-1)(j_k-1)}{n}\right),\end{equation} to carry out arbitrary high-dimensional DFTs one can simply perform 1D DFTs one dimension at a time (while fixing the indices of the other dimensions) by either 1D FFT or 1D matrix butterfly algorithms. We choose the 1D butterfly approach as our reference algorithm. For each node pair at dimension $k$ discretized into a $m_k\times n_k$ matrix, we assume that $m_kn_k=C_bn$. It has been proved in \cite{Candes_butterfly_2009,Lexing_SFT_2009} that this leads to the matrix CLR property and each 1D DFT (fixing indices in other dimensions) can be computed by the matrix butterfly algorithm in $O(n\log n)$ time with a constant butterfly rank $r_m$. Overall this approach requires $O(dn^d\log n)$ operations. 
	
	In contrast, the tensor butterfly algorithm relies on direct compression of e.g., mode-$k$ unfolding of subtensors $\ts{K}(\midx{\tau},\midx{s}_{k\leftarrow \nu})$. Consider any submatrix $\mat{K}_{sub}\in\mathbb{C}^{m_k\times n_k}$ of this unfolding matrix $\mat{K}^{(k)}$; by fixing $i_p$ and $j_p$ with $p\neq k$, its entry is simply 
	\begin{equation*}\mathrm{exp}\left(\frac{2\pi\mathrm{i}(i_k-1)(j_k-1)}{n}\right)\end{equation*} 
	scaled by a constant factor 
	\begin{equation*}\prod_{p\neq k}\mathrm{exp}\left(\frac{2\pi\mathrm{i}(i_p-1)(j_p-1)}{n}\right)
	\end{equation*} 
	of modulus 1. Therefore the tensor butterfly rank is \begin{equation}r_t=\mathrm{rank}(\mat{K}^{(k)})=\mathrm{rank}(\mat{K}_{sub})=r_m.\end{equation} The tensor CLR property thus holds true, and the tensor rank is exactly the same as the 1D butterfly algorithm per dimension. However, as we will see \cref{sec:cc}, our tensor butterfly algorithm yields a linear instead of quasi-linear CPU complexity for high-dimensional DFTs.

	\subsubsection{Complexity Analysis}\label{sec:cc}			
	Here we provide an analysis of computational complexity and memory requirement of the proposed construction algorithm (\cref{alg:tensor_bf_const}) and contraction algorithm (\cref{alg:tensor_bf_apply}), assuming that the tensor butterfly rank $r_t$ is a small constant and $d>1$. \ylrev{Recall that the $2d$-mode tensor $\ts{K}$ has size $n$ and a binary tree ($\mathcal{T}_{{t^0_k}}$ or $\mathcal{T}_{s^0_k}$) of $L$ levels along each mode $k$. $L^c=L/2$ denotes the middle level. We refer the readers to \cref{tab:notation} to recall the notations of the multi-set, $k^{th}$ set, mid-level subtensor, transfer matrix, and interpolation matrix, etc.} 
	
	At Step (1) of \cref{alg:tensor_bf_const}, each level $1\leq l\leq L^c$ has $\#\midx{s}=O(\sqrt{n}^d)$, $\#\midx{\tau}=2^{dl}$, $\#\nu=O(\sqrt{n}/2^{l})$ for each mode $k\leq d$. Each $\mat{W}^{\midx{s},k}_{\midx{\tau},\nu}$ requires $O(r_t^2)$ storage, and $O(r_t^{2d+1})$ computational time when proxy indices $\hat{\midx{\tau}}$, $\hat{\midx{s}}$ are being used. The storage requirement and computational cost for $\mat{W}^{\midx{s},k}_{\midx{\tau},\nu}$ are:
	\begin{align}
	\mathrm{mem}_W &= \sum_{l=1}^{L^c}dO(\sqrt{n}^d)2^{dl}O(\sqrt{n}/2^{l})O(r_t^2)=O(dn^dr_t^2),\\
	\mathrm{time}_W &= \sum_{l=1}^{L^c}dO(\sqrt{n}^d)2^{dl}O(\sqrt{n}/2^{l})O(r_t^{2d+1})=O(dn^dr_t^{2d+1}).
	\end{align} 
	One can easily verify that the computation and storage of $\mat{V}^{\midx{s},k}_{\midx{\tau},\nu}$ at $l=0$ is less dominant than $\mat{W}^{\midx{s},k}_{\midx{\tau},\nu}$ at $l>0$ and we skip its analysis. 
	
	At Step (2) of \cref{alg:tensor_bf_const}, we have $\#\midx{s}=O(\sqrt{n}^d)$ and $\#\midx{t}=O(\sqrt{n}^d)$, and each $\ts{K}(\overline{\midx{t}},\overline{\midx{s}})$ requires $O(r_t^{2d})$ computation time and storage units (even if it is further ZFP compressed to reduce storage requirement), which adds up to 
	\begin{align}
	\mathrm{mem}_K &= O(\sqrt{n}^d) O(\sqrt{n}^d) O(r_t^{2d}) = O(n^dr_t^{2d}),\\
	\mathrm{time}_K &= O(\sqrt{n}^d) O(\sqrt{n}^d) O(r_t^{2d}) = O(n^dr_t^{2d}).
	\end{align}  
	
	Step (3) of \cref{alg:tensor_bf_const} has similar computational cost and memory requirement to Step (1) when contracting with the intermediate matrices $\mat{P}^{\midx{t},k}_{\tau,\midx{\nu}}$, with $\mathrm{mem}_P\sim\mathrm{mem}_W$ and $\mathrm{time}_P\sim\mathrm{time}_W$. 
	
	Overall, \cref{alg:tensor_bf_const} requires
	\begin{align}
	\mathrm{mem} & =  \mathrm{mem}_W + \mathrm{mem}_K + \mathrm{mem}_P =O(n^dr_t^{2d}),\\
	\mathrm{time} &= \mathrm{time}_W + \mathrm{time}_K + \mathrm{time}_P = O(dn^dr_t^{2d+1}).
	\end{align}  
	
	Following a similar analysis, one can estimate the computational cost of \cref{alg:tensor_bf_apply} as $O(n^dr_t^{2d}n_v)$, which is essentially of the similar order as $\mathrm{mem}$ of \cref{alg:tensor_bf_const}, except an extra factor $n_v$ representing the size of the last dimension of the input tensor. 
	
	One critical observation is that the time and storage complexity of the tensor butterfly algorithm is \textit{linear} in $n^d$ with smaller ranks $r_t$, while that of the matrix butterfly algorithm is \textit{quasi-linear} in $n^d$ with much larger ranks $r_m$. This leads to a significantly superior algorithm, as will be demonstrated with the numerical results in \cref{sec:example}. That being said, one can verify that there is no difference between the two algorithms when $d=1$.

	\begin{table*} %[!tp]
		\centering	
		\begin{tabular}{|c|cc|cc|cc|}
			\hline		
			& \multicolumn{2}{c}{Factor time}  & \multicolumn{2}{|c}{Apply time} & \multicolumn{2}{|c|}{$r$} \\
			Algorithm & $d=2$ & $d=3$ & $d=2$ & $d=3$ & $d=2$ & $d=3$ \\
			\hline		
			Tensor butterfly & $n^2$ & $n^3$ & $n^2$& $n^3$& $1$ & $1$  \\	
			\hline	 		
			Matrix butterfly & $n^2\log n$ & $n^3\log n$ & $n^2\log n$& $n^3\log n$& $1$ & $1$  \\	
			\hline						
			Tucker-ID & $n^4$ & $n^{4}-n^{6*}$ & $n^4$& $n^{4}-n^{6*}$& $n$ & $n$  \\	
			\hline	 			
			QTT (Green\&Radon) & $n^{3}\log n$ & $n^{3}\log n$ & $n^{4}\log n$& $n^{5}\log n$& $n$ & $n$  \\
			\hline	 			
			QTT (DFT) & $\log n$ & $\log n$ & $n^{2}\log n$& $n^{3}\log n$& $1$ & $1$  \\ 			
			\hline		
		\end{tabular}
		\caption{CPU complexity of the tensor butterfly algorithm, matrix butterfly algorithm, Tucker-ID and QTT when applied to high-frequency Green's functions ($d=2$ represents two parallel facing unit square plates and $d=3$ represents two separated unit cubes), DFT and Radon transforms. Here we assume that tensor butterfly, matrix butterfly and Tucker-ID algorithms use proxy indices, and the QTT algorithm uses TT-cross. The big O notation is assumed. *: for $d=3$, the complexity of Tucker-ID is $n^6$ for Radon transform and DFT, and $n^4$ for Green's function.}\label{tab:cc}
	\end{table*}

	\subsubsection{Comparison with Tucker-ID and QTT}\label{sec::compare}		
	Here we make a comparison of the computational complexities of the matrix butterfly algorithm, tensor butterfly algorithm, Tucker-ID and QTT for several frequently encountered OIOs with $d=2,3$, namely Green's functions for high-frequency wave equations (where $d=2$ represents two parallel facing unit square plates and $d=3$ represents two separated unit cubes), Radon transforms (a type of Fourier integral operators), and DFT. We first summarize the computational complexities of the factorization and application of matrix and tensor butterfly algorithms in \cref{tab:cc}. Here we use $r$ to denote the maximum rank of the submatrices or (unfolding and matricization of) subtensors associated with each algorithm. In other words, we drop the subscript of $r_m$ and $r_t$ in this subsection. We note that $r=O(1)$ for butterfly algorithms, and the computational complexity for matrix and tensor butterfly algorithms is, respectively, $O(dn^d\log n)$ and $O(dn^d)$, for all OIOs considered here.     
	
	The Tucker-ID algorithm in \cref{sec:tuckerID} (even with the use of proxy indices to accelerate the factorization), always leads to $r=O(n)$ for OIOs and hence almost always $O(n^{2d})$ factorization and application complexities (see \cref{tab:cc}). One exception is perhaps the Green's function for $d=3$,  where one can easily show that 4 out of the 6 unfolding matrices have a rank of $O(n)$ and the remaining 2 have a rank of $O(1)$, leading to the $O(n^4)$ computational complexity. Overall, we remark that \ylrev{Tucker-type} decomposition algorithms are typically the least efficient tensor algorithms for OIOs. 
	
	The QTT algorithm, on the other hand, is a more subtle algorithm to compare with. Assuming that the maximum rank among all steps in QTT is $r$, we first summarize the computational complexities of the factorization and application of QTT. For factorization, we only consider the TT-cross type of algorithms, which yields the best known computational complexity among all TT-based algorithms. The computational complexity of TT-cross is $O(dr^3\log n)$ \cite{corona2017tensor,oseledets2011tensor}. Once factorized, the application cost of the QTT factorization with a full input tensor is $O(dr^2n^d\log n)$ \cite{corona2017tensor}. This complexity can be reduced to $O(dr^2r_i^2\log n)$ when the input tensor is also in the QTT format with TT rank $r_i$. However, an arbitrary input tensor can have a TT rank up to $r_i=O(n^{d/2})$ (which leads to the same application cost as contraction with a full input tensor). Therefore in our comparative study, we stick with the $O(dr^2n^d\log n)$ application complexity.  
	
	For high-frequency Green's functions and general-form Fourier integral operators (e.g. Radon transforms), the TT rank in general behaves as $r=O(n)$ \cite{corona2017tensor}, leading to a factorization cost of $O(dn^3\log n)$ and an application cost of $O(dn^{2+d}\log n)$, as detailed in \cref{tab:cc}. It is worth mentioning that, treating DFTs as a special type of Fourier integral operators, QTT can achieve $r=O(1)$ when a proper bit-reversal ordering is used \cite{chen2024direct}, leading to a factorization cost of $O(d\log n)$ and an application cost of $O(dn^{d}\log n)$, as shown in \cref{tab:cc}. In contrast, the proposed tensor butterfly algorithm can always yield $O(dn^d)$ factorization and $O(n^d)$ application costs.

	\section{Numerical Results}\label{sec:example}
	This section provides several numerical examples to demonstrate the accuracy and efficiency of the proposed tensor butterfly algorithm when applied to large-scale and high-dimensional OIOs including Green's function tensors for high-frequency Helmholtz equations (\cref{sec:example_green}), Radon transform tensors (\cref{sec:example_radon}), and high-dimensional DFTs (\cref{sec:example_DFT}). We compare our algorithm with a few existing matrix and tensor algorithms including the matrix butterfly algorithm in \cref{sec:matBF}, the Tucker-ID algorithm in \cref{sec:tuckerID}, the QTT algorithm \cite{oseledets2011tensor}, the FFT algorithm implemented in the heFFTe package \cite{heFFTe}, and the non-uniform FFT (NUFFT) algorithm implemented in the FINUFFT package \cite{finufft}. All of these algorithms except for Tucker-ID \ylrev{(sequential implementation in Fortran2008 via the ButterflyPACK package \cite{osti_1564244})} and FINUFFT \ylrev{(Python interface to the C backend with shared-memory parallelism)} are tested in distributed-memory parallelism. \ylrev{The reference binary-tree-based matrix butterfly algorithm in \cref{sec:matBF} is implemented in Fortran2008 with distributed-memory parallelism \cite{liu2021butterfly}, available in the ButterflyPACK package \cite{osti_1564244}. The proposed tensor butterfly algorithm is also available in the ButterflyPACK package with distributed-memory parallelism (which will be described in detail in a future paper). The matrix and tensor butterfly algorithms leverage ZFP to further compress the middle-level submatrices and subtensors, respectively}. It is worth noting that currently there is no single package that can both compute and apply the QTT decomposition in distributed-memory parallelism. In our tests, we perform the factorization using a distributed-memory TT code \ylrev{(fully Python)} \cite{shi2024distributed} that parallelizes a cross interpolation algorithm \cite{dolgov2020parallel}, and then we implement the distributed-memory QTT contraction via the CTF package \ylrev{(Python interface to the C++ backend)} \cite{ctf}. All experiments are performed using 4 CPU nodes of the Perlmutter machine at NERSC in Berkeley, where each node has two $64$-core AMD EPYC 7763 processors and 128GB of 2133MHz DDR4 memory. 
	
	\subsection{Green's functions for high-frequency Helmholtz equations}\label{sec:example_green}
	In this subsection, we consider the tensor discretized from 3D free-space Green's functions for high-frequency Helmholtz equations. Specifically, the tensor entry is  
	\begin{eqnarray}
	\ts{K}({\midx{i}},{\midx{j}})=\frac{\mathrm{exp}(-\mathrm{i}\omega\rho)}{\rho},~~~~ \rho=|x^{\midx{i}}-y^{\midx{j}}|,  
	\end{eqnarray}
	where $\omega$ represents the wave number. Two tests are performed: (1) A 4-way tensor representing the Green's function interaction between two parallel facing unit plates with distance 1, i.e., $x^{\midx{i}}=(\frac{i_1}{n},\frac{i_2}{n},0)$, $y^{\midx{j}}=(\frac{j_1}{n},\frac{j_2}{n},1)$, and $d=2$. (2) A 6-way tensor representing the Green's function interaction between two unit cubes with the distance between their centers set to 2, i.e., $x^{\midx{i}}=(\frac{i_1}{n},\frac{i_2}{n},\frac{i_3}{n})$, $y^{\midx{j}}=(\frac{j_1}{n},\frac{j_2}{n},\frac{j_3}{n}+2)$, and $d=3$. For both tests, the wave number is chosen such that the number of points per wave length is 4, i.e., $2\pi n/\omega=4$ or $C_p={2}/{\pi}$. We first perform compression using the tensor butterfly, Tucker-ID and QTT algorithms, and then perform application/contraction using a random input tensor $\ts{F}$. We also add results for the matrix butterfly algorithm using the corresponding matricization of $\ts{K}$ and $\ts{F}$. 
	
	\begin{figure}[!htp]
		\centering
		\vspace{-7.5pt}
		\begin{subfigure}[t]{.49\textwidth}
			\centering
			\includegraphics[width=\linewidth]{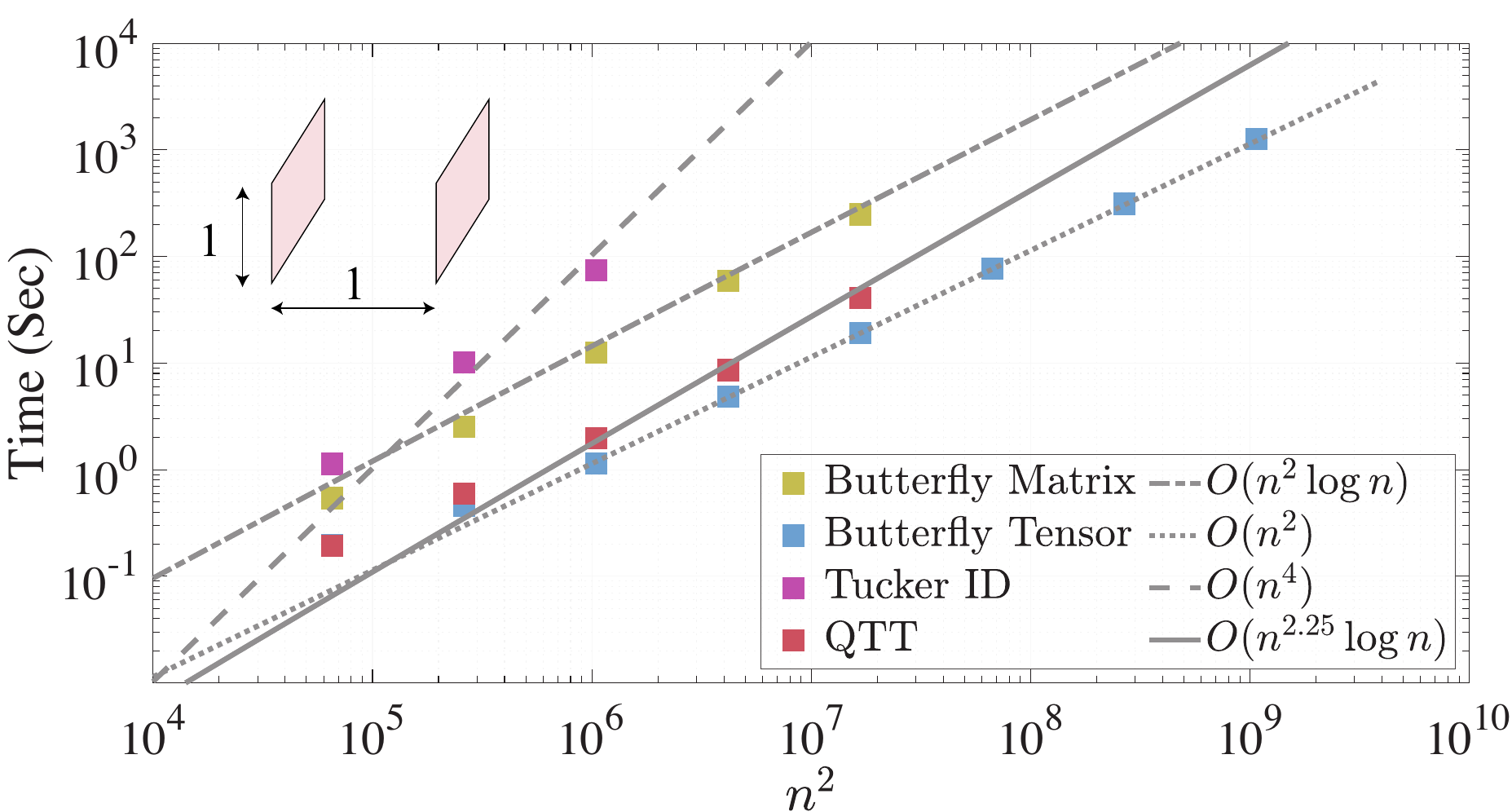}
			%\caption{Velocity model $v({\bf r})$.\label{fig:visco_acoustic_modelV}}
		\end{subfigure}
		%	\vspace{-7.5pt}
		\begin{subfigure}[t]{.49\textwidth}
			\centering
			\includegraphics[width=\linewidth]{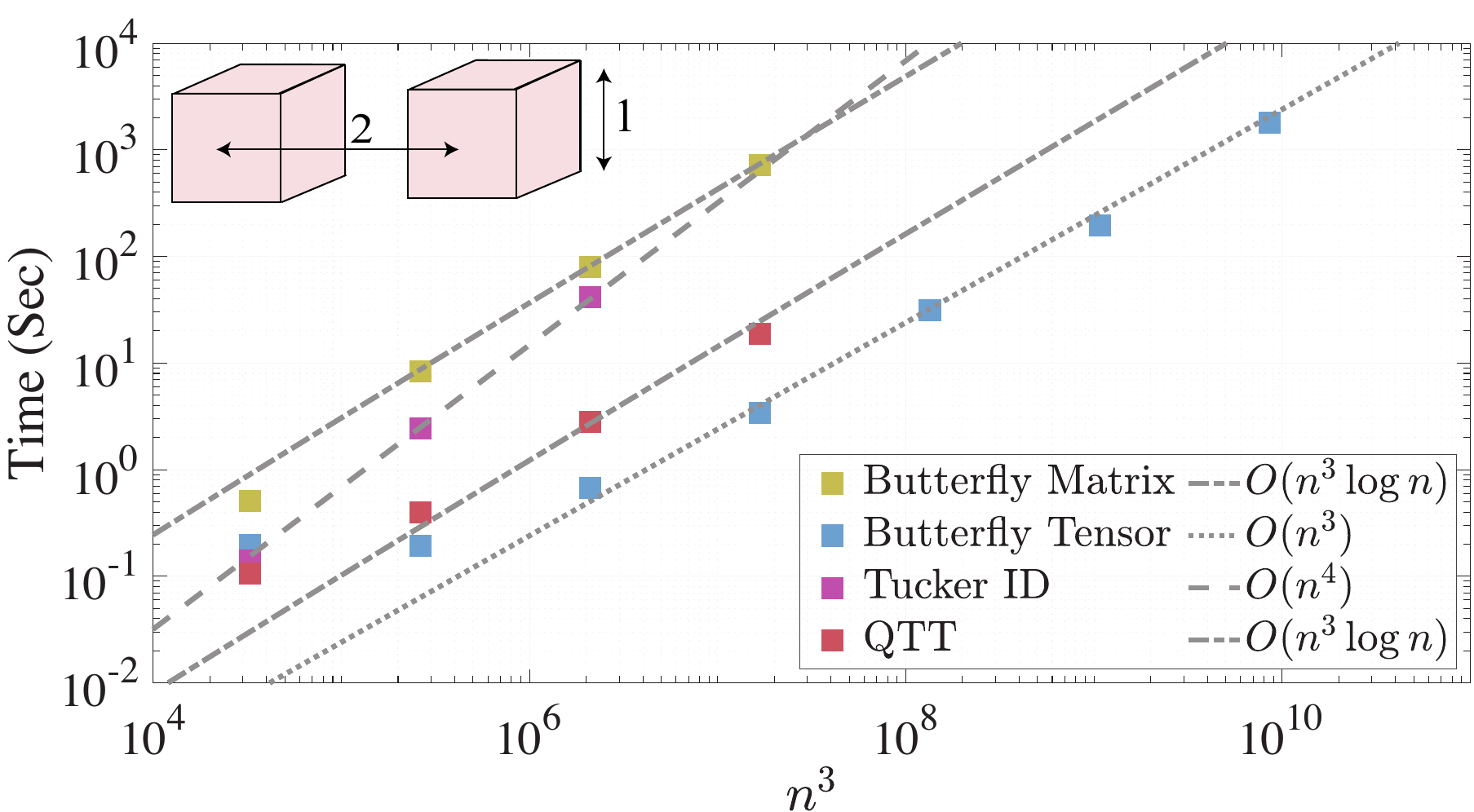}
			%	\caption{Pressure field $p({\bf r})$.\label{fig:visco_acoustic_modelP}}
		\end{subfigure}	
		\begin{subfigure}[t]{.49\textwidth}
			\centering
			\includegraphics[width=\linewidth]{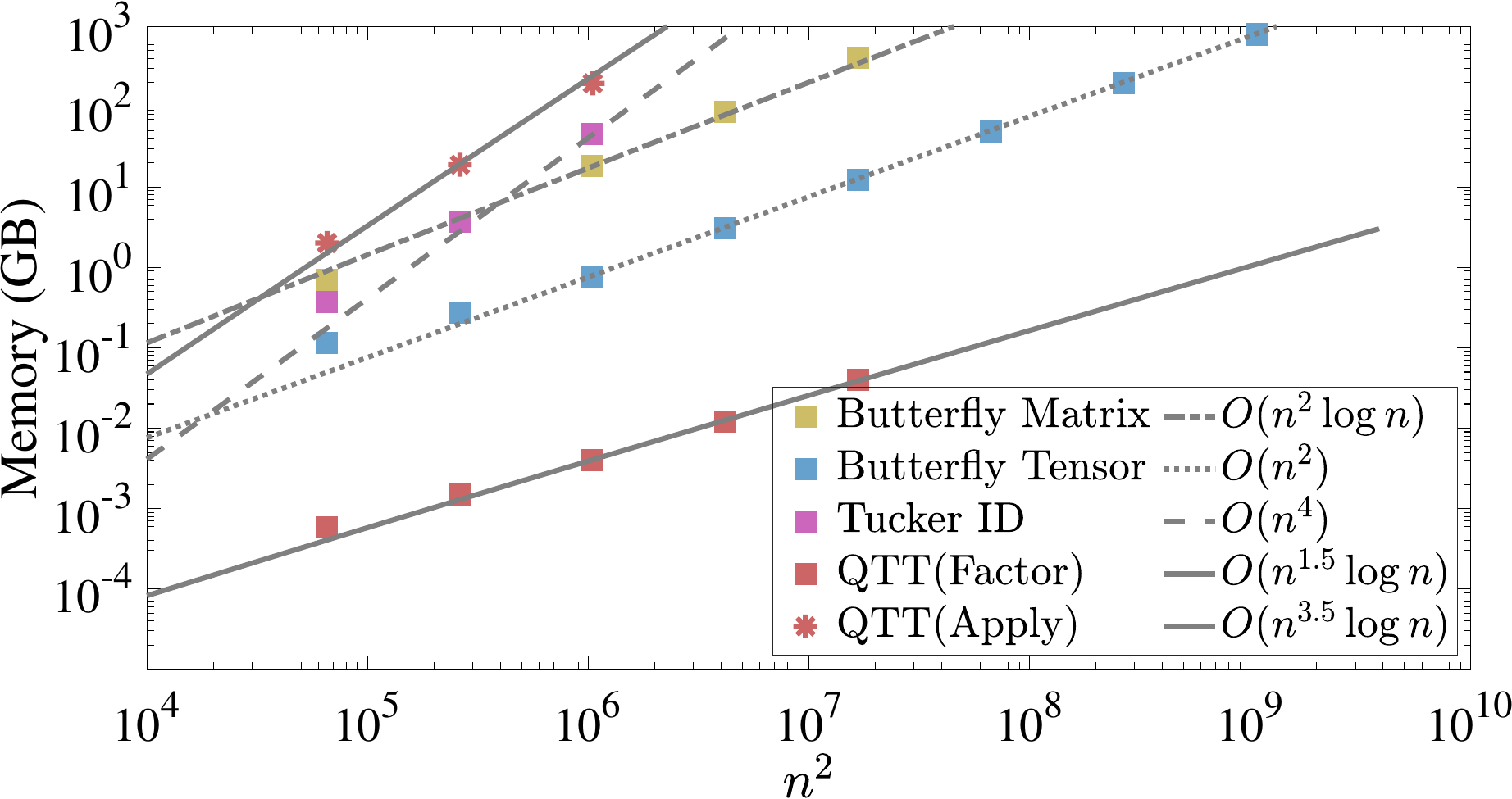}
			%\caption{Velocity model $v({\bf r})$.\label{fig:visco_acoustic_modelV}}
		\end{subfigure}
		%	\vspace{-7.5pt}
		\begin{subfigure}[t]{.49\textwidth}
			\centering
			\includegraphics[width=\linewidth]{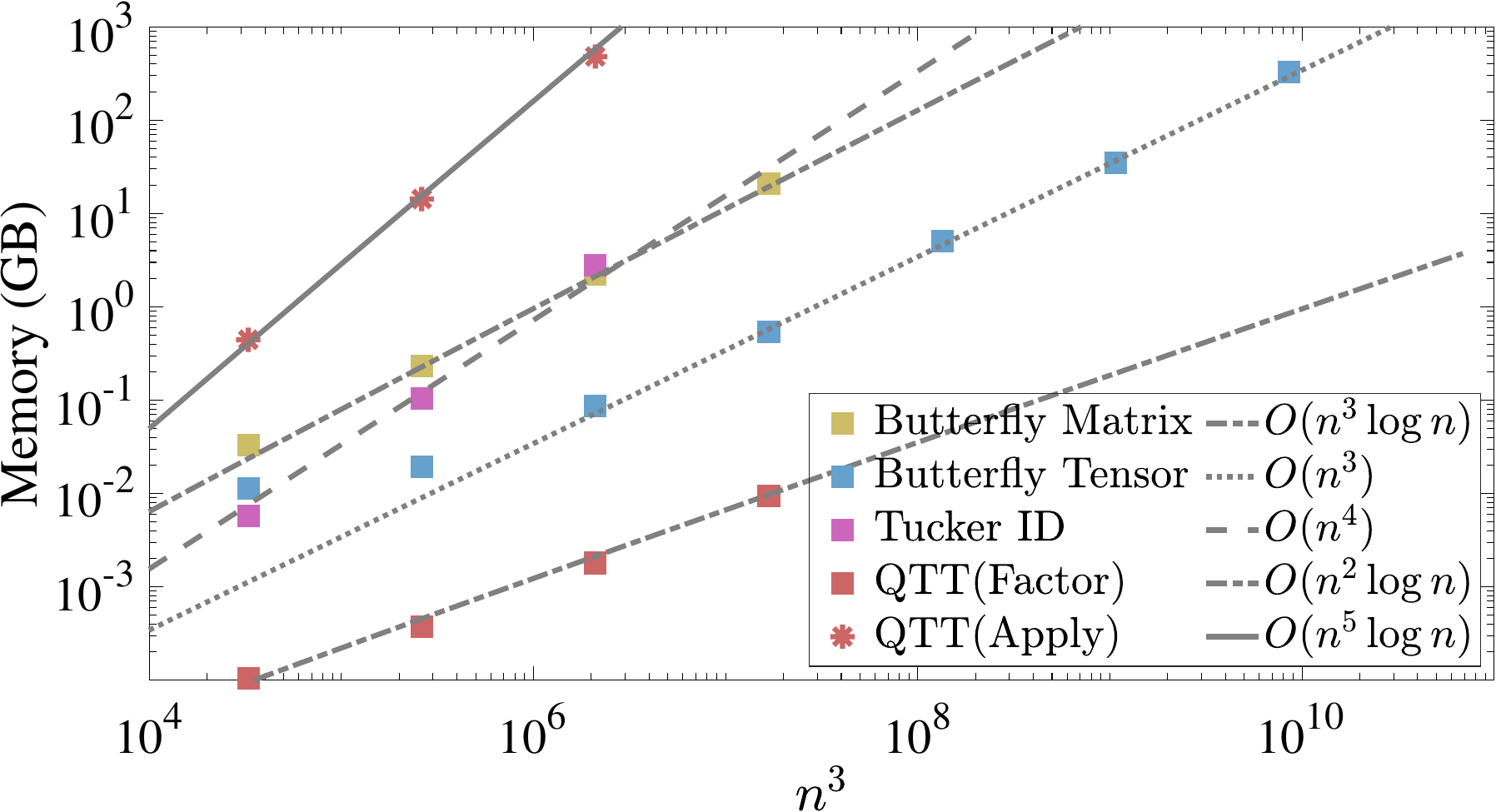}
			%	\caption{Pressure field $p({\bf r})$.\label{fig:visco_acoustic_modelP}}
		\end{subfigure}	
		\begin{subfigure}[t]{.49\textwidth}
			\centering
			\includegraphics[width=\linewidth]{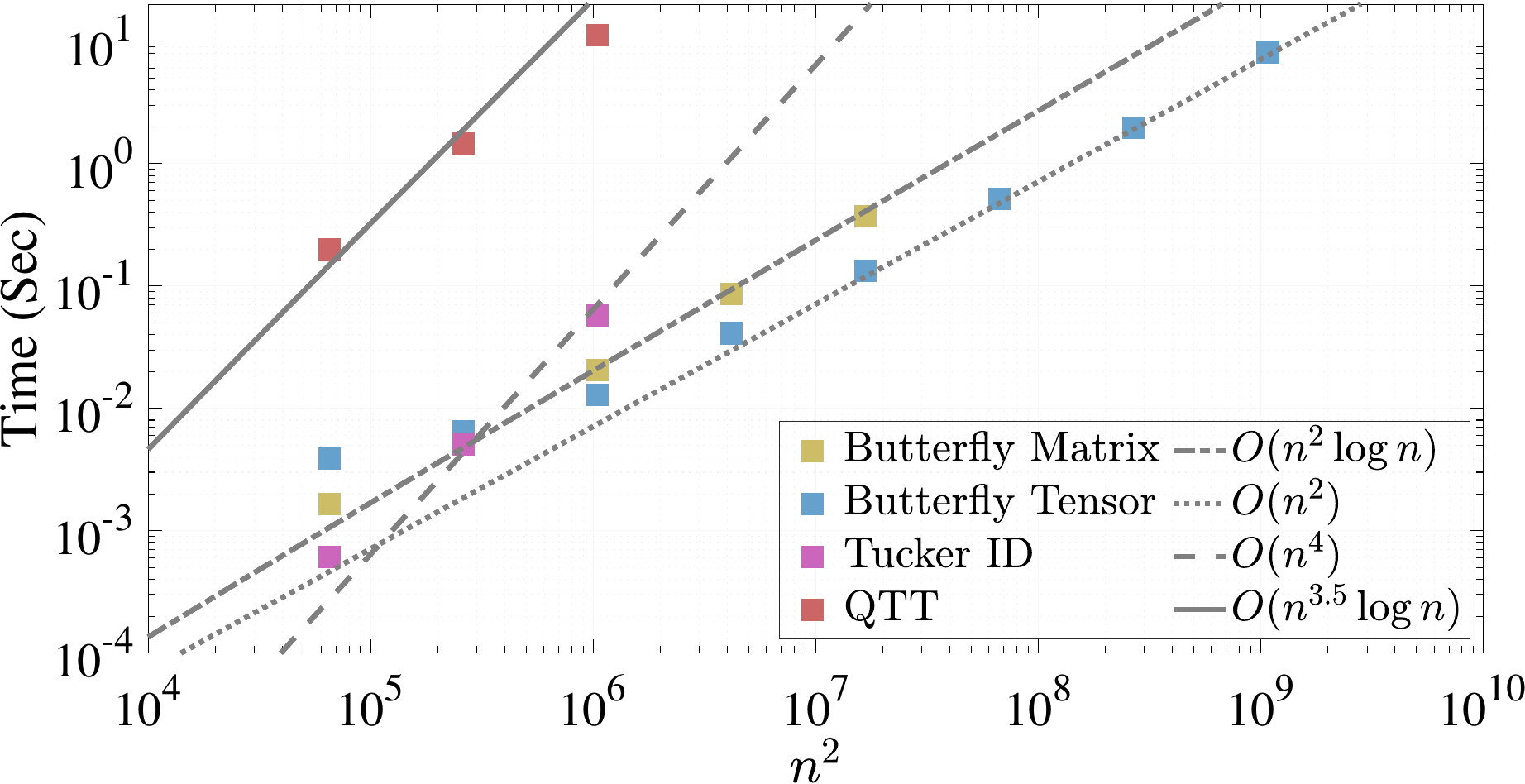}
			%\caption{Velocity model $v({\bf r})$.\label{fig:visco_acoustic_modelV}}
		\end{subfigure}
		%	\vspace{-7.5pt}
		\begin{subfigure}[t]{.49\textwidth}
			\centering
			\includegraphics[width=\linewidth]{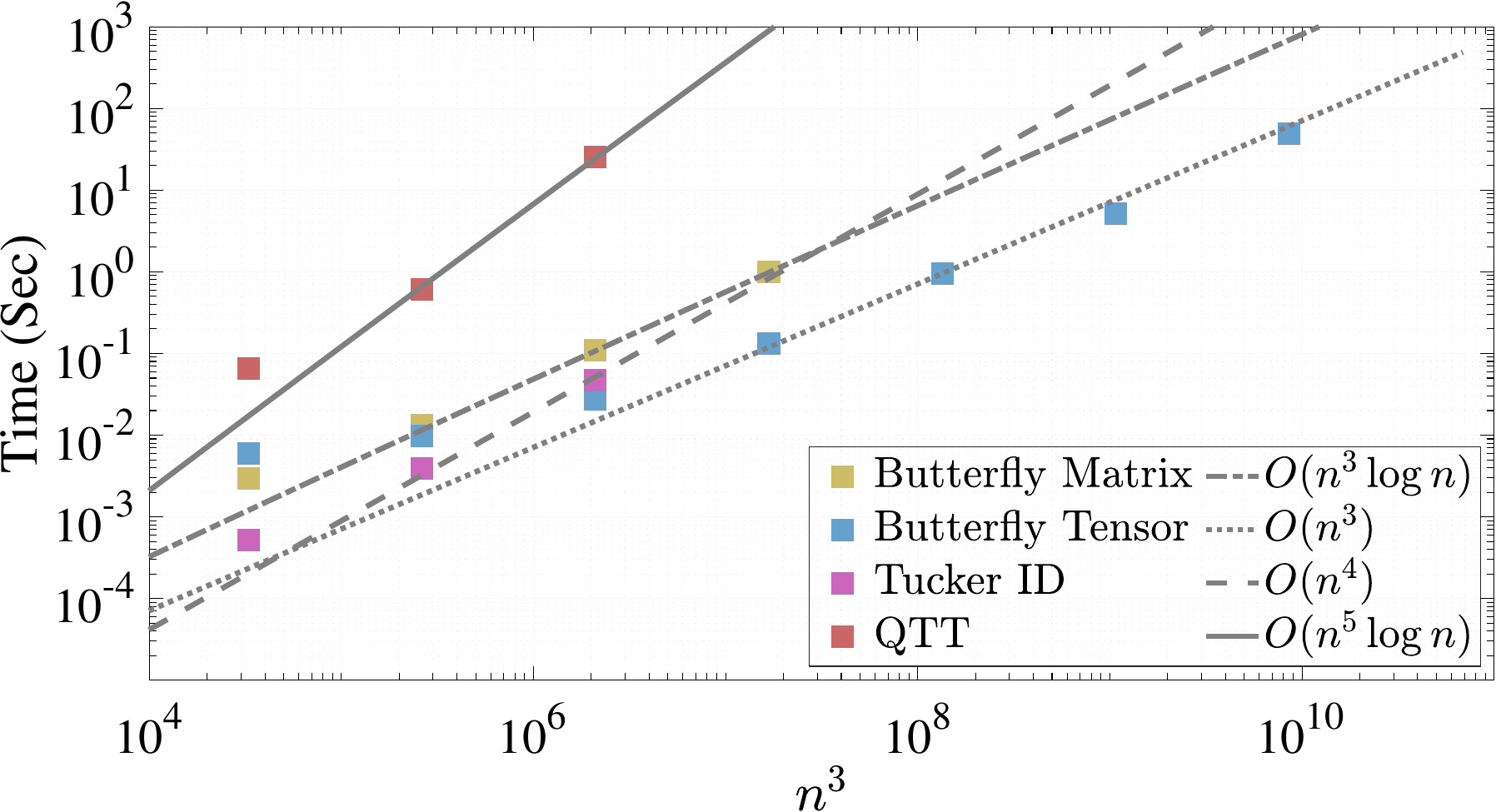}
			%	\caption{Pressure field $p({\bf r})$.\label{fig:visco_acoustic_modelP}}
		\end{subfigure}

		\vspace{-5pt}
		\caption{Helmholtz equation: Computational complexity comparison among butterfly matrix, butterfly tensor, Tucker-ID and QTT for compressing (left) a 4-way Green's function tensor for interactions between two parallel 2D plates and (right) a 6-way Green's function tensor for interactions between two 3D cubes. The geometries are discretized with 4 points per wavelength. (Top): Factor time. (Middle): Factor and apply memory. (Bottom): Apply time. The largest data points correspond to $8192$ wavelengths per direction for the 2D tests (left) and $512$ wavelengths per direction for the 3D tests (right).}
		\label{fig:ccgreen}
		%	\vspace{-20pt}
	\end{figure}	
	
	\Cref{fig:ccgreen} (left) shows the factorization time, application time and memory usage of each algorithm using a compression tolerance $\epsilon=10^{-6}$ for the parallel plate case. For QTT, we show the memory of the factorization (labeled as ``QTT(Factor)") and application (labeled as ``QTT(Apply)") separately. Note that although QTT factorization requires sub-linear memory usage, QTT contraction becomes super-linear due to the full QTT rank of the input tensor. Overall, we achieve the expected complexities listed in \cref{tab:cc} for the butterfly and Tucker-ID algorithms. For QTT, however, instead of an $O(n)$ rank scaling, we observe an $O(n^{3/4})$ rank scaling, leading to slightly better complexities compared with \cref{tab:cc}. We leave this as a future investigation. That said, the tensor butterfly algorithm achieves the linear CPU and memory complexities for both factorization and application with a much smaller prefactor compared to all the other algorithms. Remarkably, the tensor butterfly algorithm achieves a 30x memory reduction and 15x speedup, capable of handling 64x larger-sized tensors compared with the matrix butterfly algorithm.      
	
	\Cref{fig:ccgreen} (right) shows the factorization time, application time and memory usage of each algorithm using a compression tolerance $\epsilon=10^{-2}$ for the cube case. Overall, we achieve the expected complexities listed in \cref{tab:cc} for all four algorithms. The tensor butterfly algorithm achieves the linear CPU and memory complexities for both factorization and application with a much smaller prefactor compared to all the other algorithms. Remarkably, the tensor butterfly algorithm achieves a $30\times$ memory reduction and 200x speedup, capable of handling $512\times$ larger-sized tensors compared with the matrix butterfly algorithm. The largest data point $n=2048$ corresponds to 512 wavelengths per physical dimension. The results in \cref{fig:ccgreen} suggest the superiority of the tensor butterfly algorithm in solving high-frequency wave equations in 3D volumes and on 3D surfaces.      
	
	Next, we demonstrate the effect of changing compression tolerance $\epsilon$ for both test cases in \cref{tab:green}. Here the error is measured by
	\begin{eqnarray}
	\mathrm{error}=\frac{||\ts{K}\times_{d+1,d+2,\ldots,2d}\ts{F}_{e}-\ts{K}_{BF}\times_{d+1,d+2,\ldots,2d}\ts{F}_{e}||_F}{||\ts{K}\times_{d+1,d+2,\ldots,2d}\ts{F}_{e}||_F}\label{eq:error_metric}
	\end{eqnarray} 
	where $\ts{K}_{BF}$ is the tensor butterfly representation of $\ts{K}$, $\ts{F}_{e}(\midx{j})=1$ for a small set of random entries $\midx{j}$ and $0$ elsewhere. This way, $\ts{K}$ does not need to be fully formed to compute the error. \cref{tab:green} shows the minimum rank ($r_{\min}$) and maximum rank ($r$), error, factorization time, application time and memory usage of varying $\epsilon$, for $n=16384,d=2$ and $n=512,d=3$. \ylrev{We remark that the observed ranks clearly demonstrate $\Delta_\epsilon=O(\log \epsilon^{-1})$ in \cref{eqn:r_t}.} Overall, the errors are close to the prescribed tolerances and the costs increase for smaller $\epsilon$, as expected. We also note that keeping $r$ as low as possible is critical in maintaining small prefactors of the tensor butterfly algorithm, particularly for higher dimensions.

	\begin{table*}[!tp]
		\centering	
		\begin{tabular}{|c|c|c|c|c|c|c|c|}
			\hline
			$n^d$ &$\epsilon$ & $r_\mathrm{min}$ & $r$ & error & $T_f$(sec) & $T_a$ (sec) &  $Mem$ (MB)\\		
			\hline
			\hline
			$16384^2$	&1E-02 &5 &8	&1.49E-02 	 &6.83E+01 &  1.16E+00 &2.40E+04\\		
			\hline
			$16384^2$	&1E-03 &6&10	&2.19E-03 	 &1.17E+02 &  1.89E+00 &4.69E+04\\
			\hline		
			$16384^2$	&1E-04 &7 & 11	&1.84E-04 	 &1.57E+02 &  2.80E+00 &7.49E+04 \\
			\hline
			$16384^2$	&1E-05 &8 & 12	 &3.46E-05 	 &2.29E+02	 &  4.03E+00 &1.21E+05 \\
			\hline
			$16384^2$	&1E-06  &9 &13	 &9.26E-06   &3.18E+02	 &  5.92E+00 &1.96E+05 \\
			\hline
			\hline
			$512^3$	&1E-02 &2 &5	&2.01E-02 	 &1.18E+02 &  1.42E+00 &1.19E+04\\		
			\hline
			$512^3$	&1E-03 &2 &6	&1.18E-03 	 &3.46E+02 &  4.08E+00 &4.87E+04\\		
			\hline
			$512^3$	&1E-04 &2 &7	&8.39E-05 	 &6.26E+02 &  9.85E+00 &1.49E+05\\		
			\hline
			$512^3$	&1E-05 &3 &8	&9.21E-06 	 &1.25E+03 &  2.40E+01&4.07E+05\\		
			\hline		
		\end{tabular}
		\caption{The technical data for a 4-way Green's function tensor of $n=16384$ and a 6-way Green's function tensor of $n=512$ for the Helmholtz equation using the proposed tensor butterfly algorithm of varying compression tolerance $\epsilon$. The table shows the maximum rank $r$ and minimum rank $r_{\mathrm{min}}$ across all ID operations, relative error in (\ref{eq:error_metric}), factor time $T_f$, apply time $T_a$, and memory usage $Mem$.}\label{tab:green}
	\end{table*}

	\subsection{Radon transforms}\label{sec:example_radon}
	In this subsection, we consider 2D and 3D discretized Radon transforms similar to those presented in \cite{Candes_butterfly_2009}. Specifically, the tensor entry is  
	\begin{eqnarray}
	\ts{K}({\midx{i}},{\midx{j}})=\mathrm{exp}(2\pi\mathrm{i}\phi(x^{\midx{i}},y^{\midx{j}}))
	\end{eqnarray}
	with $x^{\midx{i}} = (\frac{i_1}{n},\frac{i_2}{n},\ldots, \frac{i_d}{n})$ and $y^{\midx{j}} = (j_1-\frac{n}{2},j_2-\frac{n}{2},\ldots, j_d-\frac{n}{2})$. For $d=2$, we consider 
	\begin{eqnarray}
	\phi(x,y)=x\cdot y + \sqrt{c_1^2y_1^2 + c_2^2y_2^2} \label{eq:randon2d},\\
	c_1=(2+\sin(2\pi x_{1})\sin(2\pi x_{2}))/16,\nonumber\\
	c_2=(2+\cos(2\pi x_{1})\cos(2\pi x_{2}))/16.\nonumber
	\end{eqnarray}
	For $d=3$, we consider 
	\begin{eqnarray}
	&\phi(x,y)=x\cdot y + c|y|, \label{eq:randon3d}\\
	&c=(3+\sin(2\pi x_{1})\sin(2\pi x_{2})\sin(2\pi x_{3}))/100.\nonumber
	\end{eqnarray}
	We first perform compression using the matrix butterfly, tensor butterfly, and QTT algorithms, and then perform application/contraction using a random input tensor $\ts{F}$. 
	
	\begin{figure}[!tp]
		\centering
		\vspace{-7.5pt}
		\begin{subfigure}[t]{.49\textwidth}
			\centering
			\includegraphics[width=\linewidth]{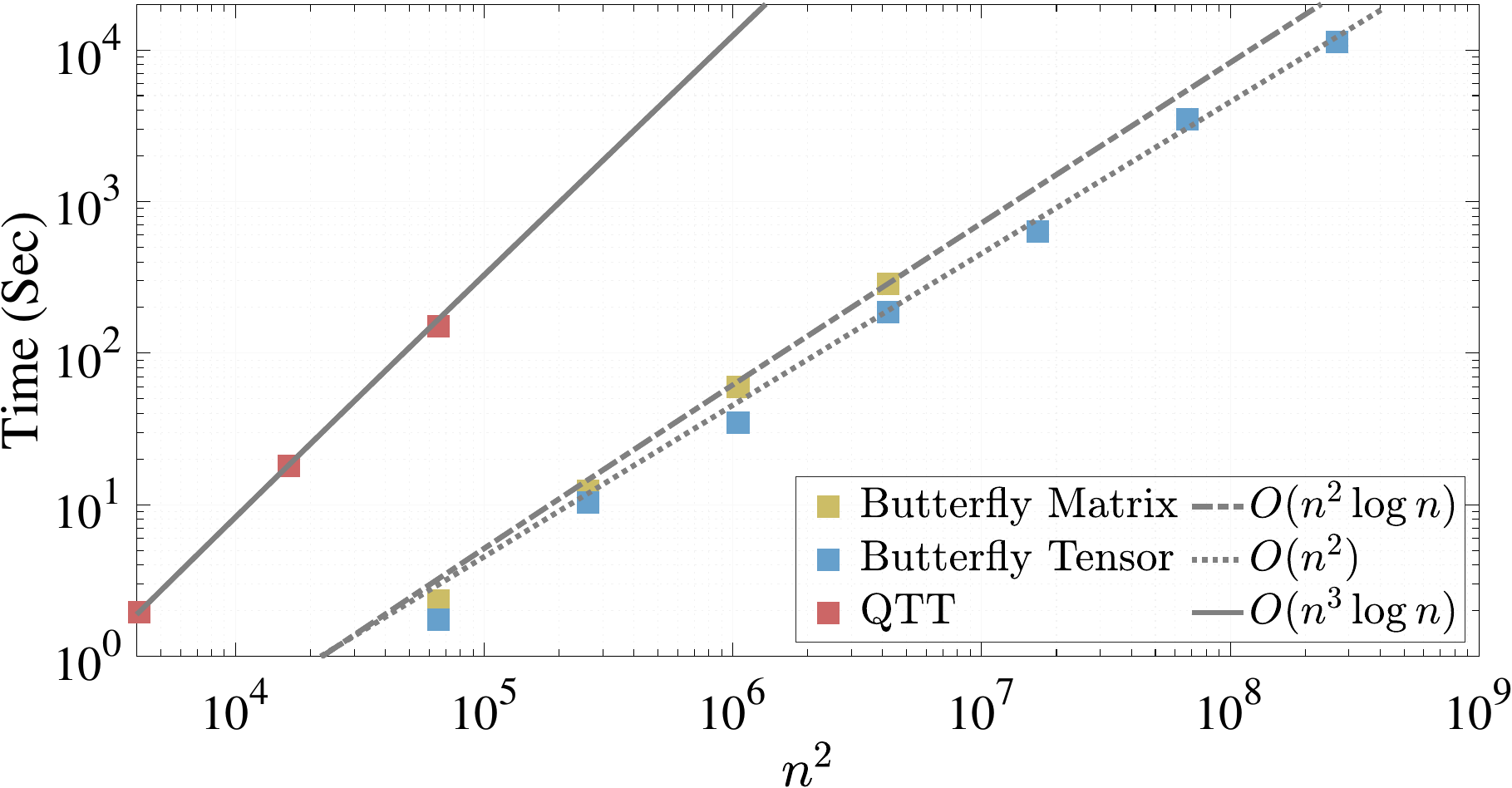}
			%\caption{Velocity model $v({\bf r})$.\label{fig:visco_acoustic_modelV}}
		\end{subfigure}
		%	\vspace{-7.5pt}
		\begin{subfigure}[t]{.49\textwidth}
			\centering
			\includegraphics[width=\linewidth]{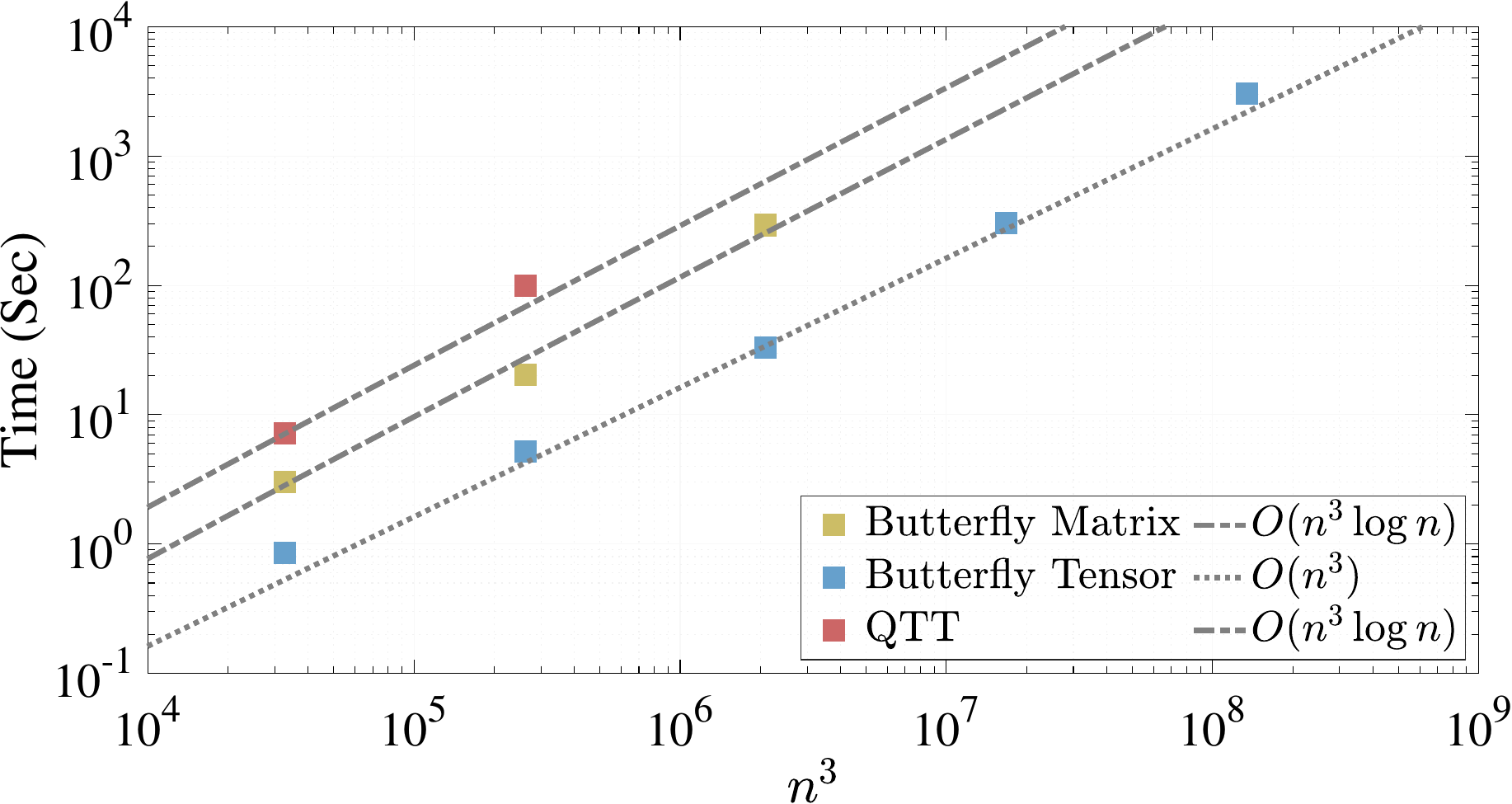}
			%	\caption{Pressure field $p({\bf r})$.\label{fig:visco_acoustic_modelP}}
		\end{subfigure}	
		\begin{subfigure}[t]{.49\textwidth}
			\centering
			\includegraphics[width=\linewidth]{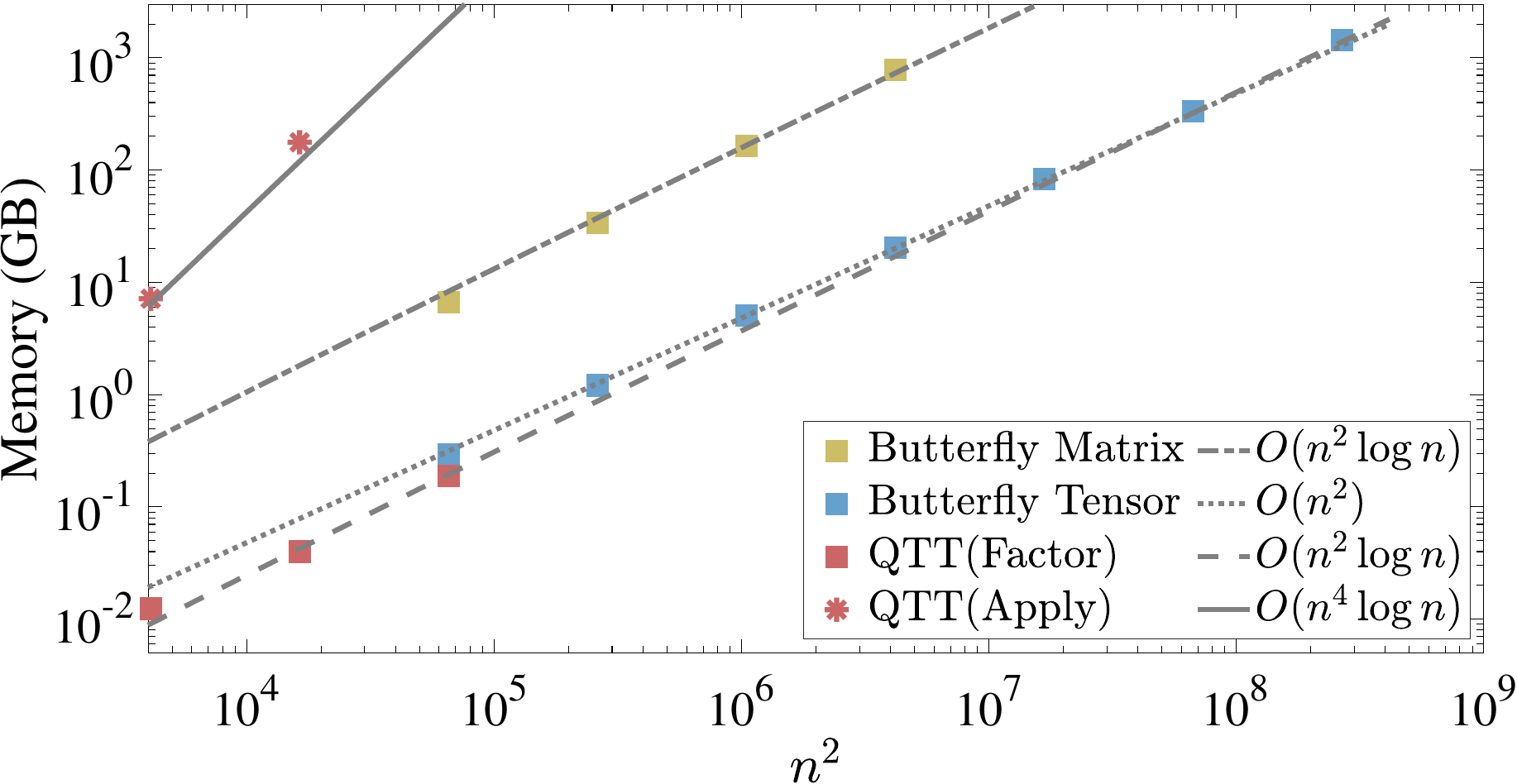}
			%\caption{Velocity model $v({\bf r})$.\label{fig:visco_acoustic_modelV}}
		\end{subfigure}
		%	\vspace{-7.5pt}
		\begin{subfigure}[t]{.49\textwidth}
			\centering
			\includegraphics[width=\linewidth]{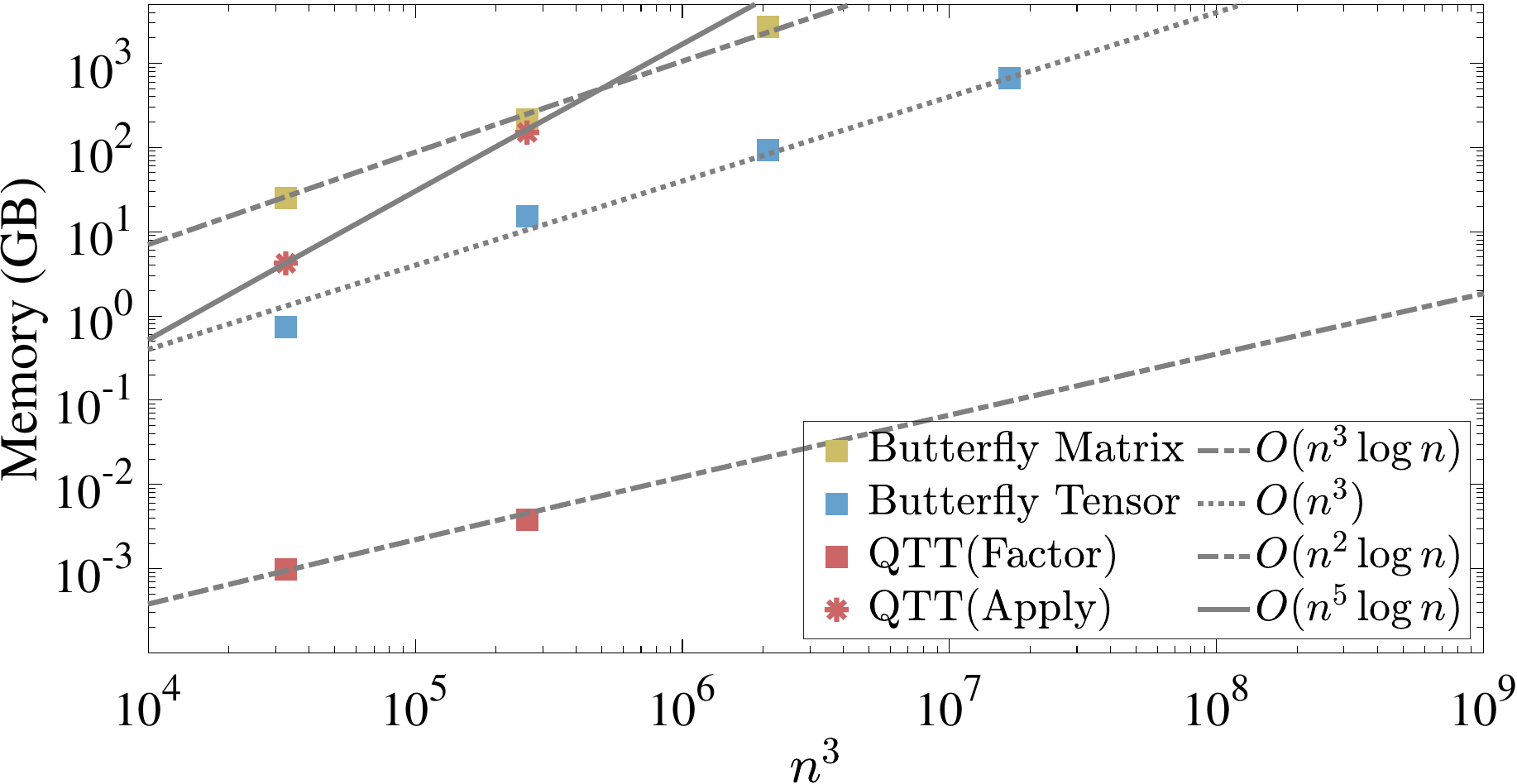}
			%	\caption{Pressure field $p({\bf r})$.\label{fig:visco_acoustic_modelP}}
		\end{subfigure}	
		\begin{subfigure}[t]{.49\textwidth}
			\centering
			\includegraphics[width=\linewidth]{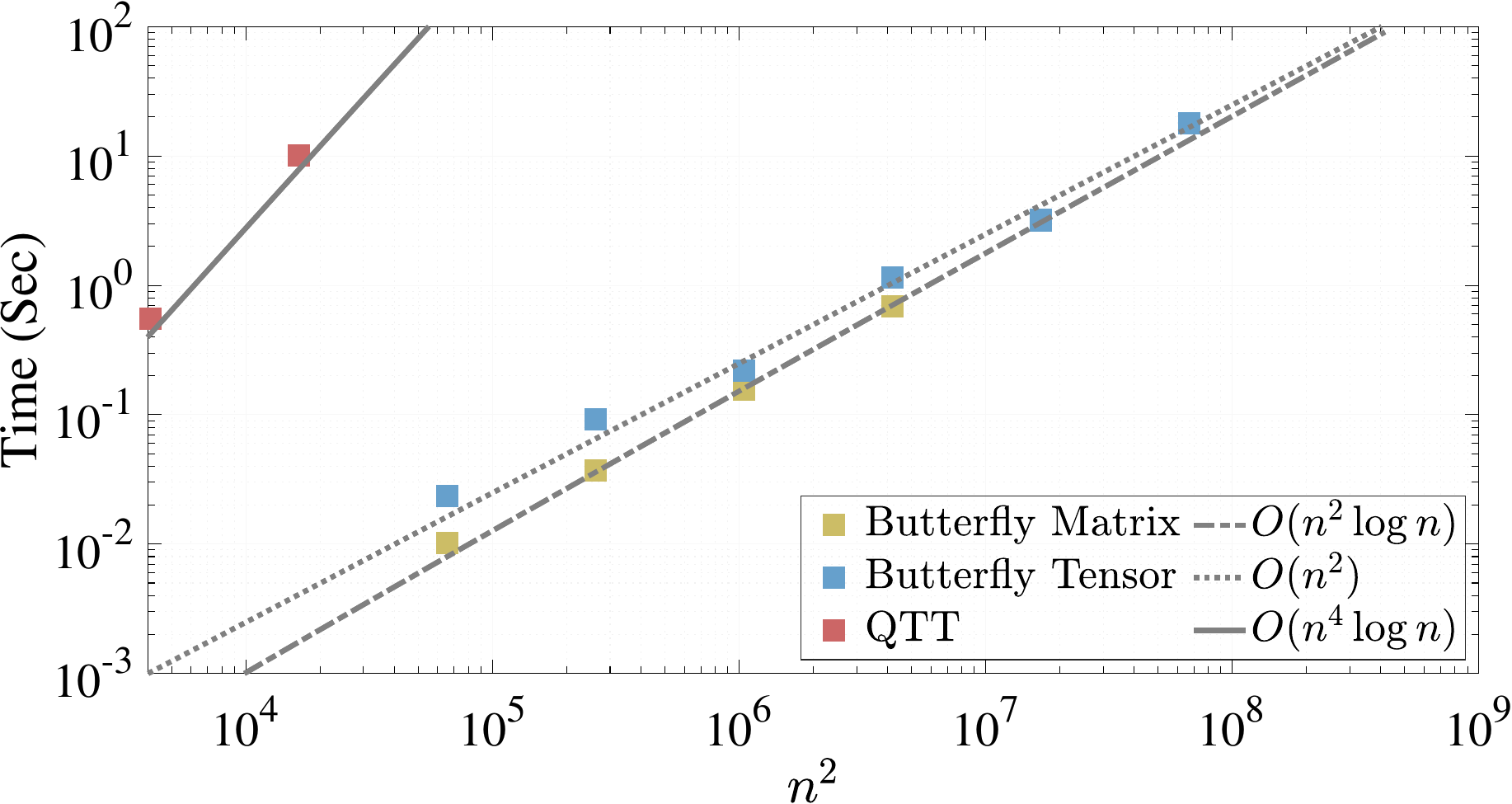}
			%\caption{Velocity model $v({\bf r})$.\label{fig:visco_acoustic_modelV}}
		\end{subfigure}
		%	\vspace{-7.5pt}
		\begin{subfigure}[t]{.49\textwidth}
			\centering
			\includegraphics[width=\linewidth]{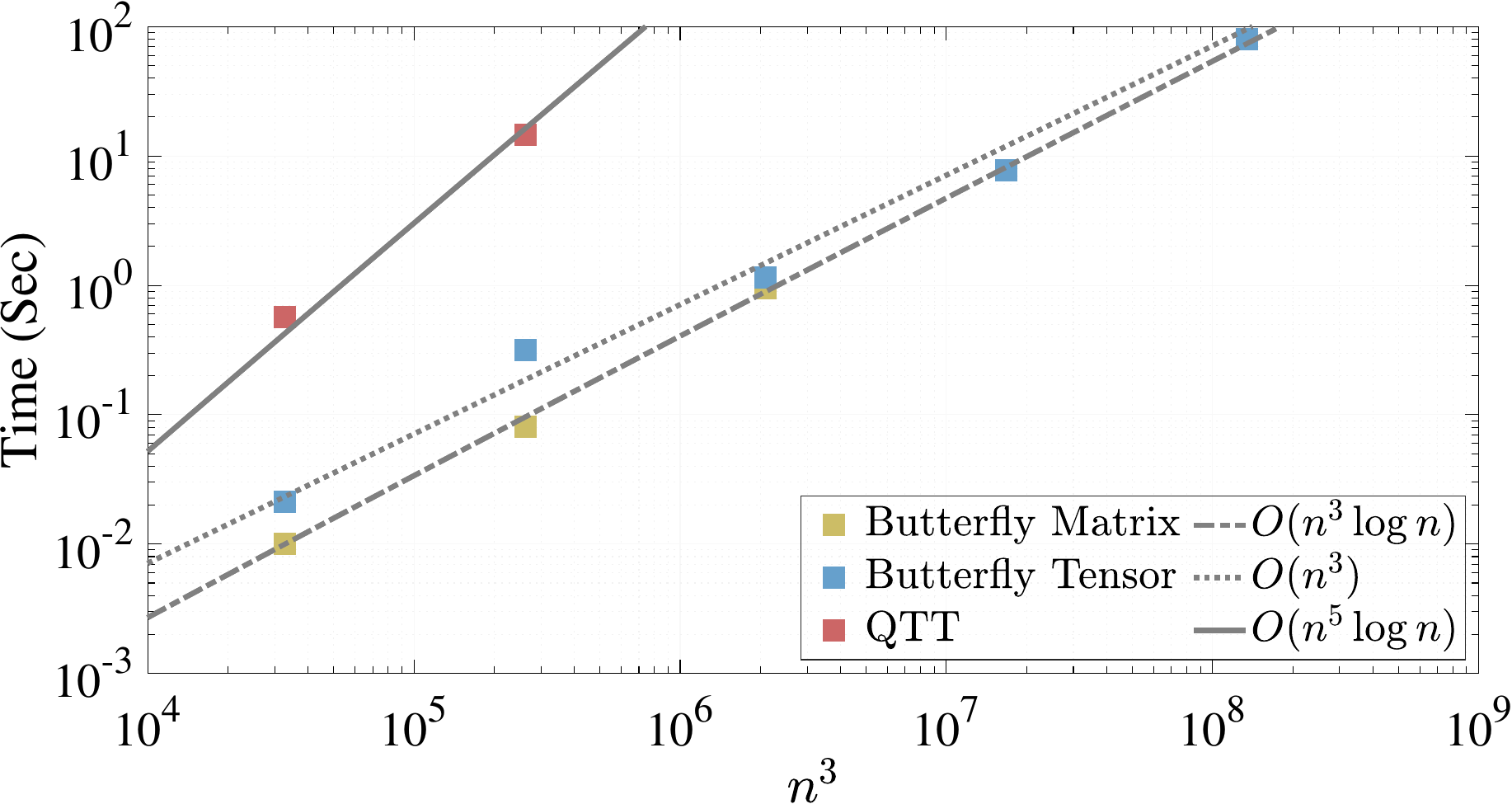}
			%	\caption{Pressure field $p({\bf r})$.\label{fig:visco_acoustic_modelP}}
		\end{subfigure}

		\vspace{-5pt}
		\caption{Radon transforms: Computational complexity comparison among butterfly matrix, butterfly tensor and QTT for compressing (left) a 2D Radon transform tensor and (right) a 3D Radon transform tensor. (Top): Factor time. (Middle): Factor and apply memory. (Bottom): Apply time. }
		\label{fig:ccradon}
		%	\vspace{-20pt}
	\end{figure}

	\Cref{fig:ccradon} shows the factorization time, application time and memory usage of each algorithm using a compression tolerance $\epsilon=10^{-3}$ for the 2D transform (left) and 3D transform (right). Overall, we achieve the expected complexities listed in \cref{tab:cc} for all three algorithms. The QTT algorithm can only obtain the first 2 or 3 data points due to its high memory usage and large QTT ranks. In comparison, the tensor butterfly algorithm achieves the linear CPU and memory complexities for both factorization and application with a much smaller prefactor compared to all the other algorithms. Note that the  Radon transform kernels in (\ref{eq:randon2d}) and (\ref{eq:randon3d}) are not translational invariant, but the tensor butterfly algorithm can still attain small ranks. As a result, the tensor butterfly algorithm can handle 64x larger-sized Radon transforms compared with the matrix butterfly algorithm, showing their superiority for solving linear inverse problems in tomography and seismic imaging.
	
	\begin{table*}[!tp]
		\centering	
		\begin{tabular}{|c|c|c|c|c|c|c|c|}
			\hline
			$n^d$ &$\epsilon$ & $r_\mathrm{min}$ & $r$ & error & $T_f$(sec) & $T_a$ (sec) &  $Mem$ (MB)\\		
			\hline
			\hline
			$2048^2$	&1E-02 &4 &18	&2.04E-02 	 &9.32E+01&  7.20E-01 &1.25E+04\\		
			\hline
			$2048^2$	&1E-03 &4&20	&1.51E-03 	 &1.61E+02 &  1.28E+00&2.40E+04\\
			\hline		
			$2048^2$	&1E-04 &4 & 22	&1.49E-04 	 &2.55E+02 &  2.05E+00&4.26E+04 \\
			\hline
			$2048^2$	&1E-05 &4 & 23	 &2.45E-05 	 &3.73E+02	 &  3.12E+00 &6.95E+04 \\
			\hline
			\hline
			$128^3$	&1E-02 & 2&	6&4.31E-02 	 &3.89E+01 &  8.57E-01&1.59E+04\\		
			\hline
			$128^3$	&1E-03 & 2&	8&1.00E-02 	 &1.31E+02 &  3.74E+00 &9.44E+04\\		
			\hline
			$128^3$	&1E-04 & 2&	9&1.68E-03 	 &2.42E+02 &  8.28E+00 &2.38E+05\\		
			\hline		
			$128^3$	&1E-05 & 2&	11&1.48E-04 	 &4.30E+02 &  2.05E+01 &6.06E+05\\		
			\hline		
		\end{tabular}
		\caption{The technical data for a 4-way Radon transform tensor of $n=2048$ in (\ref{eq:randon2d}) and a 6-way Radon transform tensor of $n=128$ in (\ref{eq:randon3d}) using the proposed tensor butterfly algorithm of varying compression tolerance $\epsilon$. The table shows the maximum rank $r$ and minimum rank $r_{\mathrm{min}}$ across all ID operations, relative error in (\ref{eq:error_metric}), factor time $T_f$, apply time $T_a$, and memory usage $Mem$..}\label{tab:Radon}
	\end{table*} 
	Next, we demonstrate the effect of changing compression tolerance $\epsilon$ for both test cases in \cref{tab:Radon} with the error defined by (\ref{eq:error_metric}). \cref{tab:Radon} shows the minimum and maximum ranks, error, factorization time, application time and memory usage of varying $\epsilon$, for $n=2048$ with $d=2$ and $n=128$ with $d=3$, respectively. Overall, the errors are close to the prescribed tolerances and the costs increase for smaller $\epsilon$, as expected. Just like the Green's function example, it is critical to keep $r$ a low constant, particularly for higher dimensions.

	\subsection{High-dimensional discrete Fourier transform}\label{sec:example_DFT}
	Finally, we consider high-dimensional DFTs defined as 
	\begin{eqnarray}
	\ts{K}({\midx{i}},{\midx{j}})=\mathrm{exp}(2\pi\mathrm{i}x^{\midx{i}}\cdot y^{\midx{j}}), 
	\end{eqnarray}
	where we choose $x^{\midx{i}} = (i_1-1,i_2-1,\ldots, i_d-1)$ and $y^{\midx{j}} = (\frac{j_1-1}{n},\frac{j_2-1}{n},\ldots, \frac{j_d-1}{n})$ for uniform DFTs, and we choose $x^{\midx{i}}$ to be random (in the sense that $x^{\midx{i}}_k\in[0,n-1]$ for $k\leq d$ is a random number) and $y^{\midx{j}} = (\frac{j_1-1}{n},\frac{j_2-1}{n},\ldots, \frac{j_d-1}{n})$ for type-2 non-uniform DFTs. For high-dimensional DFTs with $d=3,4,5,6$, we perform compression using the tensor butterfly algorithms (with the bit-reversal ordering for each dimension), and perform application/contraction using a random input tensor $\ts{F}$. In comparison, for $d=3$ we perform FFT via the heFFTe package for the uniform DFT example and NUFFT via the FINUFFT package for the type-2 non-uniform DFT example. 
	
\begin{figure}[!htp]
		\centering
		\vspace{-7.5pt}
		\begin{subfigure}[t]{.49\textwidth}
			\centering
			\includegraphics[width=\linewidth]{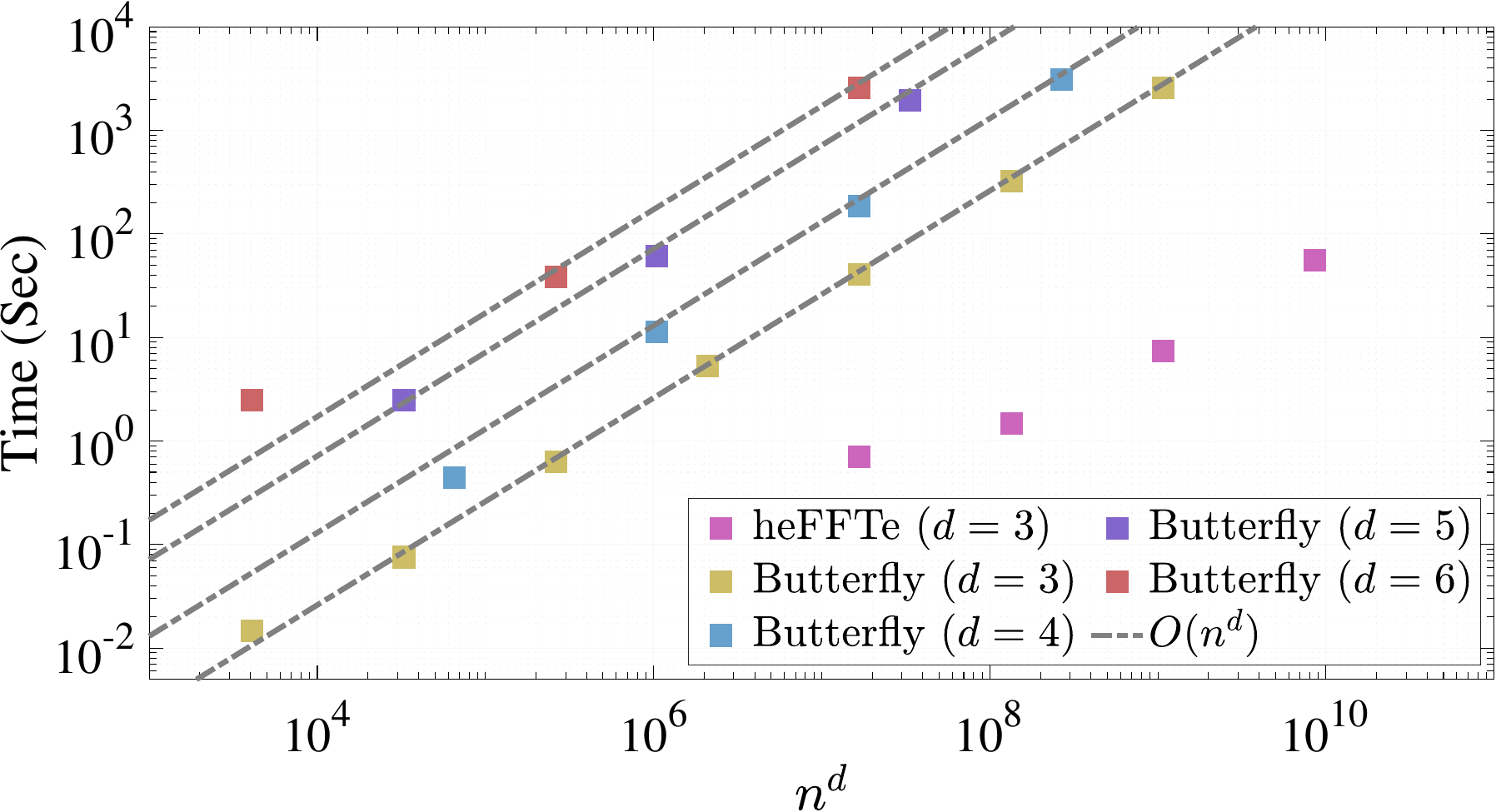}
			%\caption{Velocity model $v({\bf r})$.\label{fig:visco_acoustic_modelV}}
		\end{subfigure}
		%	\vspace{-7.5pt}
		\begin{subfigure}[t]{.49\textwidth}
			\centering
			\includegraphics[width=\linewidth]{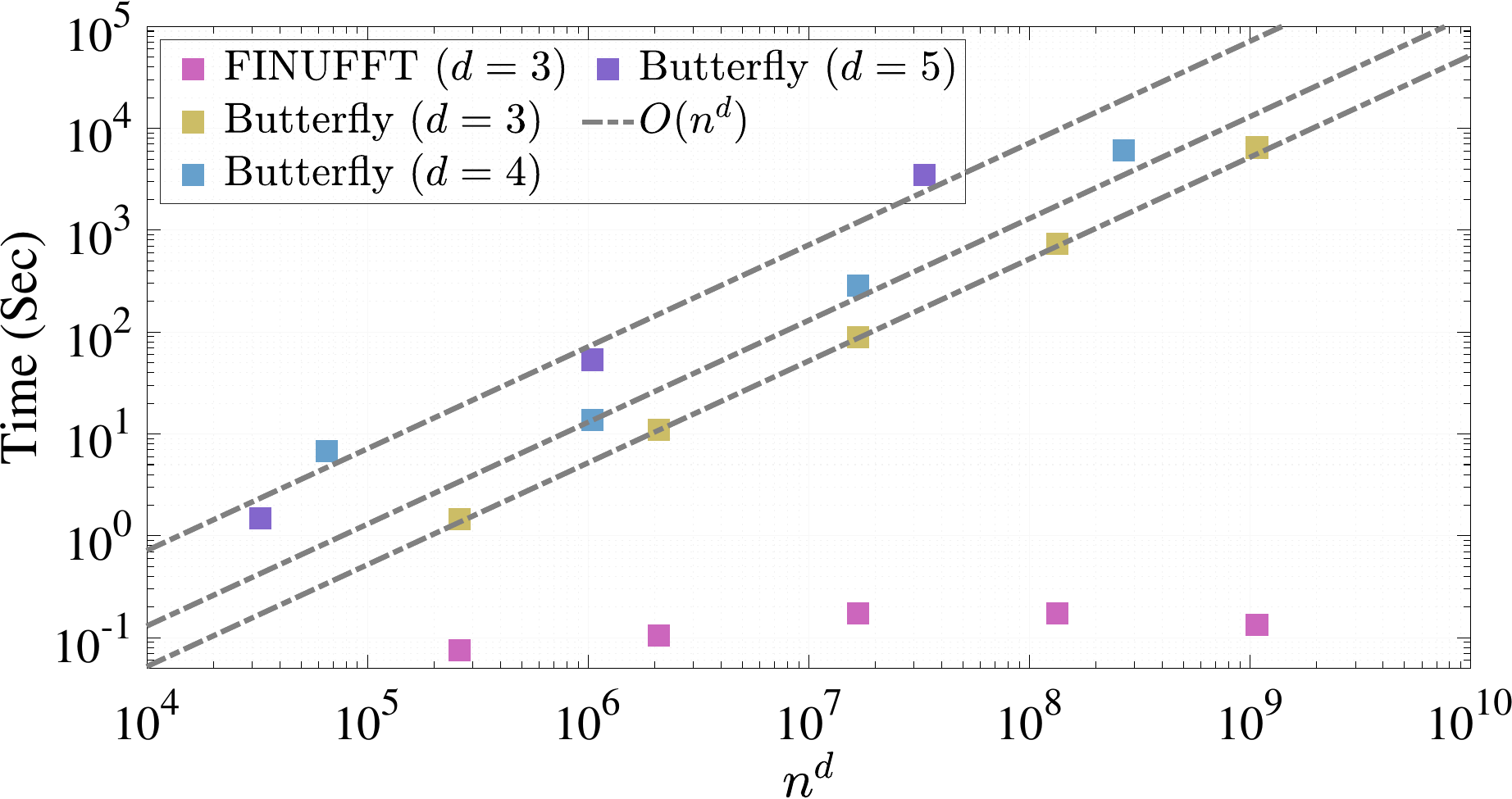}
			%	\caption{Pressure field $p({\bf r})$.\label{fig:visco_acoustic_modelP}}
		\end{subfigure}	
		\begin{subfigure}[t]{.49\textwidth}
			\centering
			\includegraphics[width=\linewidth]{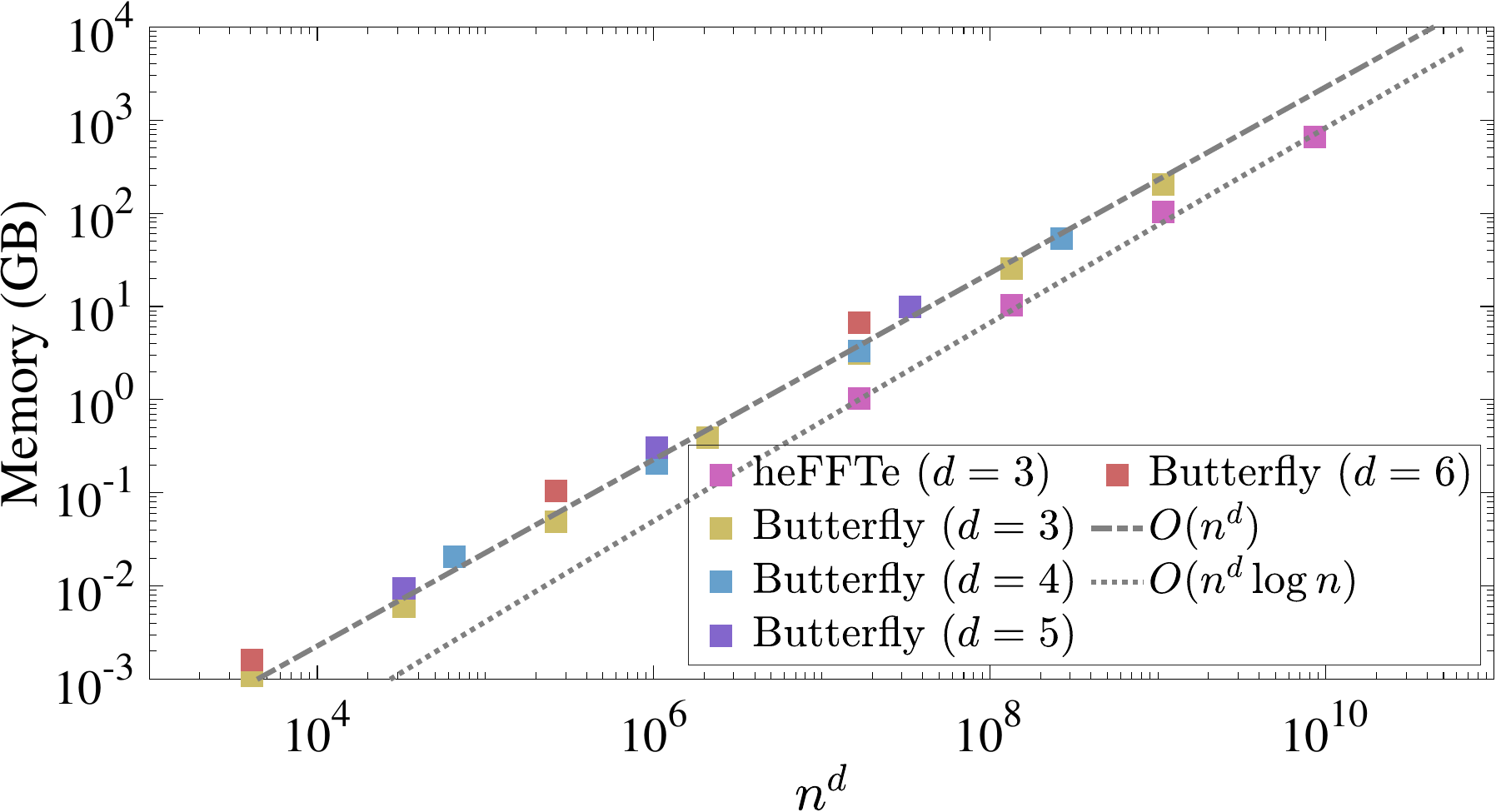}
			%\caption{Velocity model $v({\bf r})$.\label{fig:visco_acoustic_modelV}}
		\end{subfigure}
		%	\vspace{-7.5pt}
		\begin{subfigure}[t]{.49\textwidth}
			\centering
			\includegraphics[width=\linewidth]{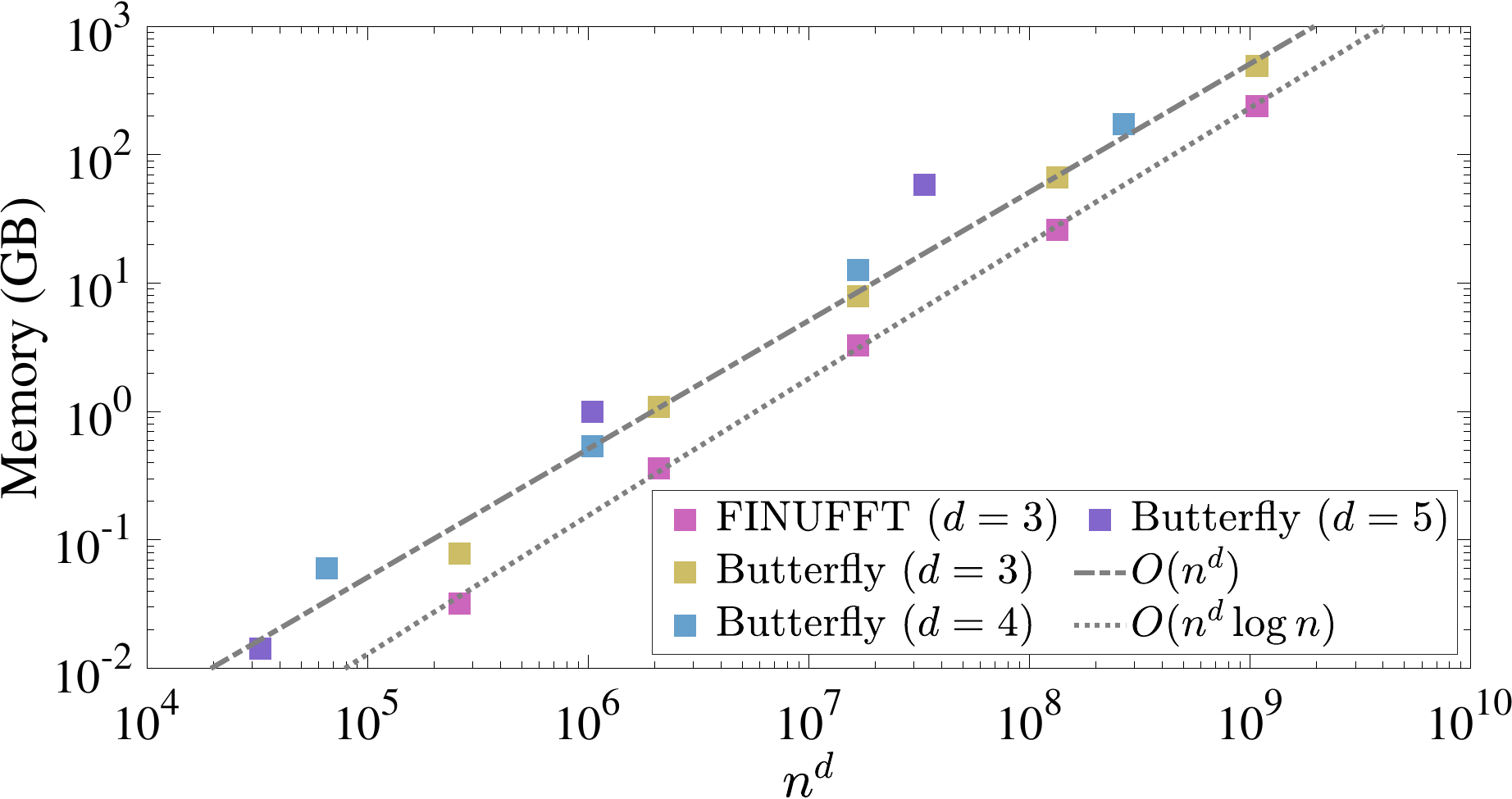}
			%	\caption{Pressure field $p({\bf r})$.\label{fig:visco_acoustic_modelP}}
		\end{subfigure}	
		\begin{subfigure}[t]{.49\textwidth}
			\centering
			\includegraphics[width=\linewidth]{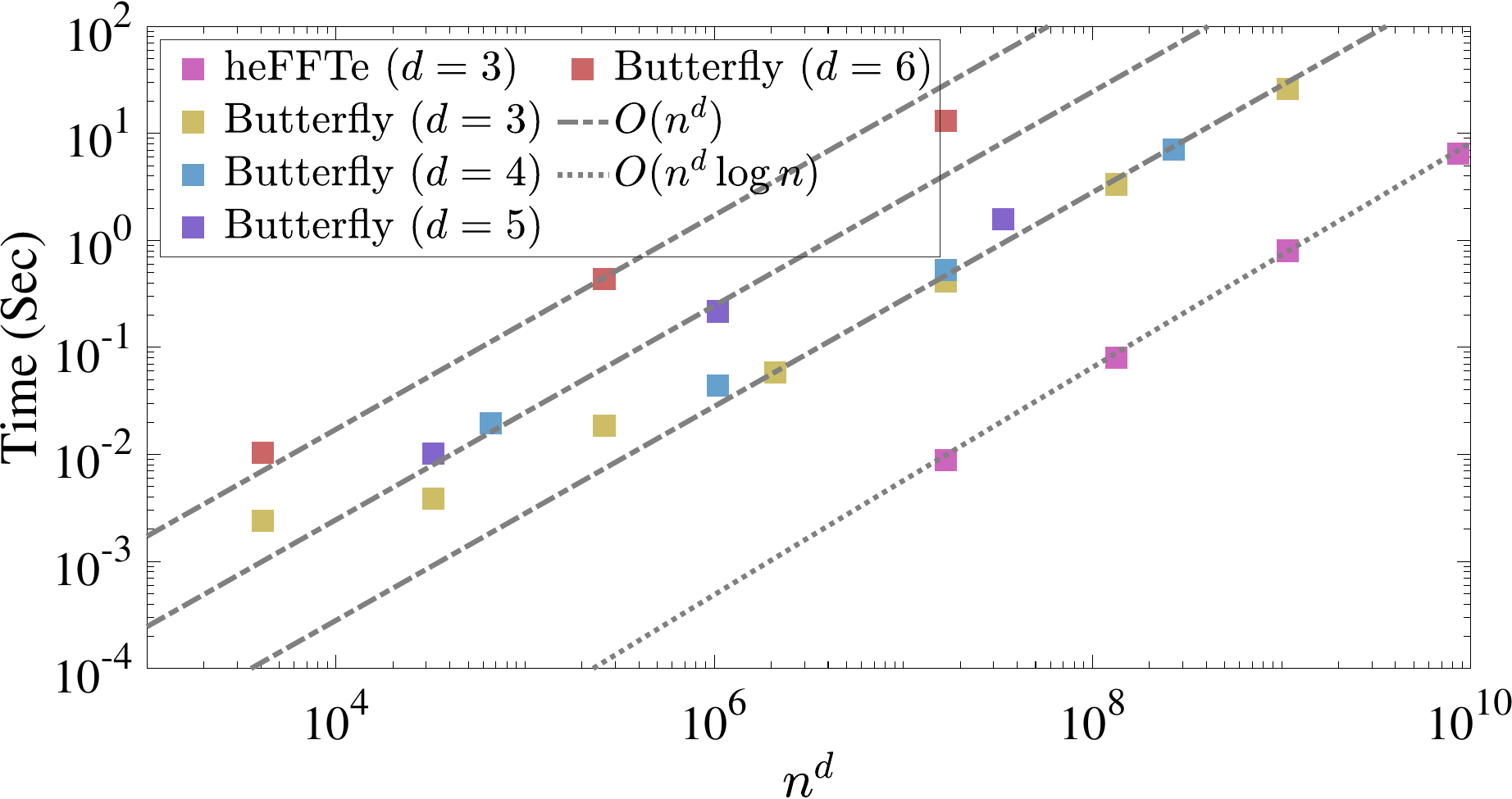}
			%\caption{Velocity model $v({\bf r})$.\label{fig:visco_acoustic_modelV}}
		\end{subfigure}
		%	\vspace{-7.5pt}
		\begin{subfigure}[t]{.49\textwidth}
			\centering
			\includegraphics[width=\linewidth]{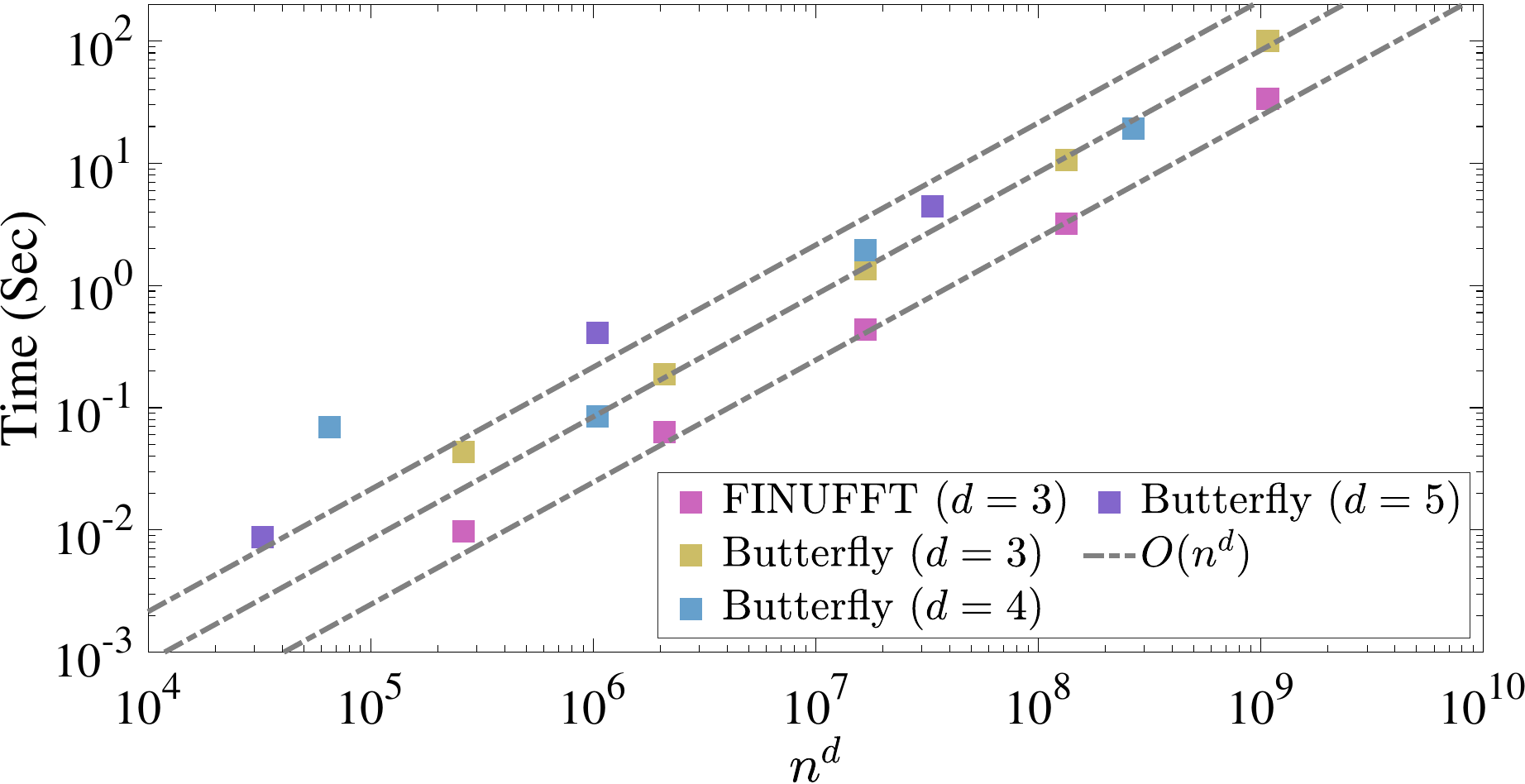}
			%	\caption{Pressure field $p({\bf r})$.\label{fig:visco_acoustic_modelP}}
		\end{subfigure}

		\vspace{-5pt}
		\caption{Fourier transforms: Computational complexity of (left) butterfly tensor and heFFTe for compressing the high-dimensional DFT tensor and (right) butterfly tensor and FINUFFT for compressing the high-dimensional NUFFT tensor. (Top): Factor time of butterfly tensor and plan creation time for heFFTe/FINUFFT. (Middle): Factor memory. (Bottom): Apply time. }
		\label{fig:DFT}
		%	\vspace{-20pt}
	\end{figure}

	\Cref{fig:DFT} shows the factorization time for the butterfly algorithm (or equivalently the plan creation time for heFFTe/FINUFFT), application time and memory usage of each algorithm using a compression tolerance $\epsilon=10^{-3}$ (for butterfly and FINUFFT) for the uniform (left) and nonuniform (right) transforms. Overall, the tensor butterfly algorithm can obtain $O(n^d)$ CPU and memory complexities compared with the $O(n^d\log n)$ complexities of FFT and NUFFT. It is also worth mentioning that QTT can attain logarithmic-complexity uniform DFTs \cite{chen2024direct} when the input tensor $\ts{F}$ is also in the QTT form with low TT ranks. However, for a general input tensor, the complexity of QTT falls back to $O(n^d\log n)$. Although the proposed tensor butterfly algorithm can obtain the best computational complexity among all existing algorithms, we observe that for the $d=3$ case, FFT or NUFFT shows a memory usage similar to the tensor butterfly algorithm but much smaller prefactors for plan creation and application time. That said, the tensor butterfly algorithm provides a unique capability to perform higher dimensional DFTs (i.e., $d\geq4$) with optimal asymptotic complexities.          
	
	\section{Conclusion\label{sec:conclude}} 
	We present a new tensor butterfly algorithm efficiently compressing and applying large-scale and high-dimensional OIOs,  such as Green's functions for wave equations and integral transforms, including Radon transforms and Fourier transforms. The tensor butterfly algorithm leverages an essential tensor CLR property to achieve both improved asymptotic computational complexities and lower leading constants. For the contraction of high-dimensional OIOs with arbitrary input tensors, the tensor butterfly algorithm achieves the optimal linear CPU and memory complexities; this is in huge contrast with both existing matrix algorithms and fast transform algorithms. The former includes the matrix butterfly algorithm, and the latter contains FFT, NUFFT, and other tensor algorithms such as Tucker-type decompositions and QTT. Nevertheless, all these algorithms exhibit higher asymptotic complexities and larger leading constants. As a result, the tensor butterfly algorithm can efficiently model high-frequency 3D Green's function interactions with over $512\times$ larger problem sizes than existing \ylrev{butterfly} algorithms; for the largest sized tensor that can be handled by existing algorithms, the tensor butterfly algorithm requires $200\times$ less CPU time and $30\times$ less memory than existing algorithms. Moreover, it can perform linear-complexity Radon transforms and DFTs with up to $d=6$ dimensions. These OIOs are frequently encountered in the solution of high-frequency wave equations, X-ray and MRI-based inverse problems, seismic imaging and signal processing; therefore, we expect the tensor butterfly algorithm developed here to be both theoretically attractive and practically useful for many applications.   
	
	The limitation of the tensor butterfly algorithm is the requirement for a tensor grid, and hence its extension for unstructured meshes will be a future work. Also, the mid-level subtensors represent a memory bottleneck and need to be compressed with more efficient algorithms.

	\section*{Acknowledgements}

	This research has been supported by the U.S. Department of Energy, Office of Science, Office of Advanced Scientific Computing Research, Mathematical Multifaceted Integrated Capability
	Centers (MMICCs) program, under Contract No. DE-AC02-05CH11231 at Lawrence Berkeley National Laboratory. This work was also supported in part by the NSF Mathematical Sciences Graduate Internship (NSF MSGI) program. This research used resources of the National Energy Research Scientific Computing Center, a DOE Office of Science User Facility supported by the Office of Science of the U.S. Department of Energy under Contract No. DE-AC02-05CH11231 using NERSC award ASCR-ERCAP0024170. Qian's research was partially supported by NSF grants 2152011 and 2309534 and an MSU SPG grant.

	\bibliographystyle{plainurl}
	\bibliography{references}

\end{document}